\documentclass[reqno, a4paper, 12pt]{amsart} %
\usepackage{epsfig}   
\usepackage{ifpdf} 
\usepackage{amssymb,amsmath,amsbsy}         
\usepackage{graphicx}                  

\usepackage{wasysym} 
\usepackage{array}

\newcolumntype{M}[1]{>{\centering\arraybackslash}m{#1}}
\newcolumntype{N}{@{}m{0pt}@{}}

\usepackage[paperwidth=8.5in,paperheight=11.0in,
  left=0.7in,right=0.7in,top=1.0in,bottom=1.0in]{geometry} 

\begin{document}

\makeatletter
\def\rank{\mathop{\operator@font rank}\nolimits}
\def\det{\mathop{\operator@font det}\nolimits}
\makeatother

\newtheorem{thm}{Theorem}[section]
\newtheorem{ex}[thm]{Example}
\newtheorem{lem}[thm]{Lemma}
\newtheorem{rmk}[thm]{Remark}
\newtheorem{defi}[thm]{Definition}


\title[feedback stabilization of double pendulum]{Feedback stabilization of double pendulum: Application to the crane systems with time-varying rope length}

\author{Robert Vrabel} 
\address{Robert Vrabel, Slovak University of Technology in Bratislava, Faculty of Materials Science and Technology, Institute of Applied Informatics, Automation and Mechatronics,  Bottova~25,  917 01 Trnava, Slovakia}
\email{robert.vrabel@stuba.sk}

\date{{\bf\today}}

\begin{abstract}
In the present paper we focus our attention on the design of the feedback-based feed-forward controller asymptotically stabilizing the double-pendulum-type crane system with the time-varying rope length in the desired end position of payload (the origin of the coordinate system). In principle, we will consider two cases, in the first case, the sway angle of payload is uncontrolled and the second case, when the sway angle of payload is controlled by an external force. Mathematical modelling in the framework of Lagrange formalism and numerical simulation in the Matlab environment indicates the substantial reduction of the transportation time to the desired end position. Another principal novelty of this paper lies in deriving and analysis of a complete mathematical model without approximating the nonlinear terms and without neglecting some structural parameters of systems for the reasons described in the Remark~\ref{remark_singularity} and Remark~\ref{rmk_approx}.  
\end{abstract}


\maketitle

\section[Introduction]{Introduction}

For overhead crane control in general, it is required that the trolley should reach the desired location as fast as possible while the payload swing should be kept as little as possible during the transferring process. However, it is extremely challenging to achieve these goals simultaneously owing to the underactuated characteristics of the crane system. Due to this reason, the development of efficient control schemes for overhead cranes has attracted wide attention from the control community. 

For overall context of crane system controls, the paper \cite{Ramli} provides a comprehensive review on modeling and the control strategies  for single-pendulum and double-pendulum-type (D-P-T) crane systems during the years from 2000 to 2016. Also the monograph \cite{Qian_book} introduces anti-sway control approaches for double-pendulum overhead cranes, including control methods, theoretical analyses and simulation results (in Matlab environment).  Passivity-based, sliding-mode-based and Fuzzy-logic-based control methods are discussed here in depth. We refer to these works as an useful preliminary introduction to the methods of crane systems control, where the constant length of rope(s) connecting the trolley with payload is considered. More detailed review on the existing literature regarding the control D-P-T crane systems is provided in the Section~\ref{section:bacground}. As far as we know, the present  paper is the first study dealing with the design of control strategy for stabilization of the D-P-T overhead crane systems at the desired end position of payload based on the complete model of system, where the time-varying length of the hoisting rope is considered and which leads to serious technical difficulties (singularities) in the analysis of their mathematical models if some of the parameters of the system, specifically, the mass moment of inertia, are neglected. In a more general framework, a double pendulum model is used in control theory to measure the effectiveness of stabilizing algorithms. Many real-life physical structures can be approximated with a double pendulum to gain more insight about the system behavior (\cite{comsol}). 
\section{Double-pendulum-type crane system dynamics}

\begin{figure}
\ifpdf\includegraphics[trim={5cm 18cm 5cm 0},clip]{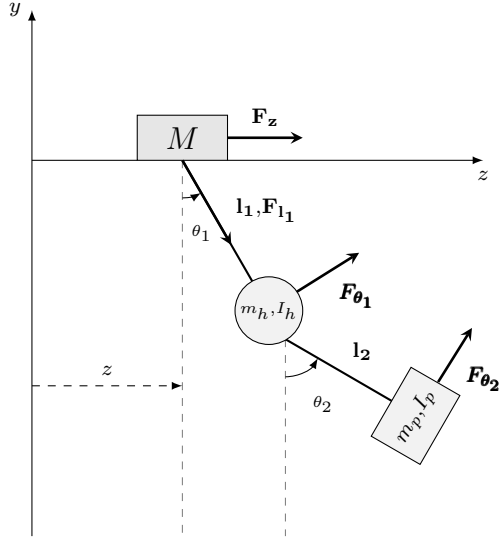}\else\includegraphics{schematic_double_pendulum.eps}\fi
\caption{Schematic diagram of the 2-D double-pendulum-type overhead crane system. The~angular torque control $F_{\theta_2}\equiv0$ in Section~\ref{without_Ft2}.} 
\label{fig:M1}        
\end{figure}

In practice, the industrial overhead cranes may exhibit double-pendulum swing phenomena, due to many factors, such as heavy payload and non-negligible hook masses. As depicted in Fig.~\ref{fig:M1}, a D-P-T overhead crane is operated to move along the bridge to and from, to transport the payload to the target position. In Fig.~\ref{fig:M1},  $z$ is the trolley position, $l_1$ and $l_2$ represent the hoisting rope length and the connecting rope length from the hook mass-center to the payload mass-center, respectively ($l_1$ is variable, $l_2$ is assumed to be constant). The center of mass of the hook and payload subtend angles $\theta_1$ and $\theta_2$ with respect to the direction of the negative $y$ axis, $F_z$ denotes the driven force applied to the trolley equipped with a winding/unwinding mechanism, $F_{l_1}$ represents the force applied to the hoisting rope, and $F_{\theta_1},$ $F_{\theta_2}$ are the angular torque controls imposed on the angles $\theta_1,$ $\theta_2,$ respectively. 

The equations of motion may be derived using Lagrangian dynamics.  The Lagrangian is $L = KE - PE$ , where $KE$ is the kinetic energy and $PE$ is the potential energy of the system. The kinetic energy is the sum of the kinetic energies,
\[
KE_{tr} = \frac12M\dot z^2
\]
of the trolley and the kinetic energies of the the center of mass of the hook and payload, each of which has a linear and a rotational component (\cite{Rafat}):
\[
KE_h + KE_p = \frac12m_h\left(\dot z_h^2+\dot y_h^2\right)+\frac12I_h\dot\theta_1^2+\frac12m_p\left(\dot z_p^2+\dot y_p^2\right)+\frac12I_p\dot\theta_2^2,
\]
where
\[
z_h=z+l_1\sin(\theta_1),\ y_h=-l_1\cos(\theta_1)
\]
and
\[
z_p=z+l_1\sin(\theta_1)+l_2\sin(\theta_2),\ y_p = -l_1\cos(\theta_1)-l_2\cos(\theta_2),
\]
and $I_h$ and $I_p$ are the moments of inertia about the center of mass. The potential energy is the sum of the potential energy of each subsystem, that is,
\[
PE =PE_{tr}+PE_h+PE_p= gMy_{tr}+gm_hy_h+gm_py_p+V_0,
\]
where $V_0$ is a suitable reference potential.
If we express $y_h$ and $y_p$ in the terms of $l_1,\theta_1$ and $l_2,\theta_2$ as above, and choose $V_0$ so that the potential energy of the system is zero for $\theta_1=\pi/2$ and $\theta_2=\pi/2,$ we obtain
\[
PE = -gm_h\left[l_1\cos(\theta_1)\right]-gm_p\left[l_1\cos(\theta_1) + l_2\cos(\theta_2)\right],
\]
where $g$ is the gravitational acceleration. Here we assume that the potential energy of the trolley is kept unchanged ($y_{tr}\equiv0$) and  we also neglect the effect of friction on the system behavior. The payload is assumed to have uniform mass density and the ropes are assumed to be massless.

The D-P-T overhead crane system under consideration consists of four independent generalized coordinates namely the trolley position, $z$, the rope length between the trolley and the hook, $l_1,$ hook angle, $\theta_1,$ and payload angle, $\theta_2.$ $M,$ $m_h,$ $m_p,$ and $l_2$ represent the trolley mass, hook mass, payload mass, and rope length between the hook and the payload respectively.

The Lagrange equation with respect to the coordinate $q_i$ is
\[
\frac{d}{dt}\left(\frac{\partial L}{\partial\dot q_i}\right)-\frac{\partial L}{\partial q_i}=F_i,
\]
where $L,$ $q_i$ ($i = 1, 2, 3, 4$) and $F_i$ represent the Lagrangian function, generalized coordinates ($q_1,$ $q_2,$ $q_3$ and $q_4$ represents $z,$ $l_1,$ $\theta_1$ and $\theta_2,$ respectively) and the forces not arising from a potential, respectively. For general theory, see, e.g. \cite{Goldstein}.

The dynamic model of the D-P-T overhead crane, according to individual generalized coordinates, is as follows: 

For a generalized coordinate $q_1=z:$
\begin{equation*}
M{\ddot z} + m_h{\ddot z} + m_p{\ddot z} + m_h\sin({\theta_1}){\ddot l_1} + m_p\sin({\theta_1}){\ddot l_1} - m_h\sin({\theta_1}){l_1}{\dot\theta_1}^2 - m_p\sin({\theta_1}){l_1}{\dot\theta_1}^2 
\end{equation*}
\begin{equation}\label{lagrange_z}
+ m_h\cos({\theta_1}){l_1}{\ddot\theta_1} + m_p\cos({\theta_1}){l_1}{\ddot\theta_1} + 2m_h\cos({\theta_1}){\dot l_1}{\dot\theta_1} + 2m_p\cos({\theta_1}){\dot l_1}{\dot\theta_1} \tag{$L_z$}
\end{equation}
\begin{equation*}
- l_2m_p\sin({\theta_2}){\dot\theta_2}^2 + l_2m_p\cos({\theta_2}){\ddot\theta_2}=F_z  
\end{equation*}
For a generalized coordinate $q_2=l_1:$
\[
m_h{\ddot l_1} + m_p{\ddot l_1} + m_h\sin({\theta_1}){\ddot z} + m_p\sin({\theta_1}){\ddot z} - gm_h\cos({\theta_1}) - gm_p\cos({\theta_1}) - m_h{l_1}{\dot\theta_1}^2
\]
\begin{equation}\label{lagrange_l_1}
 - m_p{l_1}{\dot\theta_1}^2 - l_2m_p\cos({\theta_1} - {\theta_2}){\dot\theta_2}^2 + l_2m_p\sin({\theta_1} - {\theta_2}){\ddot\theta_2}=F_{l_1} \tag{$L_{l_1}$}
\end{equation}
For a generalized coordinate $q_3=\theta_1:$
\[
I_h{\ddot\theta_1} + m_h{l^2_1}{\ddot\theta_1} + m_p{l^2_1}{\ddot\theta_1} + m_h\cos({\theta_1}){l_1}{\ddot z} + m_p\cos({\theta_1}){l_1}{\ddot z} + 2m_h{l_1}{\dot l_1}{\dot\theta_1} + 2m_p{l_1}{\dot l_1}{\dot\theta_1} + gm_h\sin({\theta_1}){l_1} 
\]
\begin{equation}\label{lagrange_theta_1}
+ gm_p\sin({\theta_1}){l_1} + l_2m_p{l_1}\sin({\theta_1} - {\theta_2}){\dot\theta_2}^2 + l_2m_p{l_1}\cos({\theta_1} - {\theta_2}){\ddot\theta_2}=F_{\theta_1} \tag{$L_{\theta_1}$}
\end{equation}
For a generalized coordinate $q_4=\theta_2:$
\[
I_p{\ddot\theta_2} + l_2^2m_p{\ddot\theta_2} + l_2m_p\sin({\theta_1} - {\theta_2}){\ddot l_1} + l_2m_p\cos({\theta_2}){\ddot z} + gl_2m_p\sin({\theta_2}) - l_2m_p{l_1}\sin({\theta_1}- {\theta_2}){\dot\theta_1}^2 
\]
\begin{equation}\label{lagrange_theta_2}
 + l_2m_p{l_1}\cos({\theta_1} - {\theta_2}){\ddot\theta_1} + 2l_2m_p\cos({\theta_1} - {\theta_2}){\dot l_1}{\dot\theta_1}=F_{\theta_2}, \tag{$L_{\theta_2}$}
\end{equation}
where for the comparison purpose of the performance of the designed state feedback-based feed-forward control, we will analyze for the last equation both cases. The first case, without considering the control force $F_{\theta_2},$ that is, the swing angle $\theta_2$ of the payload is not controlled, in Section~\ref{without_Ft2} and, second case, with considering the control force $F_{\theta_2},$ in Section~\ref{with_Ft2}.

\section{Theoretical background. Control law design}\label{section:bacground}

\begin{figure}
\ifpdf\includegraphics[trim={3cm 16cm 3cm 0},clip]{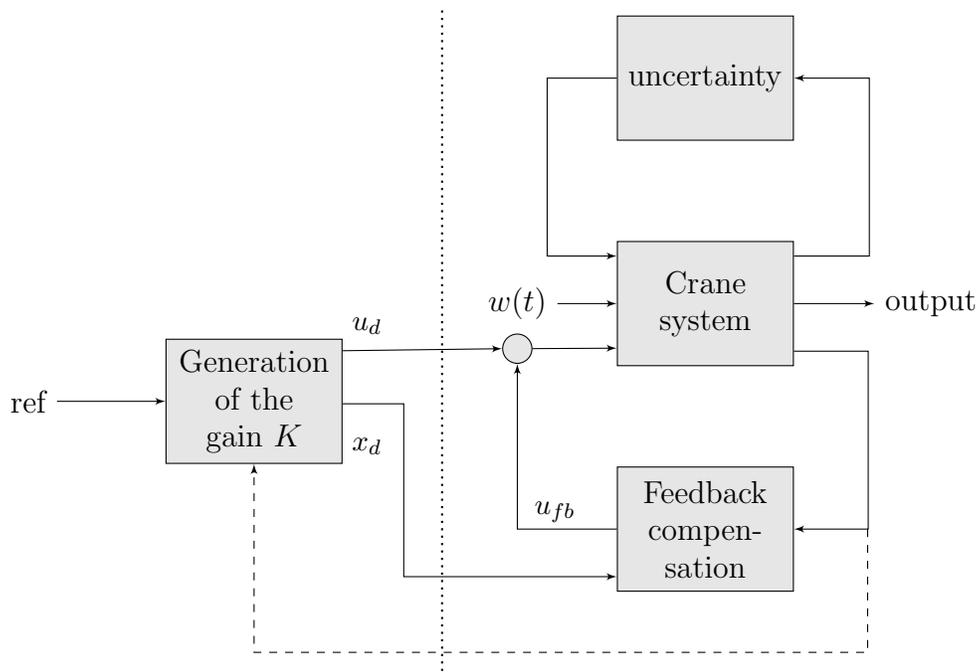}\else\includegraphics{schematic_control.eps}\fi
\caption{Schematic diagram of the control of double-pendulum-type overhead crane system. Here we can identify two subsystems, (i) feed-forward control generating the gain matrix $K$ asymptotically stabilizing the system at the desired end position and (ii) feedback loop compensating the disturbance $w(t),$  unmodeled dynamics or/and a parametric uncertainty. The input 'ref' to the control system denotes a target position of the payload, the signal $u_d$ represents the control forces $F_z,$ $F_{l_1},$ $F_{\theta_1}$ (Section~\ref{without_Ft2}) and eventually also $F_{\theta_2}$ (Section~\ref{with_Ft2}) and $x_d$ represents the desired  trajectory generated by the reference model.}
\label{fig:scheme}       
\end{figure}

 The control strategies may be in principle divided into two basic approaches, an open loop and closed-loop control  techniques, see, e.~g. \cite{sun}, \cite{Fujioka1}, \cite{Fujioka2} for open-loop  and \cite{Qian_book}, \cite{Zhang}, \cite{Qian}, \cite{Tuan}, \cite{Sun3} for closed-loop control approach, and the references therein. Each of these approaches has its advantages and disadvantages, feedback systems have the advantage of being simple, requires minimal knowledge about the process to be controlled; in particular, a mathematical model 
of the process is not required. The system measures the variables, and uses that variables to make decisions, compensates for inaccuracies in the reference model, measurement error, and unmeasured disturbances. Feed-forward systems, on the other hand, have the ability to anticipate changes in the measured variable, working proactively instead of reactively. Feed-forward systems are only as good as the reference model with which the system works. The system cannot consider an unmeasured variable when making its decisions, and these blind spots can cause control to break down. Therefore by combining feed-forward with a feedback control, we obtain a system providing a backup level of control and increases the control efficiency. The feed-forward control is responsible for the approximate course while the feedback controller adjust small deviations from reference trajectory generated by feed-forward subsystem by the external disturbance and/or unmodeled system dynamics. Fig.~\ref{fig:scheme} illustrates a block diagram of a hybrid control strategy for the crane systems in general.  

Numerous control strategies have been developed to suppress the payload sway in order to improve the safety and positioning accuracy of the D-P-T overhead crane operations.

For example, in \cite{sun}, an energy-optimal trajectory planner has been presented for double pendulum crane systems with various state constraints based on the convexification of the energy consumption function and the discretization of the system dynamics. The control problem is reformulated into a quadratic programming problem, whose solution is obtained by using a convex optimization tool.

Fujioka and Singhose in \cite{Fujioka1} and \cite{Fujioka2} compared the performances of various input-shaped
model reference control designs applied to a D-P-T crane system with uncertainty.

In the paper \cite{Zhang}, an adaptive tracking control method for double-pendulum overhead crane systems subject to tracking error limitation, parametric uncertainties and external disturbances and with S-shaped smooth function as the
trolley desired trajectory is presented. 

In \cite{Qian}, the authors developed  dynamically-weighted SIRMs "single-input-rule modules"-based fuzzy controller, designed according to human experience.

The paper \cite{Tuan} focuses on the design of robust nonlinear controllers based on both conventional and hierarchical sliding mode techniques for D-P-T overhead crane systems. 

A nonlinear adaptive antiswing control scheme for the D-P-T overhead cranes subject to uncertain/unknown parameters has been presented in \cite{Sun3}.

The most of the control methods for D-P-T overhead crane systems are mainly devoted to the design of the regulation control scheme while the trajectory generation stage is usually not considered. Also, all of the aforementioned works do not consider the coordinate $l_1$ as a state variable, which is the simplest way to avoid the singularity in the system dynamics modeling, for the details see Remark~\ref{remark_singularity} and Remark~\ref{rmk_approx} below, indicating the essential originality of our study. Although in the paper \cite{Alhazza} the length of the hoisting cable, $l_1,$ is considered variable, but only for a substantially simplified mathematical model an iterative learning control technique is proposed to generate acceleration profiles of D-P-T overhead crane maneuvers. 

In the present paper we focus our attention on the design of the feedback-based feed-forward controller asymptotically stabilizing the D-P-T crane system in the desired position of payload, which may be combined as backup control level with the methods mentioned before, according to the scheme on the Fig.\ref{fig:scheme}.

Our approach to the asymptotic stabilization of the D-P-T overhead crane system is based on the two cornerstones of the modern control theory and theory of dynamical systems.

Consider a linear time-invariant (LTI) system $\dot x=Ax+Bu,$ where $A$ and $B$ are $n\times n$ and $n\times m$ constant real matrices, respectively. A fundamental result of linear control theory is that the following three conditions are equivalent, see, e.~g. \cite{AntMich}: 
\begin{itemize}
\item[(i)] the pair $(A,B)$ is controllable;
\item[(ii)] $\rank\mathcal{C}_{(A,B)}=n,$ where $\mathcal{C}_{(A,B)}=:\left( B\  AB\  A^2B\ \cdots\ A^{n-1}B\right)$ is an $n\times mn$ Kalman's controllability matrix;
\item[(iii)] for every $n$-tuple real and/or complex conjugate numbers $\lambda_1,$ $\lambda_2,$ $\dots,$ $\lambda_n,$ there exists an $m\times n$ state feedback gain matrix $K$ such that the eigenvalues of the closed-loop system matrix $A_{cl}=A-BK$
are the desired values $\lambda_1,\lambda_2,\dots,\lambda_n.$
\end{itemize}

In general, we will study the nonlinear control system
\begin{equation*}
\dot x=G(x,u), \ t\geq 0,\ x\in\mathbb{R}^n\ u\in\mathbb{R}^m, \ \dot{} =:d/dt,
\end{equation*}
with the state feedback of the form $u=-Kx$ and for which is assumed that $x=0$ is its solution, that is, $G(0,0)=0.$ It is well-known, that if the pair $(A,B),$ where $A=G_x(0,0)$  and $B=G_u(0,0)$ are the corresponding Jacobian matrices with respect to the state and input variables, respectively, and evaluated at $(0,0)$  is controllable, then the LTI system $\dot x=(A-BK)x$ is in some neighborhood of the origin topologically equivalent, and preserving the parametrization by time, to the system $\dot x=G(x,-Kx),$ provided that the eigenvalues of the matrix $A-BK$ have non-zero real part. The precise statement about this property gives the Hartman-Grobman theorem, see, e.~g. \cite[p.~120]{Perko}.  Thus, if we choose the matrix $K$ such that all eigenvalues of $A-BK$ have negative real parts, the nonlinear system $\dot x=G(x,-Kx)$ is locally asymptotically stable in the neighborhood of $x=0.$

Here we mean the usual definition of local asymptotic stability, that is, we say, that solution $x=0$ of the system $\dot x=G(x,-Kx)$ is asymptotically stable if for every $\varepsilon>0$ there exists $\delta>0$ such that if $||x(0)||\leq\delta$ then $||x(t)||\leq\varepsilon$ for all $t\geq 0,$ and, moreover, $||x(t)||\rightarrow 0^+$ with $t\rightarrow \infty$ ($||\cdot||$ denotes an $n-$dimensional Euclidean norm).

Now we will prove the local asymptotic stability of the closed-loop system $\dot x=G(x,-Kx)$ around its equilibrium position $x=0$ by computing the lower bound of the region of attraction.

Let the gain matrix $K$ is such that the real parts of all eigenvalues of the matrix $A_{cl}=A-BK$ are  negative. Then 
\[
\dot x=G(x,-Kx)=A_{cl}x+R_1(x),
\]
where $R_1(x)$ is the first-order Taylor's remainder, obviously $||R_1(x)||=O(||x||^2)$ as $||x||\rightarrow 0^+.$ This implies that there exists a constant $\sigma=\sigma(K)>0$ such that $||R_1(x)||\leq\sigma(K)||x||$ for all $x$ with $||x||\leq\Gamma(K).$  Let us consider as a Lyapunov function candidate $V(x)=x^TPx,$ where symmetric and  positive definite matrix $P$ is a solution of Lyapunov equation 
\[
PA_{cl}+A^T_{cl}P=-Q(K)
\]
for appropriate choice of the symmetric and positive definite matrix $Q,$ which can be solved as an optimization problem with regards to the gain matrix $K.$ Let the constant $\sigma(K)$ is such that 
\[
\lambda_{\min}(Q(K))>2\sigma(K)\lambda_{\max}(P),
\]
where $\lambda_{\min}$ and $\lambda_{\max}$ denote the minimal and maximal eigenvalue of the matrix, respectively.
Then along the trajectories of the system $\dot x=G(x,-Kx)=:g_K(x)$ we have
\[
\dot V(x(t))=x^T(t)Pg_K(x(t))+g^T_K(x(t))Px(t)
\]
\[
=x^TP[A_{cl}x+R_1(x)]+[x^TA^T_{cl}+R^T_1(x)]Px
\]
\[
=x^T(PA_{cl}+A^T_{cl}P)x+2x^TPR_1(x)=-x^TQx+2x^TPR_1(x).
\]
Because for each $n\times n$ symmetric and positive definite real matrix $C$ is
\[
\lambda_{\min}(C)||x||^2\leq x^TCx\leq\lambda_{\max}(C)||x||^2,\ x\in\mathbb{R}^n,
\]
we get that
\[
\dot V(x(t))\leq-\left( \lambda_{\min}(Q(K))-2\sigma(K)\lambda_{\max}(P) \right)||x||^2<0
\]
for all $x$ satisfying $||x||\leq\Gamma(K),$ $x\neq0,$ which implies local asymptotic stability of the zero solution of the closed-loop system $\dot x=G(x,-Kx).$

The general strategy in stabilizing the system under consideration at the desired end position is as follows:
\begin{itemize}
\item First, we show that the linear part of mathematical model of the D-P-T crane system derived above is controllable in the Kalman's sense;
\item Secondly, we establish the range of the permissible eigenvalues $\lambda_i$ of the closed-loop system to be the system asymptotically stable in the work area of the D-P-T crane taking into account its technical limitations, for example, the maximum permitted velocities of the trolley and hoist device, jerk, {\it etc.} are limited by their individual technical parameters and this range of the permissible eigenvalues can be calculated and tabulated from the results of the simulation experiments with known weights and applied forces.
\end{itemize}
To keep the situation clear and without loss of generality we will hereafter assume that the desired end position of payload is $x=0,$ which corresponds in Fig.~\ref{fig:M1}  to $z=0,$ $y=-l_2,$ $\theta_1=0,$ $\theta_2=0,$ and also all velocities ($\dot z,$ $\dot l_1,$ $\dot\theta_1,$ and $\dot\theta_2$) are the zeros.

In the following section we will analyze the underactuated crane system with $F_{\theta_2}\equiv 0$ in the fourth equation derived from the Lagrange equation. Subsequently, in Section~\ref{with_Ft2}, we will deal with the same system but with considering the swing angle control of payload.

\section{Crane system without considering the angular torque control $F_{\theta_2}$}\label{without_Ft2}

Solving the equations (\ref{lagrange_z}), (\ref{lagrange_l_1}), (\ref{lagrange_theta_1}) and (\ref{lagrange_theta_2}) with regard to the variable $\ddot z,$ $\ddot l_1,$ $\ddot\theta_1,$ $\ddot\theta_2,$
and substituting $x_1,$ $x_2,$ $x_3,$ $x_4,$
$x_5,$ $x_6,$ $x_7$ and $x_8$ instead of $z,$ $l_1,$ $\theta_1,$ $\theta_2,$ $\dot z,$ $\dot l_1,$ $\dot\theta_1$ and $\dot\theta_2,$ respectively, and with $F_{\theta_2}\equiv0,$ we obtain the functions $G_i,$ $i=5,6,7,8.$
First let us denote
\[
\Delta=: I_{h}I_{p}{m^2_{h}}+I_{h}I_{p}{m^2_{p}}+I_{h}M{l^2_{2}}{m^2_{p}}+2I_{p}M{m^2_{h}}{x^2_{2}}+2I_{p}M{m^2_{p}}{x^2_{2}}+2I_{h}I_{p}Mm_{h}+2I_{h}I_{p}Mm_{p}+I_{h}{l^2_{2}}m_{h}{m^2_{p}}
\]
\[
+I_{h}{l^2_{2}}{m^2_{h}}m_{p}+2I_{h}I_{p}m_{h}m_{p}+I_{h}I_{p}{m^2_{h}}\cos\left(2x_{3}\right)+I_{h}I_{p}{m^2_{p}}\cos\left(2x_{3}\right)+I_{h}M{l^2_{2}}{m^2_{p}}\cos\left(2x_{3}-2x_{4}\right)
\]
\[
+2I_{h}I_{p}m_{h}m_{p}\cos\left(2x_{3}\right)+2M{l^2_{2}}m_{h}{m^2_{p}}{x^2_{2}}+2M{l^2_{2}}{m^2_{h}}m_{p}{x^2_{2}}+2I_{h}M{l^2_{2}}m_{h}m_{p}
\]
\[
+4I_{p}Mm_{h}m_{p}{x^2_{2}}+I_{h}{l^2_{2}}m_{h}{m^2_{p}}\cos\left(2x_{3}\right)+I_{h}{l^2_{2}}{m^2_{h}}m_{p}\cos\left(2x_{3}\right).
\]
Then
\[
G_5=
 -\frac{1}{\Delta}\times\Bigg[2{F_{\theta_1}}I_{p}{m^2_{h}}x_{2}\cos\left(x_{3}\right)-2{F_z}I_{p}{m^2_{h}}{x^2_{2}}-2{F_z}I_{p}{m^2_{p}}{x^2_{2}}
\]
\[
-2{F_z}I_{h}I_{p}m_{h}-2{F_z}I_{h}I_{p}m_{p}-{F_z}I_{h}{l^2_{2}}{m^2_{p}}+2{F_{\theta_1}}I_{p}{m^2_{p}}x_{2}\cos\left(x_{3}\right)+I_{h}I_{p}g{m^2_{h}}\sin\left(2x_{3}\right)
\]
\[
+I_{h}I_{p}g{m^2_{p}}\sin\left(2x_{3}\right)+{F_{l_1}}I_{h}{l^2_{2}}{m^2_{p}}\sin\left(x_{3}\right)+2{F_{l_1}}I_{p}{m^2_{h}}{x^2_{2}}\sin\left(x_{3}\right)+2{F_{l_1}}I_{p}{m^2_{p}}{x^2_{2}}\sin\left(x_{3}\right)
\]
\[
-2{F_z}{l^2_{2}}m_{h}{m^2_{p}}{x^2_{2}}-2{F_z}{l^2_{2}}{m^2_{h}}m_{p}{x^2_{2}}+2{F_{l_1}}I_{h}I_{p}m_{h}\sin\left(x_{3}\right)+2{F_{l_1}}I_{h}I_{p}m_{p}\sin\left(x_{3}\right)
\]
\[
-2{F_z}I_{h}{l^2_{2}}m_{h}m_{p}-4{F_z}I_{p}m_{h}m_{p}{x^2_{2}}-{F_{l_1}}I_{h}{l^2_{2}}{m^2_{p}}\sin\left(x_{3}-2x_{4}\right)-{F_z}I_{h}{l^2_{2}}{m^2_{p}}\cos\left(2x_{3}-2x_{4}\right)
\]
\[
+2{F_{\theta_1}}{l^2_{2}}m_{h}{m^2_{p}}x_{2}\cos\left(x_{3}\right)+2{F_{\theta_1}}{l^2_{2}}{m^2_{h}}m_{p}x_{2}\cos\left(x_{3}\right)+I_{h}{l^3_{2}}m_{h}{m^2_{p}}{x^2_{8}}\sin\left(2x_{3}-x_{4}\right)
\]
\[
+4{F_{\theta_1}}I_{p}m_{h}m_{p}x_{2}\cos\left(x_{3}\right)+I_{h}g{l^2_{2}}m_{h}{m^2_{p}}\sin\left(2x_{3}\right)+I_{h}g{l^2_{2}}{m^2_{h}}m_{p}\sin\left(2x_{3}\right)
\]
\[
+2{F_{l_1}}{l^2_{2}}m_{h}{m^2_{p}}{x^2_{2}}\sin\left(x_{3}\right)+2{F_{l_1}}{l^2_{2}}{m^2_{h}}m_{p}{x^2_{2}}\sin\left(x_{3}\right)-I_{h}{l^3_{2}}m_{h}{m^2_{p}}{x^2_{8}}\sin\left(x_{4}\right)
\]
\[
+2I_{h}I_{p}gm_{h}m_{p}\sin\left(2x_{3}\right)+2{F_{l_1}}I_{h}{l^2_{2}}m_{h}m_{p}\sin\left(x_{3}\right)+I_{h}I_{p}l_{2}{m^2_{p}}{x^2_{8}}\sin\left(2x_{3}-x_{4}\right)
\]
\[
+4I_{h}I_{p}{m^2_{h}}x_{6}x_{7}\cos\left(x_{3}\right)+4I_{h}I_{p}{m^2_{p}}x_{6}x_{7}\cos\left(x_{3}\right)+4{F_{l_1}}I_{p}m_{h}m_{p}{x^2_{2}}\sin\left(x_{3}\right)-I_{h}I_{p}l_{2}{m^2_{p}}{x^2_{8}}\sin\left(x_{4}\right)
\]
\[
+4I_{h}{l^2_{2}}m_{h}{m^2_{p}}x_{6}x_{7}\cos\left(x_{3}\right)+4I_{h}{l^2_{2}}{m^2_{h}}m_{p}x_{6}x_{7}\cos\left(x_{3}\right)+I_{h}I_{p}l_{2}m_{h}m_{p}{x^2_{8}}\sin\left(2x_{3}-x_{4}\right)
\]
\[
+8I_{h}I_{p}m_{h}m_{p}x_{6}x_{7}\cos\left(x_{3}\right)-I_{h}I_{p}l_{2}m_{h}m_{p}{x^2_{8}}\sin\left(x_{4}\right)\Bigg],
\]

\[
G_6=
 \frac{1}{\Delta}\times\Bigg[2{F_{l_1}}I_{h}I_{p}M+{F_{l_1}}I_{h}{l^2_{2}}{m^2_{p}}+{F_{l_1}}I_{p}{m^2_{h}}{x^2_{2}}+{F_{l_1}}I_{p}{m^2_{p}}{x^2_{2}}
\]
\[
+2{F_{l_1}}I_{h}I_{p}m_{h}+2{F_{l_1}}I_{h}I_{p}m_{p}+2I_{h}I_{p}g{m^2_{h}}\cos\left(x_{3}\right)+2I_{h}I_{p}g{m^2_{p}}\cos\left(x_{3}\right)+I_{h}I_{p}{m^2_{h}}x_{2}{x^2_{7}}+I_{h}I_{p}{m^2_{p}}x_{2}{x^2_{7}}
\]
\[
+{F_{l_1}}M{l^2_{2}}{m^2_{p}}{x^2_{2}}+2I_{p}M{m^2_{h}}{x_{2}}^3{x^2_{7}}+2I_{p}M{m^2_{p}}{x_{2}}^3{x^2_{7}}-{F_z}I_{h}{l^2_{2}}{m^2_{p}}\sin\left(x_{3}\right)+2{F_{l_1}}I_{h}M{l^2_{2}}m_{p}
\]
\[
+{F_{\theta_1}}I_{p}{m^2_{h}}x_{2}\sin\left(2x_{3}\right)+{F_{\theta_1}}I_{p}{m^2_{p}}x_{2}\sin\left(2x_{3}\right)-2{F_z}I_{p}{m^2_{h}}{x^2_{2}}\sin\left(x_{3}\right)
\]
\[
-2{F_z}I_{p}{m^2_{p}}{x^2_{2}}\sin\left(x_{3}\right)+2{F_{l_1}}I_{p}Mm_{h}{x^2_{2}}+2{F_{l_1}}I_{p}Mm_{p}{x^2_{2}}+{F_{l_1}}{l^2_{2}}m_{h}{m^2_{p}}{x^2_{2}}
\]
\[
+{F_{l_1}}{l^2_{2}}{m^2_{h}}m_{p}{x^2_{2}}-2{F_z}I_{h}I_{p}m_{h}\sin\left(x_{3}\right)-2{F_z}I_{h}I_{p}m_{p}\sin\left(x_{3}\right)+2{F_{l_1}}I_{h}{l^2_{2}}m_{h}m_{p}
\]
\[
+2{F_{l_1}}I_{p}m_{h}m_{p}{x^2_{2}}+{F_z}I_{h}{l^2_{2}}{m^2_{p}}\sin\left(x_{3}-2x_{4}\right)-{F_{l_1}}I_{h}{l^2_{2}}{m^2_{p}}\cos\left(2x_{4}\right)
\]
\[
-{F_{l_1}}I_{p}{m^2_{h}}{x^2_{2}}\cos\left(2x_{3}\right)-{F_{l_1}}I_{p}{m^2_{p}}{x^2_{2}}\cos\left(2x_{3}\right)+I_{h}Mg{l^2_{2}}{m^2_{p}}\cos\left(x_{3}-2x_{4}\right)
\]
\[
+4I_{h}I_{p}gm_{h}m_{p}\cos\left(x_{3}\right)+I_{h}I_{p}l_{2}{m^2_{p}}{x^2_{8}}\cos\left(x_{3}-x_{4}\right)+2I_{h}I_{p}m_{h}m_{p}x_{2}{x^2_{7}}+I_{h}{l^3_{2}}m_{h}{m^2_{p}}{x^2_{8}}\cos\left(x_{3}+x_{4}\right)
\]
\[
+I_{h}I_{p}{m^2_{h}}x_{2}{x^2_{7}}\cos\left(2x_{3}\right)+I_{h}I_{p}{m^2_{p}}x_{2}{x^2_{7}}\cos\left(2x_{3}\right)+2M{l^2_{2}}m_{h}{m^2_{p}}{x_{2}}^3{x^2_{7}}+2M{l^2_{2}}{m^2_{h}}m_{p}{x_{2}}^3{x^2_{7}}
\]
\[
+{F_{\theta_1}}{l^2_{2}}m_{h}{m^2_{p}}x_{2}\sin\left(2x_{3}\right)+{F_{\theta_1}}{l^2_{2}}{m^2_{h}}m_{p}x_{2}\sin\left(2x_{3}\right)-2{F_z}{l^2_{2}}m_{h}{m^2_{p}}{x^2_{2}}\sin\left(x_{3}\right)
\]
\[
-2{F_z}{l^2_{2}}{m^2_{h}}m_{p}{x^2_{2}}\sin\left(x_{3}\right)+2{F_{l_1}}M{l^2_{2}}m_{h}m_{p}{x^2_{2}}+2I_{h}M{l^3_{2}}{m^2_{p}}{x^2_{8}}\cos\left(x_{3}-x_{4}\right)
\]
\[
+4I_{p}Mm_{h}m_{p}{x_{2}}^3{x^2_{7}}-2{F_z}I_{h}{l^2_{2}}m_{h}m_{p}\sin\left(x_{3}\right)+{F_{\theta_1}}M{l^2_{2}}{m^2_{p}}x_{2}\sin\left(2x_{3}-2x_{4}\right)
\]
\[
+2{F_{\theta_1}}I_{p}m_{h}m_{p}x_{2}\sin\left(2x_{3}\right)-4{F_z}I_{p}m_{h}m_{p}{x^2_{2}}\sin\left(x_{3}\right)+I_{h}{l^3_{2}}m_{h}{m^2_{p}}{x^2_{8}}\cos\left(x_{3}-x_{4}\right)
\]
\[
+I_{h}I_{p}l_{2}{m^2_{p}}{x^2_{8}}\cos\left(x_{3}+x_{4}\right)-{F_{l_1}}{l^2_{2}}m_{h}{m^2_{p}}{x^2_{2}}\cos\left(2x_{3}\right)-{F_{l_1}}{l^2_{2}}{m^2_{h}}m_{p}{x^2_{2}}\cos\left(2x_{3}\right)
\]
\[
+I_{h}Mg{l^2_{2}}{m^2_{p}}\cos\left(x_{3}\right)-{F_{l_1}}M{l^2_{2}}{m^2_{p}}{x^2_{2}}\cos\left(2x_{3}-2x_{4}\right)-2{F_{l_1}}I_{p}m_{h}m_{p}{x^2_{2}}\cos\left(2x_{3}\right)
\]
\[
+2I_{p}Mg{m^2_{h}}{x^2_{2}}\cos\left(x_{3}\right)+2I_{p}Mg{m^2_{p}}{x^2_{2}}\cos\left(x_{3}\right)+I_{h}M{l^2_{2}}{m^2_{p}}x_{2}{x^2_{7}}+2I_{h}I_{p}Mgm_{h}\cos\left(x_{3}\right)
\]
\[
+2I_{h}I_{p}Mgm_{p}\cos\left(x_{3}\right)+2I_{h}g{l^2_{2}}m_{h}{m^2_{p}}\cos\left(x_{3}\right)+2I_{h}g{l^2_{2}}{m^2_{h}}m_{p}\cos\left(x_{3}\right)
\]
\[
+2I_{h}I_{p}{m^2_{h}}x_{6}x_{7}\sin\left(2x_{3}\right)+2I_{h}I_{p}{m^2_{p}}x_{6}x_{7}\sin\left(2x_{3}\right)+2I_{h}I_{p}Mm_{h}x_{2}{x^2_{7}}+2I_{h}I_{p}Mm_{p}x_{2}{x^2_{7}}
\]
\[
+I_{h}{l^2_{2}}m_{h}{m^2_{p}}x_{2}{x^2_{7}}+I_{h}{l^2_{2}}{m^2_{h}}m_{p}x_{2}{x^2_{7}}+2I_{h}M{l^2_{2}}{m^2_{p}}x_{6}x_{7}\sin\left(2x_{3}-2x_{4}\right)+2I_{h}I_{p}Ml_{2}m_{p}{x^2_{8}}\cos\left(x_{3}-x_{4}\right)
\]
\[
+4I_{h}I_{p}m_{h}m_{p}x_{6}x_{7}\sin\left(2x_{3}\right)+I_{h}I_{p}l_{2}m_{h}m_{p}{x^2_{8}}\cos\left(x_{3}-x_{4}\right)+I_{h}{l^2_{2}}m_{h}{m^2_{p}}x_{2}{x^2_{7}}\cos\left(2x_{3}\right)
\]
\[
+I_{h}{l^2_{2}}{m^2_{h}}m_{p}x_{2}{x^2_{7}}\cos\left(2x_{3}\right)+I_{h}M{l^2_{2}}{m^2_{p}}x_{2}{x^2_{7}}\cos\left(2x_{3}-2x_{4}\right)+2I_{h}I_{p}m_{h}m_{p}x_{2}{x^2_{7}}\cos\left(2x_{3}\right)
\]
\[
+2M{l^3_{2}}m_{h}{m^2_{p}}{x^2_{2}}{x^2_{8}}\cos\left(x_{3}-x_{4}\right)+I_{h}I_{p}l_{2}m_{h}m_{p}{x^2_{8}}\cos\left(x_{3}+x_{4}\right)+2Mg{l^2_{2}}m_{h}{m^2_{p}}{x^2_{2}}\cos\left(x_{3}\right)
\]
\[
+2Mg{l^2_{2}}{m^2_{h}}m_{p}{x^2_{2}}\cos\left(x_{3}\right)+2I_{h}Mg{l^2_{2}}m_{h}m_{p}\cos\left(x_{3}\right)+2I_{h}{l^2_{2}}m_{h}{m^2_{p}}x_{6}x_{7}\sin\left(2x_{3}\right)
\]
\[
+2I_{h}{l^2_{2}}{m^2_{h}}m_{p}x_{6}x_{7}\sin\left(2x_{3}\right)+4I_{p}Mgm_{h}m_{p}{x^2_{2}}\cos\left(x_{3}\right)+2I_{h}M{l^2_{2}}m_{h}m_{p}x_{2}{x^2_{7}}
\]
\[
+2I_{p}Ml_{2}{m^2_{p}}{x^2_{2}}{x^2_{8}}\cos\left(x_{3}-x_{4}\right)+2I_{p}Ml_{2}m_{h}m_{p}{x^2_{2}}{x^2_{8}}\cos\left(x_{3}-x_{4}\right)\Bigg],
\]

\[
G_7=
 \frac{1}{\Delta}\times\Bigg[{F_{\theta_1}}I_{p}{m^2_{h}}+{F_{\theta_1}}I_{p}{m^2_{p}}+{F_{\theta_1}}M{l^2_{2}}{m^2_{p}}+2{F_{\theta_1}}I_{p}Mm_{h}+2{F_{\theta_1}}I_{p}Mm_{p}
\]
\[
+{F_{\theta_1}}{l^2_{2}}m_{h}{m^2_{p}}+{F_{\theta_1}}{l^2_{2}}{m^2_{h}}m_{p}+2{F_{\theta_1}}I_{p}m_{h}m_{p}+{F_{\theta_1}}I_{p}{m^2_{h}}\cos\left(2x_{3}\right)+{F_{\theta_1}}I_{p}{m^2_{p}}\cos\left(2x_{3}\right)
\]
\[
+2{F_{\theta_1}}I_{p}m_{h}m_{p}\cos\left(2x_{3}\right)-2{F_z}I_{p}{m^2_{h}}x_{2}\cos\left(x_{3}\right)-2{F_z}I_{p}{m^2_{p}}x_{2}\cos\left(x_{3}\right)
\]
\[
+{F_{l_1}}I_{p}{m^2_{h}}x_{2}\sin\left(2x_{3}\right)+{F_{l_1}}I_{p}{m^2_{p}}x_{2}\sin\left(2x_{3}\right)+2{F_{\theta_1}}M{l^2_{2}}m_{h}m_{p}
\]
\[
+{F_{\theta_1}}{l^2_{2}}m_{h}{m^2_{p}}\cos\left(2x_{3}\right)+{F_{\theta_1}}{l^2_{2}}{m^2_{h}}m_{p}\cos\left(2x_{3}\right)+{F_{\theta_1}}M{l^2_{2}}{m^2_{p}}\cos\left(2x_{3}-2x_{4}\right)
\]
\[
-2{F_z}{l^2_{2}}m_{h}{m^2_{p}}x_{2}\cos\left(x_{3}\right)-2{F_z}{l^2_{2}}{m^2_{h}}m_{p}x_{2}\cos\left(x_{3}\right)-4{F_z}I_{p}m_{h}m_{p}x_{2}\cos\left(x_{3}\right)
\]
\[
-4I_{p}M{m^2_{h}}x_{2}x_{6}x_{7}-4I_{p}M{m^2_{p}}x_{2}x_{6}x_{7}+{F_{l_1}}{l^2_{2}}m_{h}{m^2_{p}}x_{2}\sin\left(2x_{3}\right)+{F_{l_1}}{l^2_{2}}{m^2_{h}}m_{p}x_{2}\sin\left(2x_{3}\right)
\]
\[
+{F_{l_1}}M{l^2_{2}}{m^2_{p}}x_{2}\sin\left(2x_{3}-2x_{4}\right)+2{F_{l_1}}I_{p}m_{h}m_{p}x_{2}\sin\left(2x_{3}\right)-2I_{p}Mg{m^2_{h}}x_{2}\sin\left(x_{3}\right)
\]
\[
-2I_{p}Mg{m^2_{p}}x_{2}\sin\left(x_{3}\right)-4M{l^2_{2}}m_{h}{m^2_{p}}x_{2}x_{6}x_{7}-4M{l^2_{2}}{m^2_{h}}m_{p}x_{2}x_{6}x_{7}-2M{l^3_{2}}m_{h}{m^2_{p}}x_{2}{x^2_{8}}\sin\left(x_{3}-x_{4}\right)
\]
\[
-8I_{p}Mm_{h}m_{p}x_{2}x_{6}x_{7}-2Mg{l^2_{2}}m_{h}{m^2_{p}}x_{2}\sin\left(x_{3}\right)-2Mg{l^2_{2}}{m^2_{h}}m_{p}x_{2}\sin\left(x_{3}\right)
\]
\[
-4I_{p}Mgm_{h}m_{p}x_{2}\sin\left(x_{3}\right)-2I_{p}Ml_{2}{m^2_{p}}x_{2}{x^2_{8}}\sin\left(x_{3}-x_{4}\right)-2I_{p}Ml_{2}m_{h}m_{p}x_{2}{x^2_{8}}\sin\left(x_{3}-x_{4}\right)\Bigg]
\]
and
\[
G_8=
-\frac{1}{\Delta}\times\Bigg[\sin\left(x_{3}-x_{4}\right)\bigg(l_{2}m_{p}\left(2{F_{l_1}}Mm_{p}{x^2_{2}}+2{F_{l_1}}I_{h}M+{F_{l_1}}I_{h}m_{p}\right)
+l_{2}m_{h}m_{p}\left(2{F_{l_1}}M{x^2_{2}}+{F_{l_1}}I_{h}\right)\bigg)
\]
\[
+l_{2}m_{p}\bigg[I_{h}Ml_{2}m_{p}\sin\left(2x_{3}-2x_{4}\right){x^2_{8}}-{F_{l_1}}I_{h}m_{p}\sin\left(x_{3}+x_{4}\right)+{F_z}I_{h}m_{p}\cos\left(x_{4}\right)
\]
\[
+{F_z}I_{h}m_{p}\cos\left(2x_{3}-x_{4}\right)+I_{h}Mgm_{p}\sin\left(2x_{3}-x_{4}\right)+I_{h}Mgm_{p}\sin\left(x_{4}\right)
\]
\[
+2{F_{\theta_1}}Mm_{p}x_{2}\cos\left(x_{3}-x_{4}\right)+4I_{h}Mm_{p}x_{6}x_{7}\cos\left(x_{3}-x_{4}\right)\bigg]
\]
\[
+l_{2}m_{h}m_{p}\bigg({F_z}I_{h}\cos\left(2x_{3}-x_{4}\right)
-{F_{l_1}}I_{h}\sin\left(x_{3}+x_{4}\right)
\]
\[
+{F_z}I_{h}\cos\left(x_{4}\right)+I_{h}Mg\sin\left(x_{4}\right)
+2{F_{\theta_1}}Mx_{2}\cos\left(x_{3}-x_{4}\right)
\]
\[
+I_{h}Mg\sin\left(2x_{3}-x_{4}\right)+4I_{h}Mx_{6}x_{7}\cos\left(x_{3}-x_{4}\right)\bigg)\Bigg].
\]
Substituting into $G_i,$ $i=5,6,7,8,$ $u_1$ instead of $F_z,$  $u_3$ instead of $F_{\theta_1},$ and setting
\[
F_{l_1} = -\frac{1}{\Theta_{l_1}}\times\Bigg[2I_{h}I_{p}g{m^2_{h}}\cos\left(x_{3}\right)-u_{2}+2I_{h}I_{p}g{m^2_{p}}\cos\left(x_{3}\right)+I_{h}Mg{l^2_{2}}{m^2_{p}}\cos\left(x_{3}-2x_{4}\right)
\]
\[
+4I_{h}I_{p}gm_{h}m_{p}\cos\left(x_{3}\right)+I_{h}Mg{l^2_{2}}{m^2_{p}}\cos\left(x_{3}\right)+2I_{h}I_{p}Mgm_{h}\cos\left(x_{3}\right)
\]
\[
+2I_{h}I_{p}Mgm_{p}\cos\left(x_{3}\right)+2I_{h}g{l^2_{2}}m_{h}{m^2_{p}}\cos\left(x_{3}\right)+2I_{h}g{l^2_{2}}{m^2_{h}}m_{p}\cos\left(x_{3}\right)
\]
\begin{equation}\label{Fl1}
+2I_{h}Mg{l^2_{2}}m_{h}m_{p}\cos\left(x_{3}\right)\Bigg],
\end{equation}
where
\[
\Theta_{l_1}=:I_{h}\left(2I_{p}m_{h}+2I_{p}m_{p}+2I_{p}M+{l^2_{2}}{m^2_{p}}+2{l^2_{2}}m_{h}m_{p}-{l^2_{2}}{m^2_{p}}\cos\left(2x_{4}\right)+2M{l^2_{2}}m_{p}\right),
\]
we get the control system $\dot x=G(x,u),$ $G(x,u)=(x_5,x_6,x_7,x_8,G_5,G_6,G_7,G_8)^T,$ $u=(u_1,u_2,u_3),$ satisfying $G(0,0)=0.$ The control force $F_{l_1}$ is a solution of linear algebraic equation where on the left-hand side are collected the terms from the function $G_6$ that are nonzero for $x=0,$ $F_z=0,$ $F_{\theta_1}=0,$ (also $F_{\theta_2}=0$ in $\hat G_6$ in Section~\ref{with_Ft2}) and on the right-hand side is $u_2.$ Such $F_{l_1}$ will ensure that $G(0,0)=0$ ($\hat G(0,0)=0$ in Section~\ref{with_Ft2}). These eight, highly coupled equations of motion, capture full system dynamics.
\begin{rmk}\label{forces_asymptotics}
\rm Obviously, $u_1(t)=F_z(t)\rightarrow 0$ and $u_3(t)=F_{\theta_1}(t)\rightarrow 0$ for $t\rightarrow \infty$ because $u(t)=-Kx(t)\rightarrow 0$ for a suitably chosen matrix $K$ and initial state $x(0),$ $||x(0)||\leq\Gamma(K).$ To obtain the asymptotics of the control force $F_{l_1}(t),$ we set $x=0$ and $u_2=0$ into (\ref{Fl1}), and we have for $t\rightarrow \infty$ that
\[
F_{l_1}(t)\rightarrow \frac{-g\Bigg[2I_{p}\left(m_{h}+m_{p}\right)\left(M+m_{h}+m_{p}\right)+2{l^2_{2}}m_{p}\left(M+m_{h}\right)\left(m_{h}+m_{p}\right)\Bigg]}{2I_{p}\left(M+m_{h}+m_{p}\right)+2{l^2_{2}}m_{p}\left(M+m_{h}\right)}
\]
\[
=-g\left(m_{h}+m_{p}\right),
\]
which reflects the clear fact that control force $F_{l_1}$ must compensate the hook and payload weight. 
\end{rmk}
\begin{rmk}\label{remark_singularity}
\rm Here, it is worth mentioning the key moment of our study, namely, it would be much harder to analyze the just derived reference model if we did not take into account the mass moment of inertia of the hook, $I_h.$ Indeed, for $I_h\equiv0,$ the denominator $\Delta$ in $G_i,$ $i=5,6,7,8$ is of the order $O(x_2^2),$ that is, the control system $\dot x=G(x,u)$ would be singular in the sense that $|G_i|\rightarrow\infty$ for $x_2(=l_1)\rightarrow 0^+.$ This is probably the main reason why, to our best knowledge, there is no paper dealing analytically with the complete D-P-T crane system with the time-varying rope length between the trolley and hook, $l_1$ in our notation.
\end{rmk}
Now we check the controllability of the linear part. The corresponding jacobians are
\begingroup\makeatletter\def\f@size{8}\check@mathfonts
\[
G_x(0,0)=\left(\begin{array}{cccccccc} 0 & 0 & 0 & 0 & 1 & 0 & 0 & 0\\ 0 & 0 & 0 & 0 & 0 & 1 & 0 & 0\\ 0 & 0 & 0 & 0 & 0 & 0 & 1 & 0\\ 0 & 0 & 0 & 0 & 0 & 0 & 0 & 1\\ 0 & 0 & 0 & \frac{g{l^2_{2}}{m^2_{p}}}{I_{p}\Sigma_m+{l^2_{2}}m_{h}m_{p}+M{l^2_{2}}m_{p}} & 0 & 0 & 0 & 0\\ 0 & 0 & 0 & 0 & 0 & 0 & 0 & 0\\ 0 & 0 & 0 & 0 & 0 & 0 & 0 & 0\\ 0 & 0 & 0 & -\frac{gl_{2}m_{p}\Sigma_m}{I_{p}\Sigma_m+{l^2_{2}}m_{h}m_{p}+M{l^2_{2}}m_{p}} & 0 & 0 & 0 & 0 \end{array}\right)
\]
and
\[
G_u(0,0)=\left(\begin{array}{ccc} 0 & 0 & 0\\ 0 & 0 & 0\\ 0 & 0 & 0\\ 0 & 0 & 0\\ \frac{m_{p}{l^2_{2}}+I_{p}}{I_{p}\Sigma_m+{l^2_{2}}m_{h}m_{p}+M{l^2_{2}}m_{p}} & 0 & 0\\ 0 & \frac{1}{2I_{h}\left(m_{h}+m_{p}\right)\left(I_{p}\Sigma_m+{l^2_{2}}m_{h}m_{p}+M{l^2_{2}}m_{p}\right)} & 0\\ 0 & 0 & \frac{1}{I_{h}}\\ -\frac{l_{2}m_{p}}{I_{p}\Sigma_m+{l^2_{2}}m_{h}m_{p}+M{l^2_{2}}m_{p}} & 0 & 0 \end{array}\right),
\]
\endgroup
where $\Sigma_m=:M+m_{h}+m_{p}.$

\begin{figure} 
   \centerline{
    \hbox{
     \psfig{file=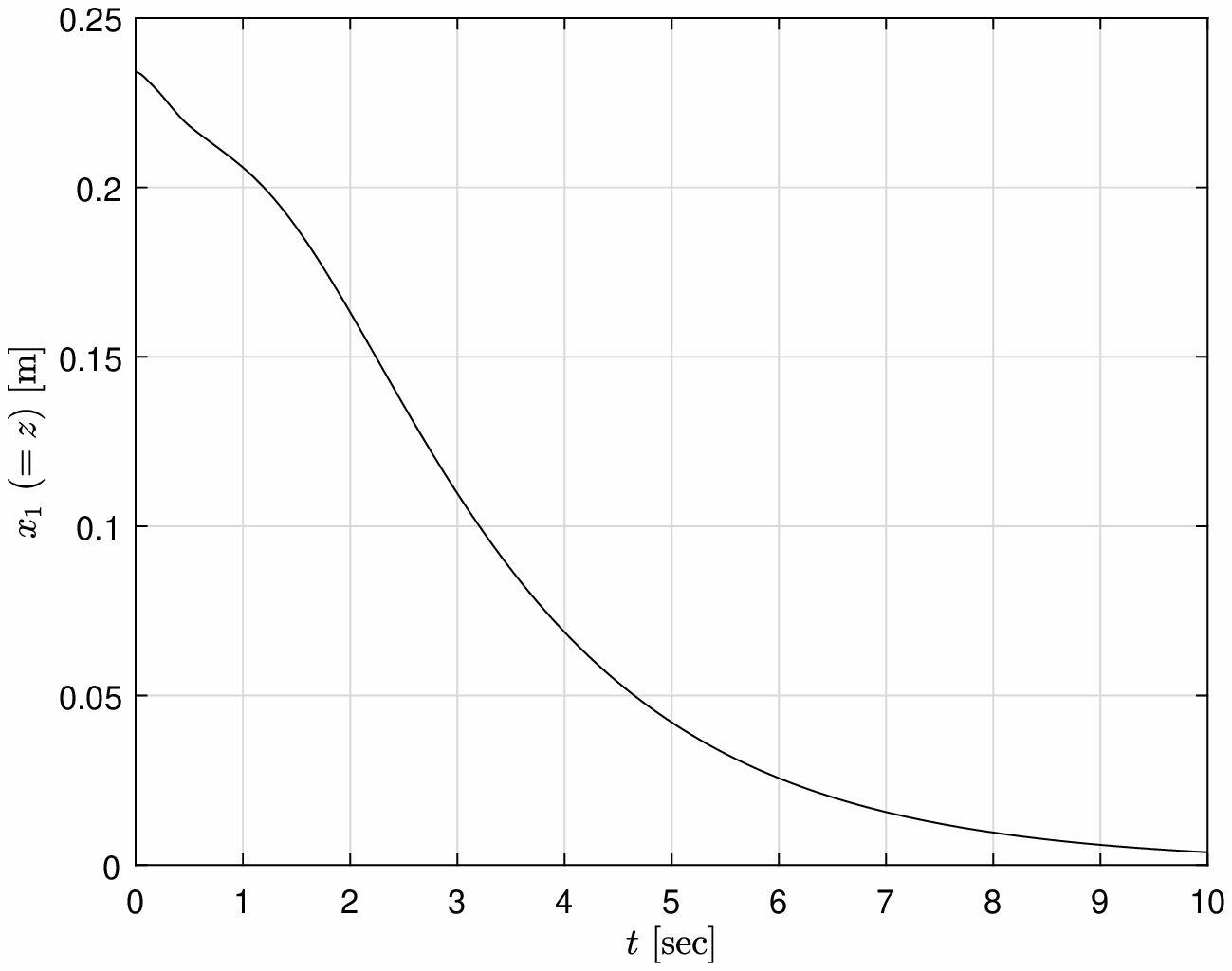,width=5.cm, clip=}
     \hspace{1.cm}
     \psfig{file=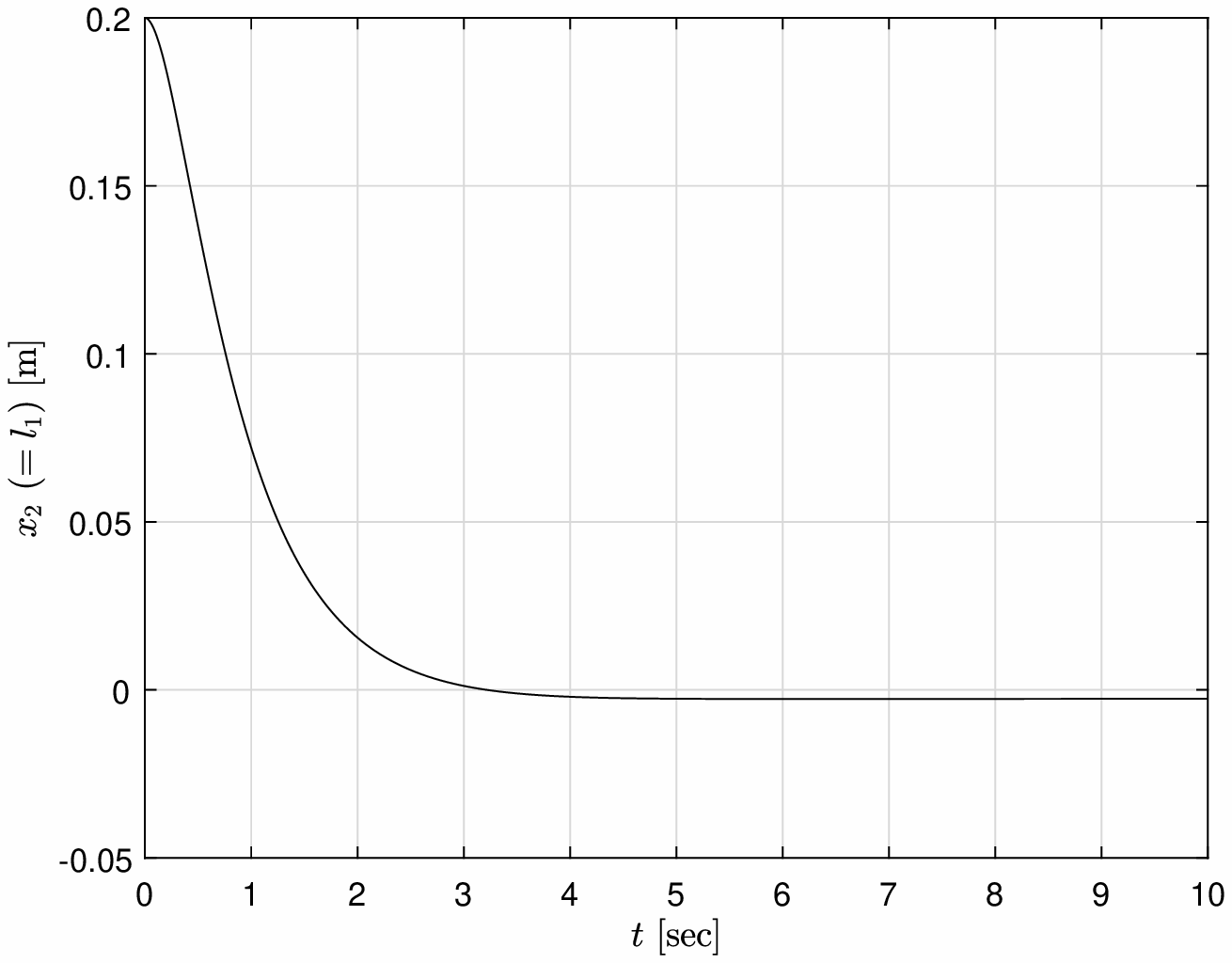,width=5cm,clip=}
    }
   }
   \vspace{0.5cm}
   \centerline{
    \hbox{
     \psfig{file=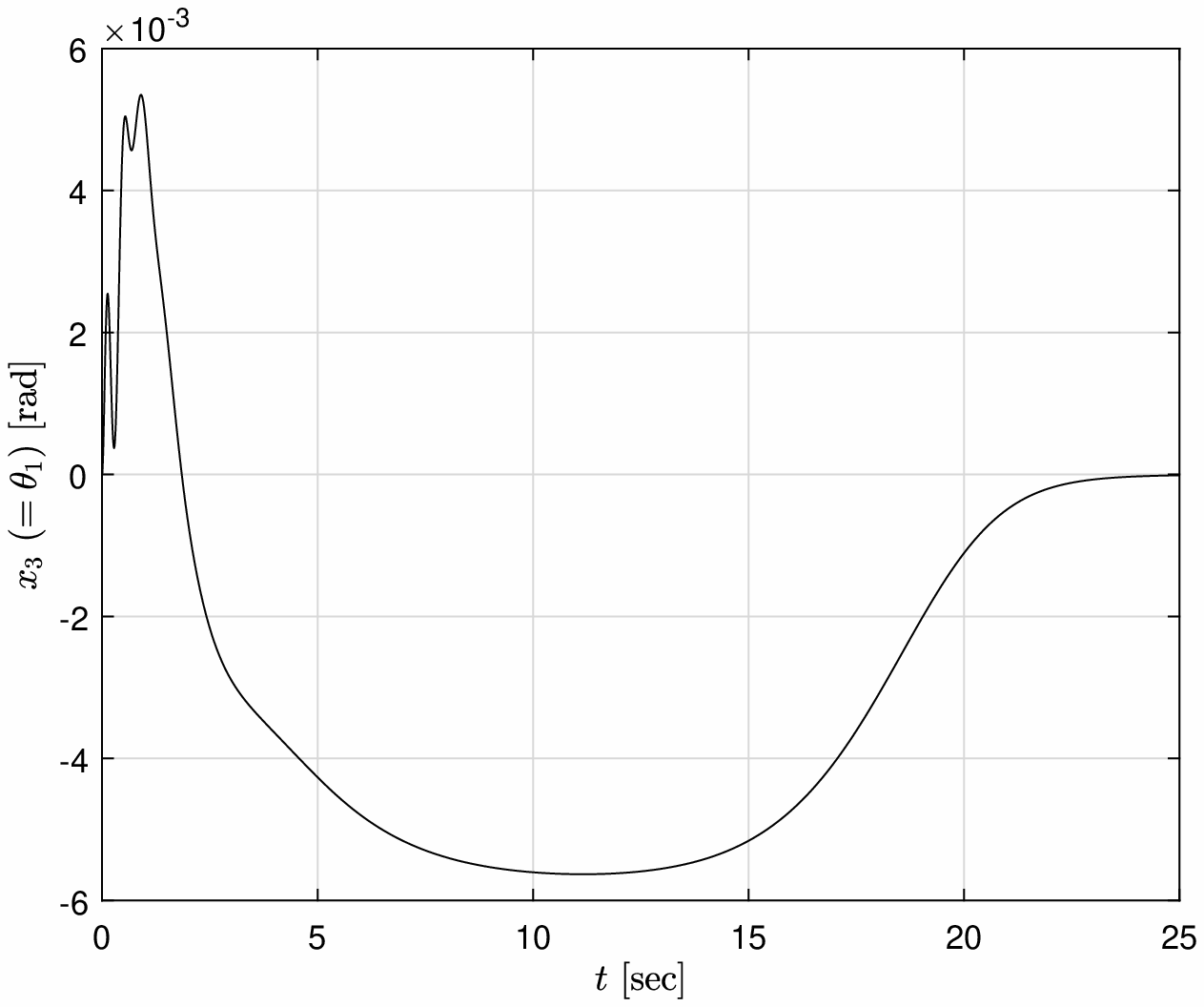,width=5.cm,clip=}
     \hspace{1.cm}
     \psfig{file=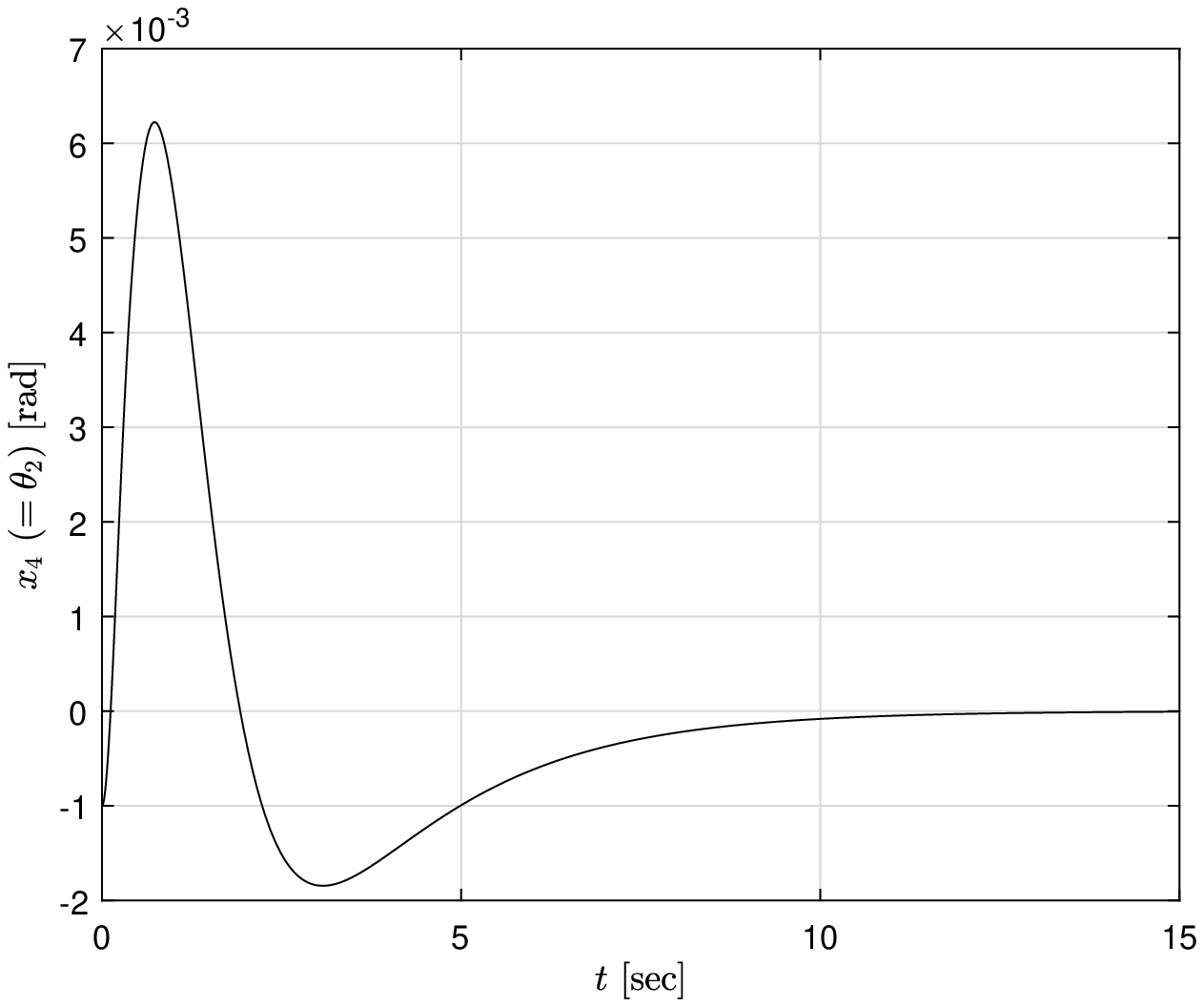,width=5.cm,clip=}
    }
   }
   \vspace{0.5cm}
   \centerline{
    \hbox{
     \psfig{file=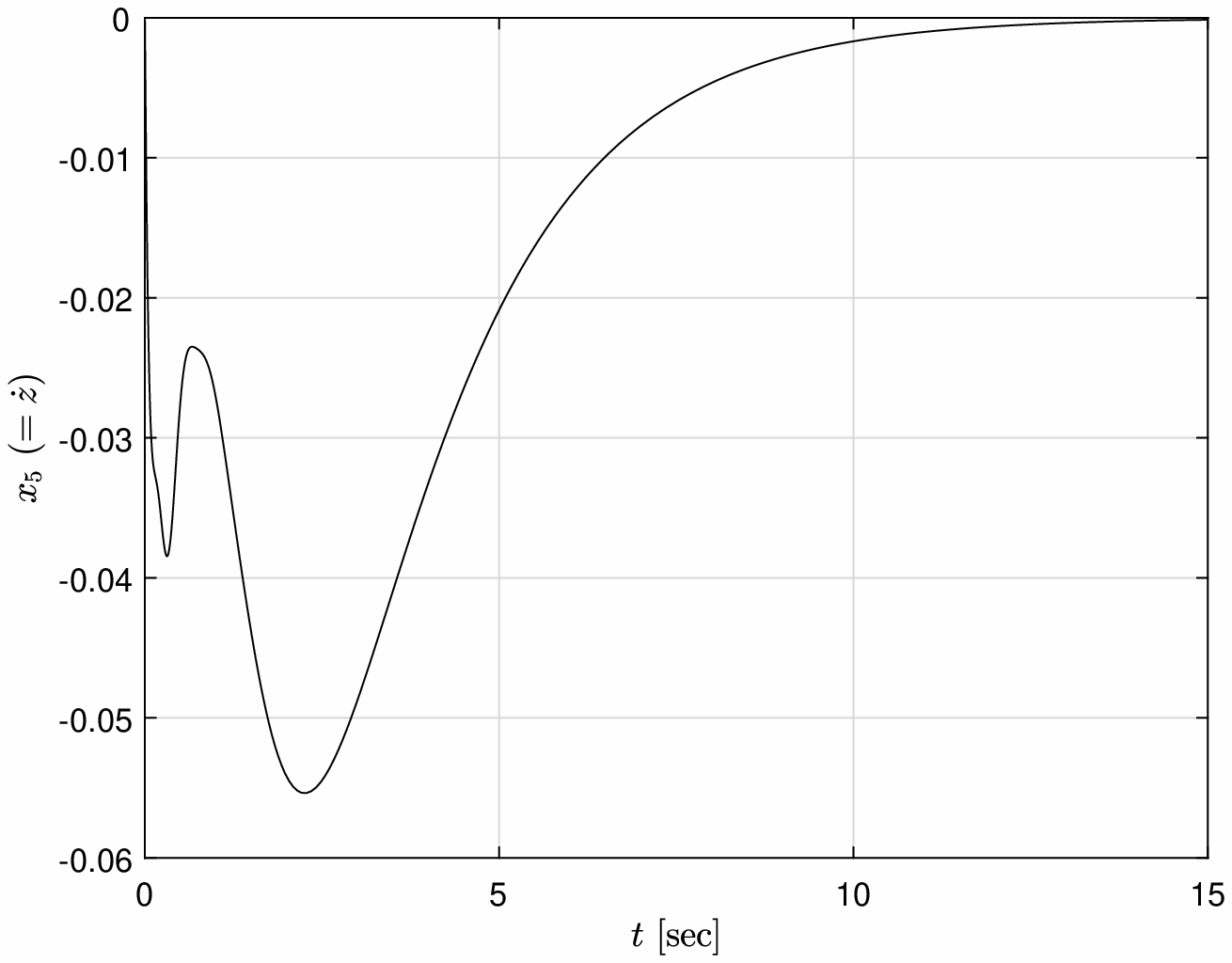,width=5.cm,clip=}
     \hspace{1.cm}
     \psfig{file=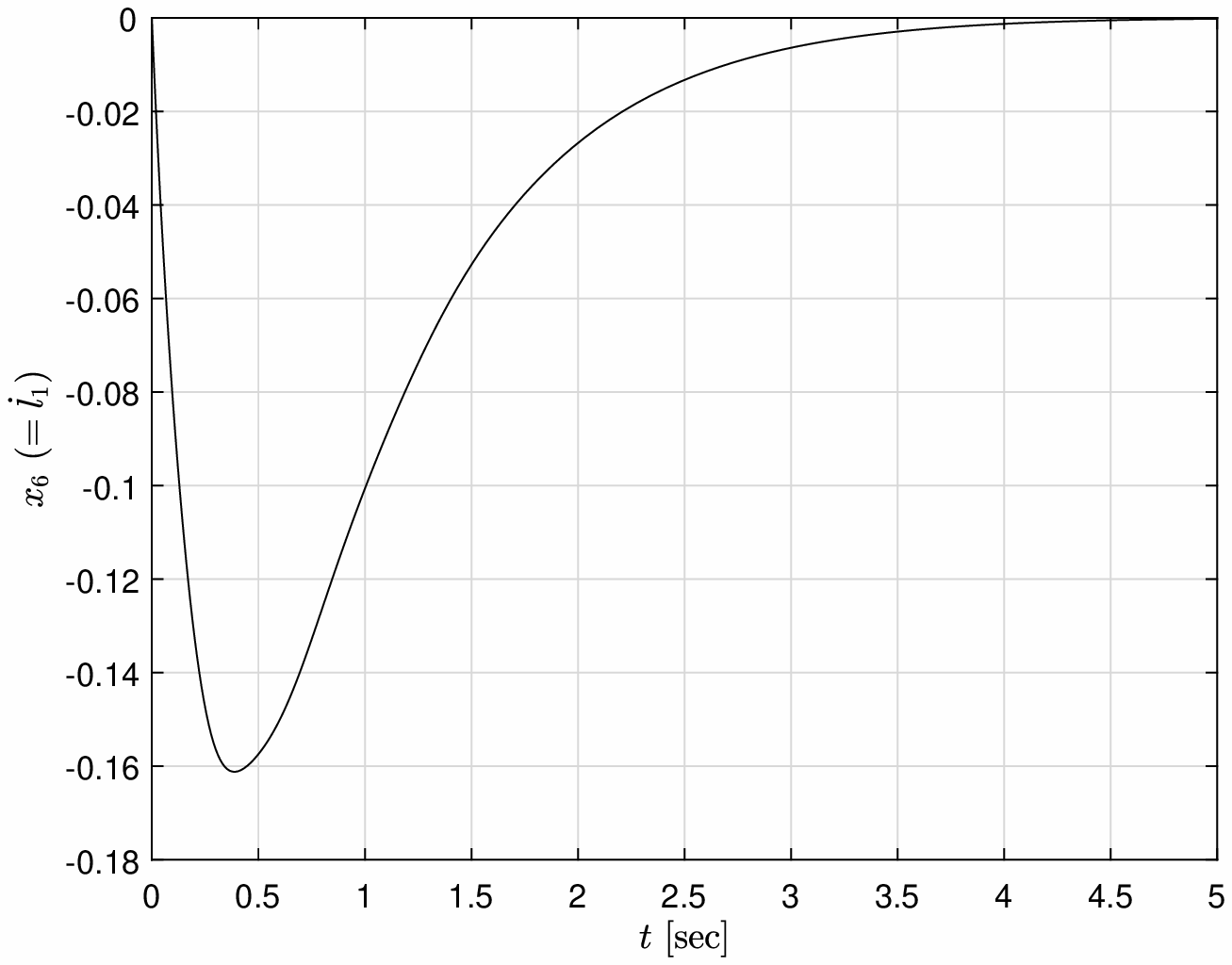,width=5.cm,clip=}
         }
   }
   \vspace{0.5cm}
   \centerline{
    \hbox{
     \psfig{file=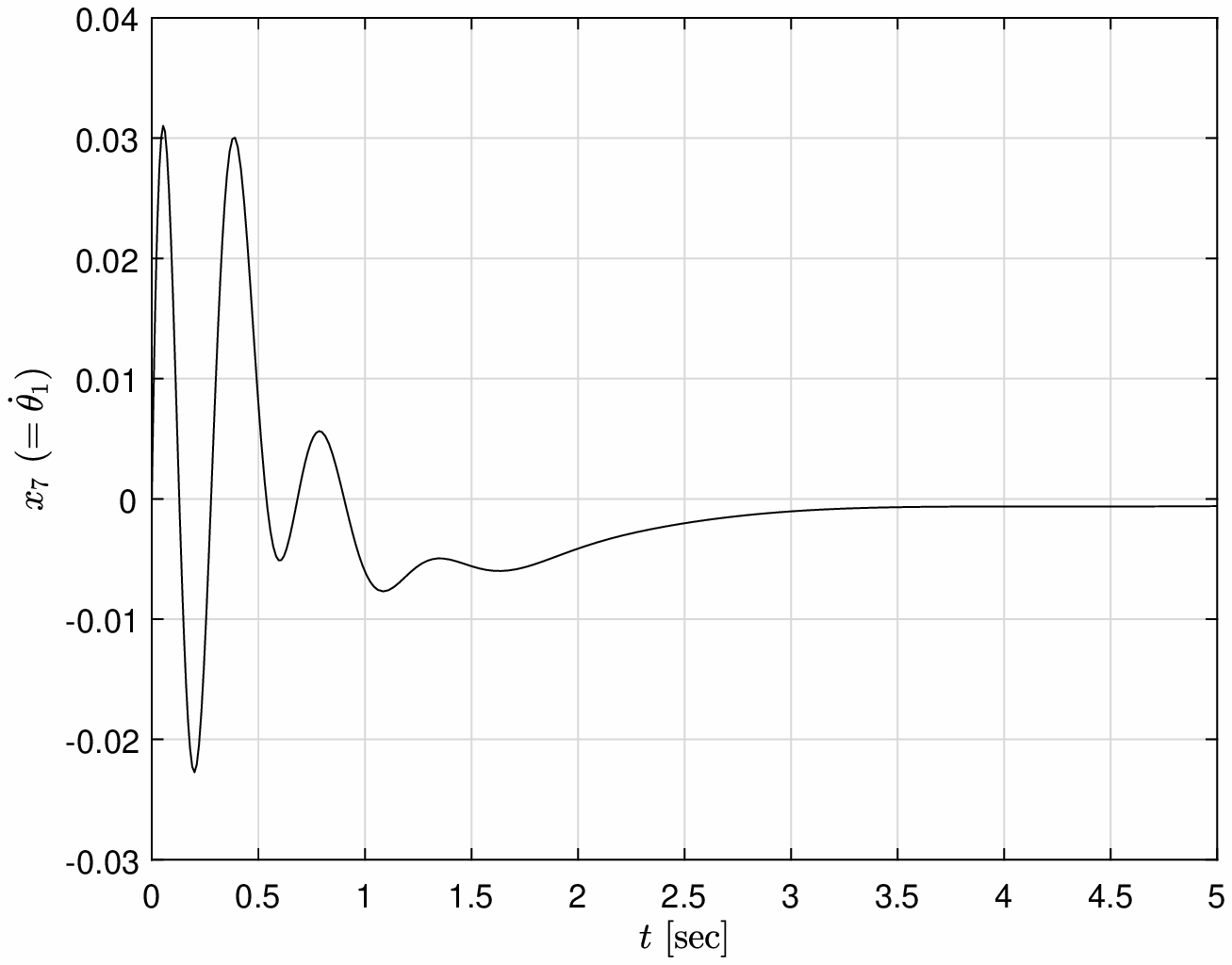,width=5.cm,clip=}
     \hspace{1.cm}
     \psfig{file=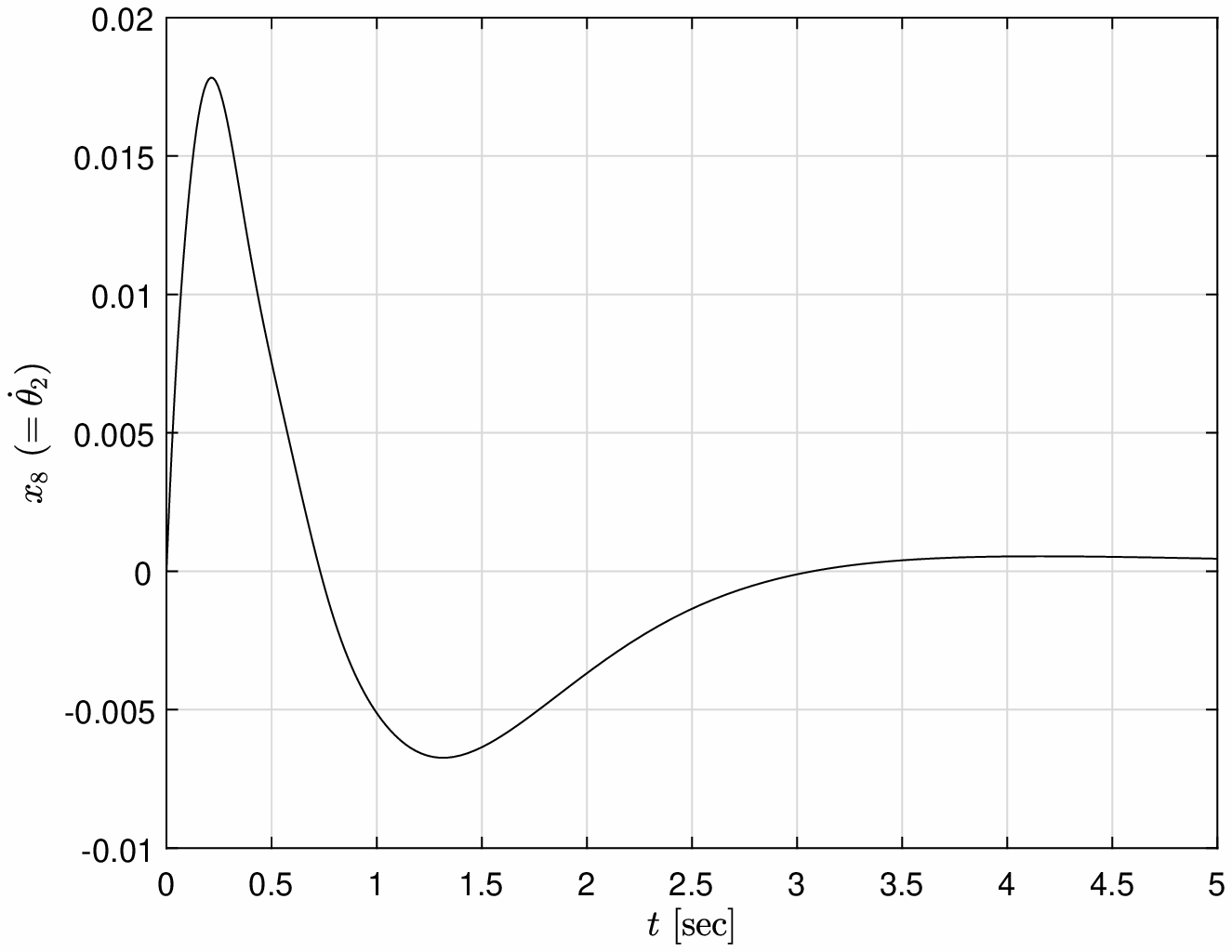,width=5.cm,clip=}
         }
   }
\caption{The solutions $ x_i,$ $i=1,\dots,8$ of the reference model {\bf without considering} the angular torque control $F_{\theta_2}$ with the initial state $(0.234\ 0.2\ 0\ -0.001\ 0\ 0\ 0\ 0)^T.$ Corresponding control forces $F_z,$ $F_{l_1}$ and $F_{\theta_1}$ are depicted in the Fig.\ref{control_forces_without_Ft2}.}
\label{solutions_xi_without_Ft2}
\end{figure}

\begin{figure} 
   \centerline{
    \hbox{
     \psfig{file=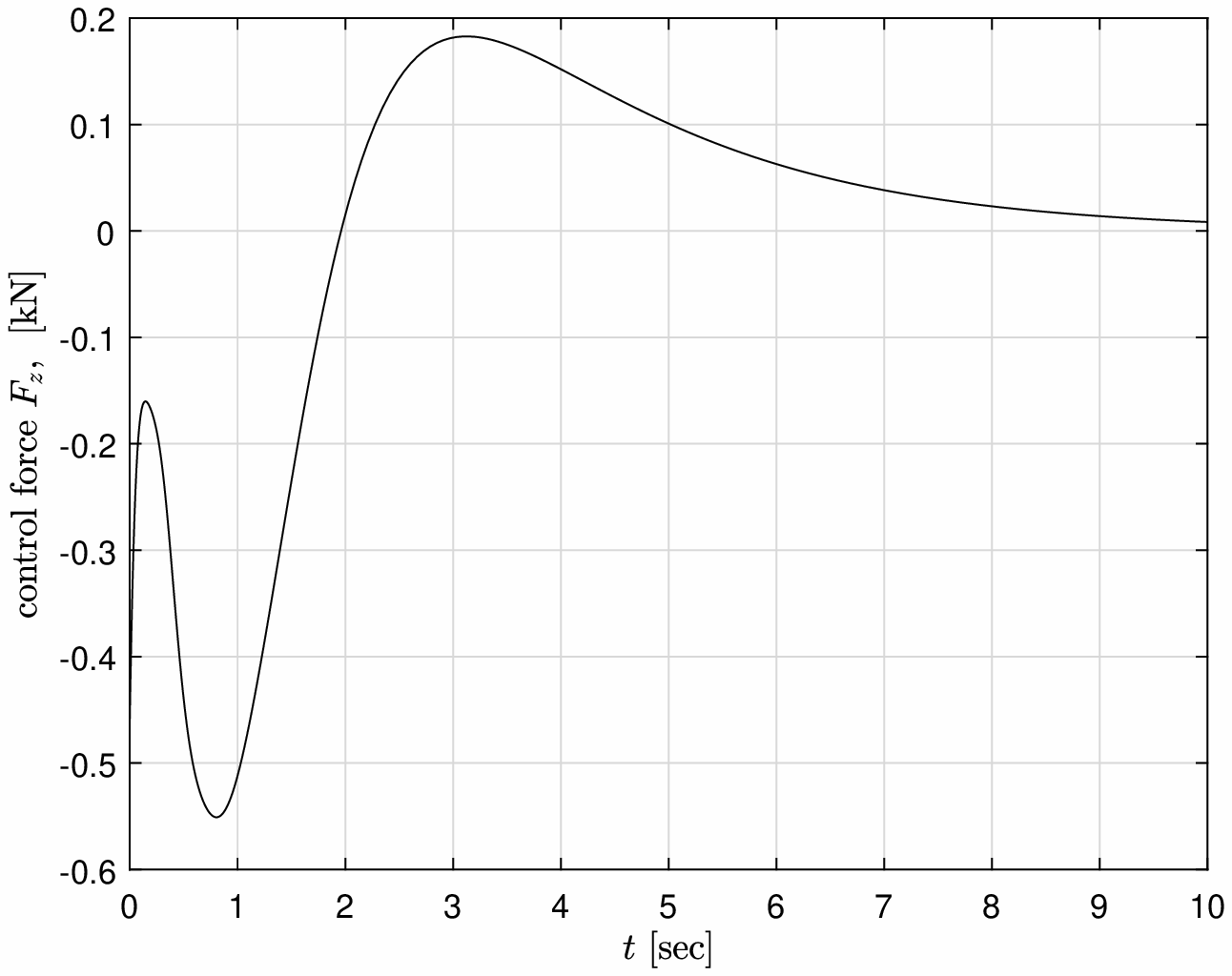,width=5.cm, clip=}
     \hspace{1.cm}
     \psfig{file=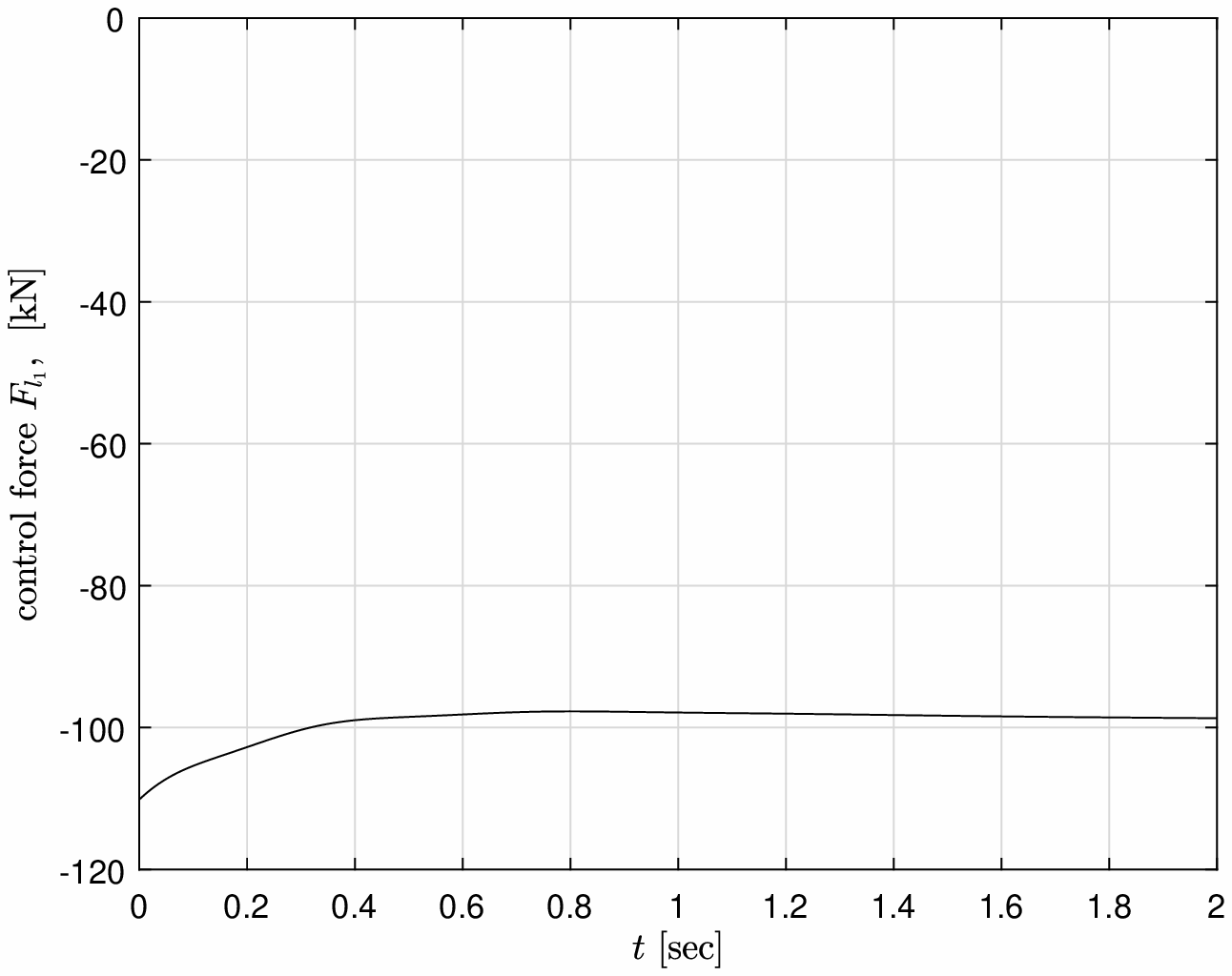,width=5cm,clip=}
    }
   }
   \vspace{0.5cm}
   \centerline{
    \hbox{
     \psfig{file=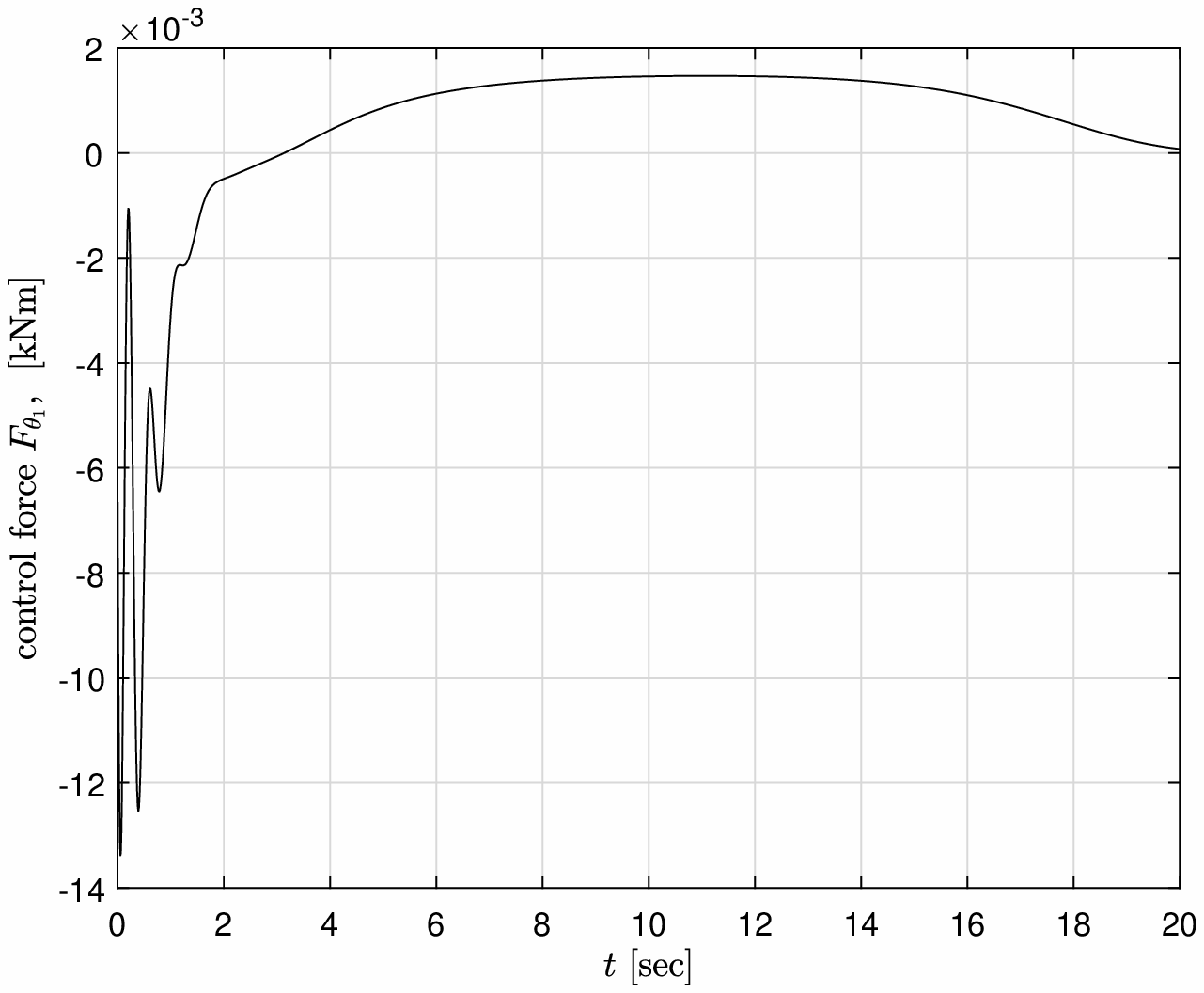,width=5.cm,clip=}
         }
   }
\caption{The corresponding control forces to Fig.~\ref{solutions_xi_without_Ft2}; $F_z(t)\rightarrow 0,$ $F_{\theta_1}(t)\rightarrow 0$ and $F_{l_1}(t)\rightarrow -g(m_h+m_p)$ for $t\rightarrow \infty$ as follows from Remark~\ref{forces_asymptotics}.}
\label{control_forces_without_Ft2}
\end{figure}

\begin{figure} 
   \centerline{
    \hbox{
     \psfig{file=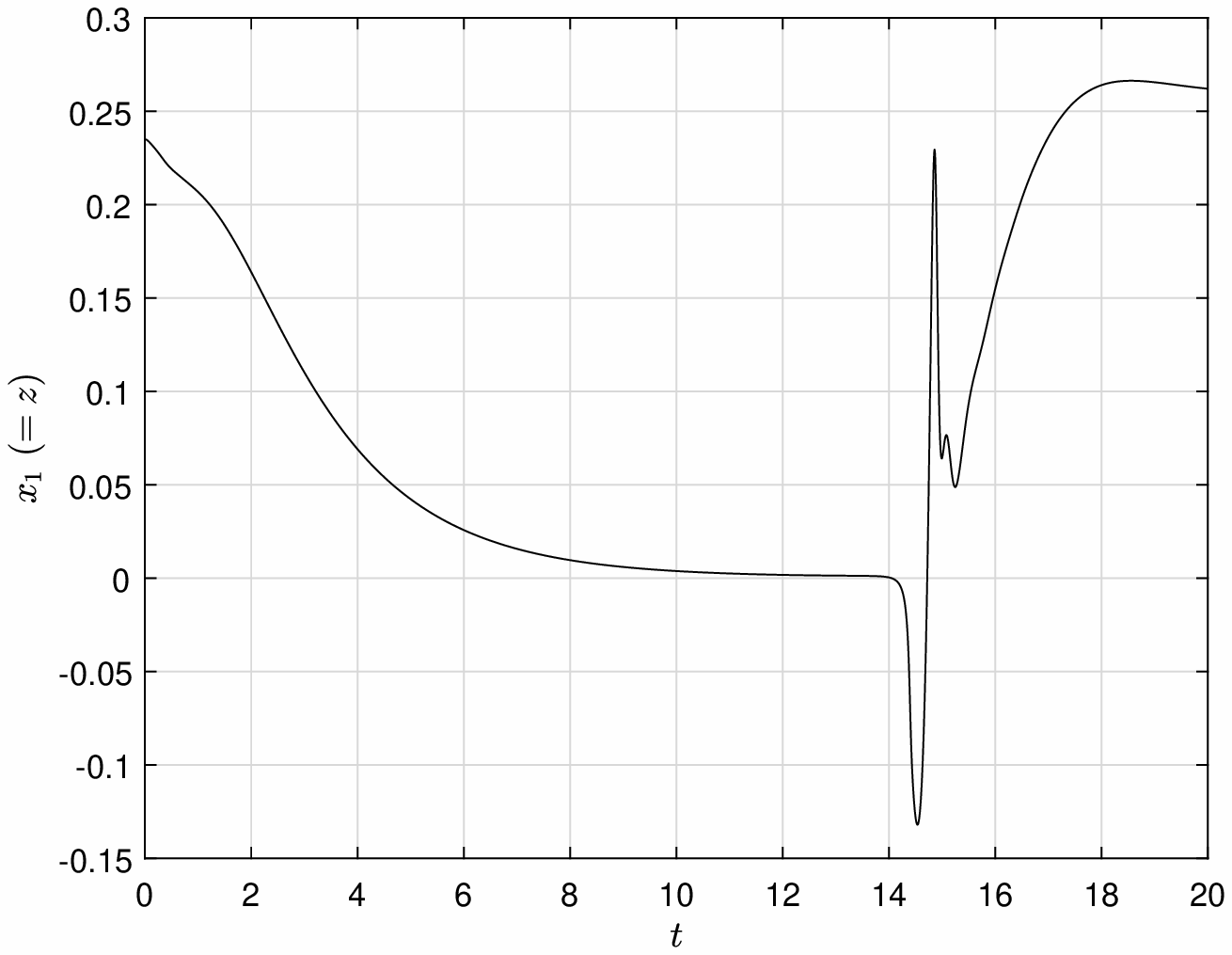,width=5.cm, clip=}
     \hspace{1.cm}
     \psfig{file=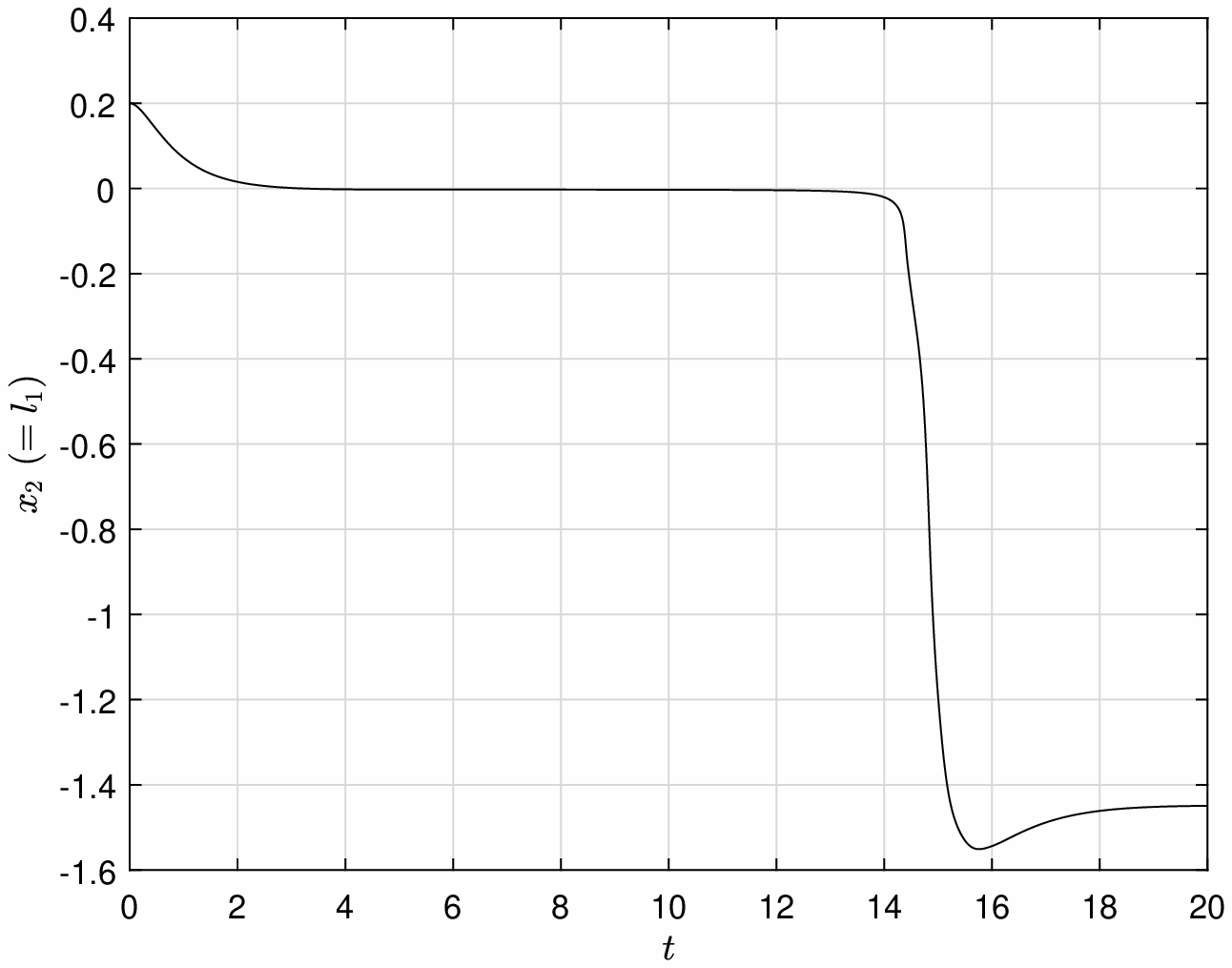,width=5cm,clip=}
    }
   }
   \vspace{0.5cm}
   \centerline{
    \hbox{
     \psfig{file=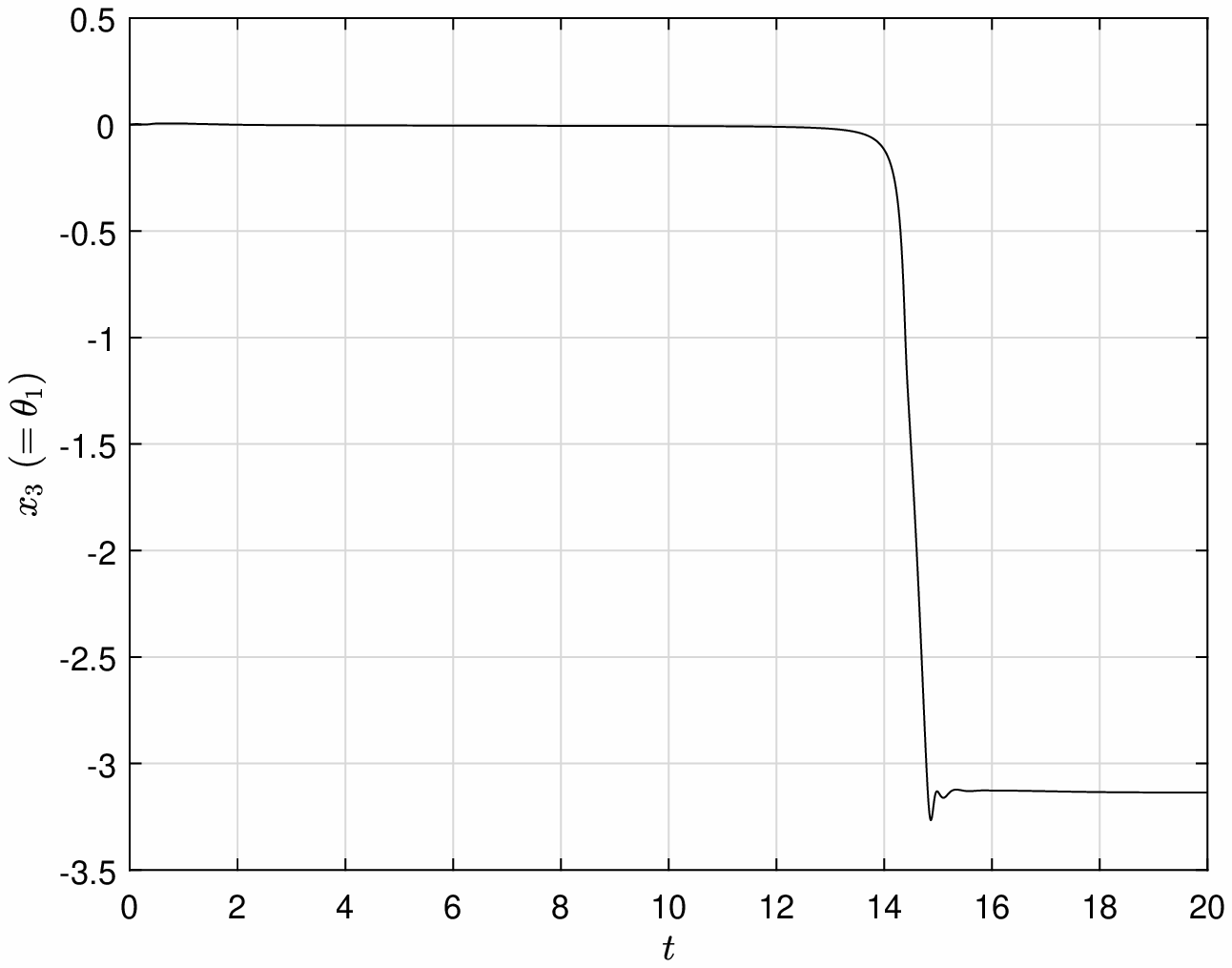,width=5.cm,clip=}
     \hspace{1.cm}
     \psfig{file=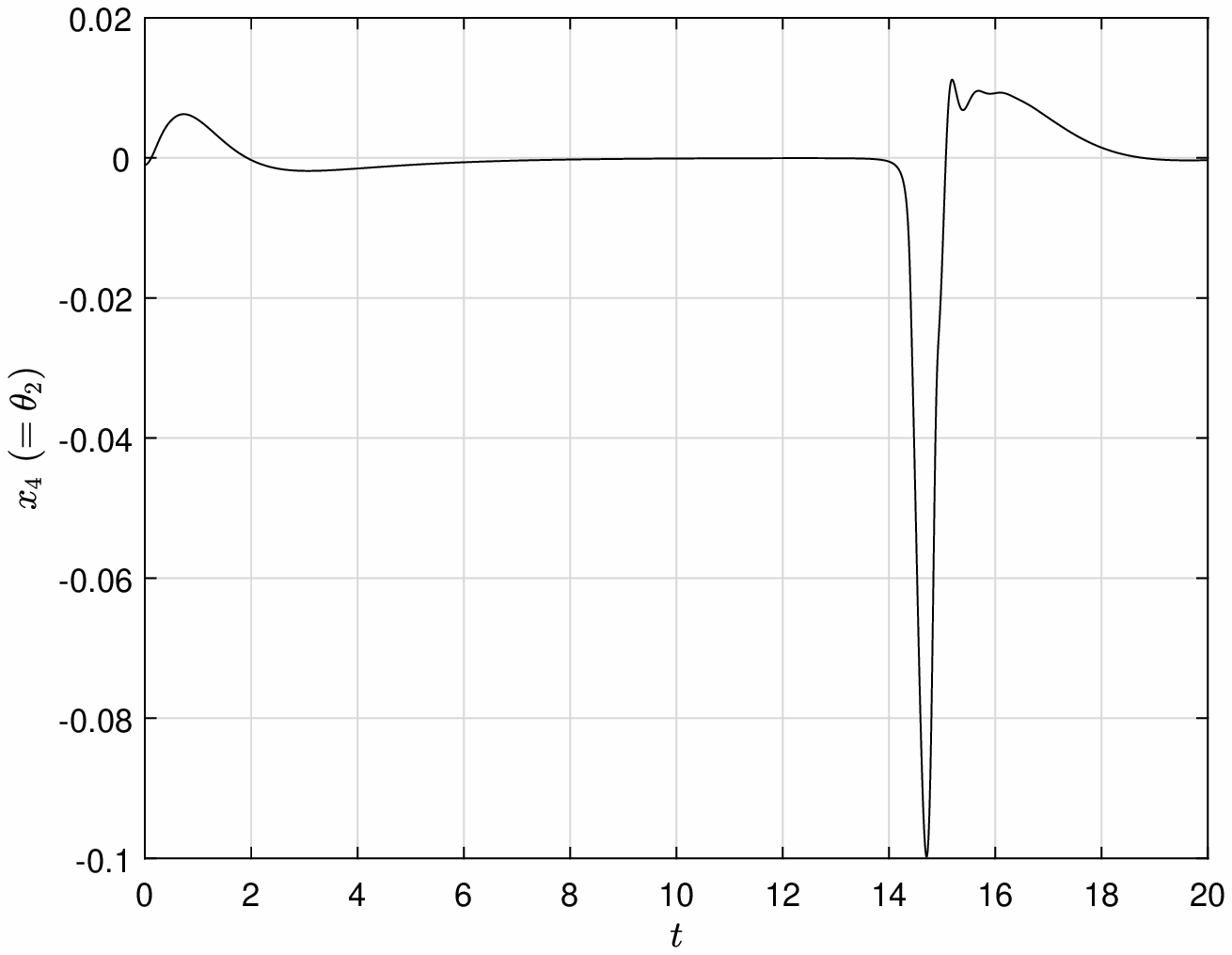,width=5.cm,clip=}
    }
   }
\caption{The solutions $ x_i,$ $i=1,\dots,4$ of the reference model {\bf without considering} the angular torque control $F_{\theta_2}$ with the initial state $(0.235\ 0.2\ 0\ -0.001\ 0\ 0\ 0\ 0)^T$ demonstrating the sensitivity of the system on the initial data change.}
\label{solutions_xi_without_Ft2_n}
\end{figure}

The Kalman's test of controllability of the pair $(A,B),$ $A=G_x(0,0),$ $B=G_u(0,0)$ gives
\[
\rank\left(B,AB,A^2B,\dots, A^7B\right)=8,
\]
which implies the controllability of the linear part of the control system $\dot x=G(x,u).$ For the simulation purposes, let us consider for the permissible eigenvalues of $A_{cl}=A-BK$ the set 
\[
p = [-0.5 -1 -1.5 -2 -2.5 -3 -3.5 -4].
\]
For the purpose of numerical simulation in Matlab the following data will be used: 
\begin{itemize}
\item[] $M=0.2\ [\times10^3\,\mathrm{kg}],\ $ $m_p=10\ [\times10^3\,\mathrm{kg}],\ $ $m_h=0.1\ [\times10^3\,\mathrm{kg}],\ $ $I_p=4\ [\times10^3\,\mathrm{kg\, m^2}],\ $ $I_h = 0.05\ [\times10^3\,\mathrm{kg\, m^2}],\ $  $g = 9.81\ [\mathrm{m\, s^{-2}}],$ $ l_2=2\ [\mathrm{m}],$
\end{itemize}
and the Matlab output of pole placement command {\tt K = place(A,B,p)} gives the gain matrix $K:$ 
\[
K=\left(
\begin{array}{rrrrrrrr} 
 2.0574  &  0.1052  &  0.1224 &  42.3471  &  6.7649 &   0.0526  &  0.0488 &  -7.7172 \\
 64.5917 & 214.6134  & -93.7833 & -985.9842 & 173.1363 & 236.2118  &-36.9947 & 180.4709 \\
 0.0840  & -0.0704 &   0.3034 &  -0.9495 &   0.2016 &  -0.0283  &  0.2554  &  0.2631
\end{array}
\right)
\]
and so the state feedback controller $(u_1,u_2,u_3)^T=-K(x_1,\dots,x_8)^T$ defines the control forces for the reference system and for the feed-forward controller. The routine {\tt place} uses the algorithm of \cite{KaNiDo} which, for multi-input systems, optimizes the choice of eigenvectors for a robust solution. The corresponding solutions $ x_i,$ $i=1,\dots,8$ are depicted on the Fig.~\ref{solutions_xi_without_Ft2} and the control forces $F_z,$ $F_{l_1}$ and $F_{\theta_1}$ on the Fig.~\ref{control_forces_without_Ft2}. In this context, the Fig.~\ref{solutions_xi_without_Ft2_n} illustrates the sensitivity of the system on the initial data change, which can be suppressed by using the control force $F_{\theta_2},$ analyzed in the following section.

\clearpage 
\newpage

\section{Crane system with considering the angular torque control $F_{\theta_2}$}\label{with_Ft2}
In this section we will consider the D-P-T crane system described by Lagrange equations (\ref{lagrange_z}), (\ref{lagrange_l_1}), (\ref{lagrange_theta_1}) and (\ref{lagrange_theta_2}), where the the swing angle of payload is controlled by force $F_{\theta_2}.$

Now, solving the equations  with regard to the variable $\ddot z,$ $\ddot l_1,$ $\ddot\theta_1,$ $\ddot\theta_2,$
and substituting $x_1,$ $x_2,$ $x_3,$ $x_4,$ $x_5,$ $x_6,$ $x_7$ and $x_8$ instead of $z,$ $l_1,$ $\theta_1,$ $\theta_2,$ $\dot z,$ $\dot l_1,$ $\dot\theta_1$ and $\dot\theta_2,$ respectively,  we obtain the functions $\hat G_i,$ $i=5,6,7,8:$
\[
\hat G_5=
 -\frac{1}{\Delta}\times\Bigg[{F_{\theta_2}}I_{h}l_{2}{m^2_{p}}\cos\left(x_{4}\right)-2{F_z}I_{p}{m^2_{h}}{x^2_{2}}-2{F_z}I_{p}{m^2_{p}}{x^2_{2}}
\]
\[
-2{F_z}I_{h}I_{p}m_{h}-2{F_z}I_{h}I_{p}m_{p}-{F_z}I_{h}{l^2_{2}}{m^2_{p}}+2{F_{\theta_1}}I_{p}{m^2_{h}}x_{2}\cos\left(x_{3}\right)
\]
\[
+2{F_{\theta_1}}I_{p}{m^2_{p}}x_{2}\cos\left(x_{3}\right)
+I_{h}I_{p}g{m^2_{h}}\sin\left(2x_{3}\right)+I_{h}I_{p}g{m^2_{p}}\sin\left(2x_{3}\right)+{F_{l_1}}I_{h}{l^2_{2}}{m^2_{p}}\sin\left(x_{3}\right)
\]
\[
+2{F_{l_1}}I_{p}{m^2_{h}}{x^2_{2}}\sin\left(x_{3}\right)+2{F_{l_1}}I_{p}{m^2_{p}}{x^2_{2}}\sin\left(x_{3}\right)-2{F_z}{l^2_{2}}m_{h}{m^2_{p}}{x^2_{2}}
\]
\[
-2{F_z}{l^2_{2}}{m^2_{h}}m_{p}{x^2_{2}}+2{F_{l_1}}I_{h}I_{p}m_{h}\sin\left(x_{3}\right)+2{F_{l_1}}I_{h}I_{p}m_{p}\sin\left(x_{3}\right)-2{F_z}I_{h}{l^2_{2}}m_{h}m_{p}
\]
\[
+{F_{\theta_2}}I_{h}l_{2}{m^2_{p}}\cos\left(2x_{3}-x_{4}\right)-4{F_z}I_{p}m_{h}m_{p}{x^2_{2}}-{F_{l_1}}I_{h}{l^2_{2}}{m^2_{p}}\sin\left(x_{3}-2x_{4}\right)
\]
\[
-{F_z}I_{h}{l^2_{2}}{m^2_{p}}\cos\left(2x_{3}-2x_{4}\right)+2{F_{\theta_1}}{l^2_{2}}m_{h}{m^2_{p}}x_{2}\cos\left(x_{3}\right)
\]
\[
+2{F_{\theta_1}}{l^2_{2}}{m^2_{h}}m_{p}x_{2}\cos\left(x_{3}\right)+I_{h}{l^3_{2}}m_{h}{m^2_{p}}{x^2_{8}}\sin\left(2x_{3}-x_{4}\right)+{F_{\theta_2}}I_{h}l_{2}m_{h}m_{p}\cos\left(x_{4}\right)
\]
\[
+4{F_{\theta_1}}I_{p}m_{h}m_{p}x_{2}\cos\left(x_{3}\right)+I_{h}g{l^2_{2}}m_{h}{m^2_{p}}\sin\left(2x_{3}\right)+I_{h}g{l^2_{2}}{m^2_{h}}m_{p}\sin\left(2x_{3}\right)
\]
\[
+2{F_{l_1}}{l^2_{2}}m_{h}{m^2_{p}}{x^2_{2}}\sin\left(x_{3}\right)+2{F_{l_1}}{l^2_{2}}{m^2_{h}}m_{p}{x^2_{2}}\sin\left(x_{3}\right)
\]
\[
-I_{h}{l^3_{2}}m_{h}{m^2_{p}}{x^2_{8}}\sin\left(x_{4}\right)+2I_{h}I_{p}gm_{h}m_{p}\sin\left(2x_{3}\right)+2{F_{l_1}}I_{h}{l^2_{2}}m_{h}m_{p}\sin\left(x_{3}\right)
\]
\[
+I_{h}I_{p}l_{2}{m^2_{p}}{x^2_{8}}\sin\left(2x_{3}-x_{4}\right)+4I_{h}I_{p}{m^2_{h}}x_{6}x_{7}\cos\left(x_{3}\right)+4I_{h}I_{p}{m^2_{p}}x_{6}x_{7}\cos\left(x_{3}\right)
\]
\[
+4{F_{l_1}}I_{p}m_{h}m_{p}{x^2_{2}}\sin\left(x_{3}\right)+{F_{\theta_2}}I_{h}l_{2}m_{h}m_{p}\cos\left(2x_{3}-x_{4}\right)-I_{h}I_{p}l_{2}{m^2_{p}}{x^2_{8}}\sin\left(x_{4}\right)
\]
\[
+4I_{h}{l^2_{2}}m_{h}{m^2_{p}}x_{6}x_{7}\cos\left(x_{3}\right)+4I_{h}{l^2_{2}}{m^2_{h}}m_{p}x_{6}x_{7}\cos\left(x_{3}\right)+I_{h}I_{p}l_{2}m_{h}m_{p}{x^2_{8}}\sin\left(2x_{3}-x_{4}\right)
\]
\[
+8I_{h}I_{p}m_{h}m_{p}x_{6}x_{7}\cos\left(x_{3}\right)-I_{h}I_{p}l_{2}m_{h}m_{p}{x^2_{8}}\sin\left(x_{4}\right)\Bigg],
\]

\[
\hat G_6=
\frac{1}{\Delta}\times\Bigg[2{F_{l_1}}I_{h}I_{p}M+{F_{l_1}}I_{h}{l^2_{2}}{m^2_{p}}+{F_{l_1}}I_{p}{m^2_{h}}{x^2_{2}}+{F_{l_1}}I_{p}{m^2_{p}}{x^2_{2}}
\]
\[
+2{F_{l_1}}I_{h}I_{p}m_{h}+2{F_{l_1}}I_{h}I_{p}m_{p}+2I_{h}I_{p}g{m^2_{h}}\cos\left(x_{3}\right)+2I_{h}I_{p}g{m^2_{p}}\cos\left(x_{3}\right)
\]
\[
+I_{h}I_{p}{m^2_{h}}x_{2}{x^2_{7}}+I_{h}I_{p}{m^2_{p}}x_{2}{x^2_{7}}-{F_{\theta_2}}I_{h}l_{2}{m^2_{p}}\sin\left(x_{3}-x_{4}\right)
\]
\[
+{F_{l_1}}M{l^2_{2}}{m^2_{p}}{x^2_{2}}+2I_{p}M{m^2_{h}}{x_{2}}^3{x^2_{7}}+2I_{p}M{m^2_{p}}{x_{2}}^3{x^2_{7}}-{F_z}I_{h}{l^2_{2}}{m^2_{p}}\sin\left(x_{3}\right)
\]
\[
+2{F_{l_1}}I_{h}M{l^2_{2}}m_{p}+{F_{\theta_1}}I_{p}{m^2_{h}}x_{2}\sin\left(2x_{3}\right)+{F_{\theta_1}}I_{p}{m^2_{p}}x_{2}\sin\left(2x_{3}\right)
\]
\[
-2{F_z}I_{p}{m^2_{h}}{x^2_{2}}\sin\left(x_{3}\right)-2{F_z}I_{p}{m^2_{p}}{x^2_{2}}\sin\left(x_{3}\right)+2{F_{l_1}}I_{p}Mm_{h}{x^2_{2}}+2{F_{l_1}}I_{p}Mm_{p}{x^2_{2}}
\]
\[
+{F_{l_1}}{l^2_{2}}m_{h}{m^2_{p}}{x^2_{2}}+{F_{l_1}}{l^2_{2}}{m^2_{h}}m_{p}{x^2_{2}}-2{F_z}I_{h}I_{p}m_{h}\sin\left(x_{3}\right)
\]
\[
-2{F_z}I_{h}I_{p}m_{p}\sin\left(x_{3}\right)+2{F_{l_1}}I_{h}{l^2_{2}}m_{h}m_{p}+2{F_{l_1}}I_{p}m_{h}m_{p}{x^2_{2}}
\]
\[
+{F_z}I_{h}{l^2_{2}}{m^2_{p}}\sin\left(x_{3}-2x_{4}\right)-{F_{l_1}}I_{h}{l^2_{2}}{m^2_{p}}\cos\left(2x_{4}\right)-{F_{l_1}}I_{p}{m^2_{h}}{x^2_{2}}\cos\left(2x_{3}\right)
\]
\[
-{F_{l_1}}I_{p}{m^2_{p}}{x^2_{2}}\cos\left(2x_{3}\right)+{F_{\theta_2}}I_{h}l_{2}{m^2_{p}}\sin\left(x_{3}+x_{4}\right)+I_{h}Mg{l^2_{2}}{m^2_{p}}\cos\left(x_{3}-2x_{4}\right)
\]
\[
+4I_{h}I_{p}gm_{h}m_{p}\cos\left(x_{3}\right)+I_{h}I_{p}l_{2}{m^2_{p}}{x^2_{8}}\cos\left(x_{3}-x_{4}\right)+2I_{h}I_{p}m_{h}m_{p}x_{2}{x^2_{7}}
\]
\[
-2{F_{\theta_2}}Ml_{2}{m^2_{p}}{x^2_{2}}\sin\left(x_{3}-x_{4}\right)-2{F_{\theta_2}}I_{h}Ml_{2}m_{p}\sin\left(x_{3}-x_{4}\right)+I_{h}{l^3_{2}}m_{h}{m^2_{p}}{x^2_{8}}\cos\left(x_{3}+x_{4}\right)
\]
\[
+I_{h}I_{p}{m^2_{h}}x_{2}{x^2_{7}}\cos\left(2x_{3}\right)+I_{h}I_{p}{m^2_{p}}x_{2}{x^2_{7}}\cos\left(2x_{3}\right)-{F_{\theta_2}}I_{h}l_{2}m_{h}m_{p}\sin\left(x_{3}-x_{4}\right)
\]
\[
+2M{l^2_{2}}m_{h}{m^2_{p}}{x_{2}}^3{x^2_{7}}+2M{l^2_{2}}{m^2_{h}}m_{p}{x_{2}}^3{x^2_{7}}+{F_{\theta_1}}{l^2_{2}}m_{h}{m^2_{p}}x_{2}\sin\left(2x_{3}\right)
\]
\[
+{F_{\theta_1}}{l^2_{2}}{m^2_{h}}m_{p}x_{2}\sin\left(2x_{3}\right)-2{F_z}{l^2_{2}}m_{h}{m^2_{p}}{x^2_{2}}\sin\left(x_{3}\right)
\]
\[
-2{F_z}{l^2_{2}}{m^2_{h}}m_{p}{x^2_{2}}\sin\left(x_{3}\right)+2{F_{l_1}}M{l^2_{2}}m_{h}m_{p}{x^2_{2}}+2I_{h}M{l^3_{2}}{m^2_{p}}{x^2_{8}}\cos\left(x_{3}-x_{4}\right)
\]
\[
+4I_{p}Mm_{h}m_{p}{x_{2}}^3{x^2_{7}}-2{F_z}I_{h}{l^2_{2}}m_{h}m_{p}\sin\left(x_{3}\right)+{F_{\theta_1}}M{l^2_{2}}{m^2_{p}}x_{2}\sin\left(2x_{3}-2x_{4}\right)
\]
\[
+2{F_{\theta_1}}I_{p}m_{h}m_{p}x_{2}\sin\left(2x_{3}\right)-4{F_z}I_{p}m_{h}m_{p}{x^2_{2}}\sin\left(x_{3}\right)+I_{h}{l^3_{2}}m_{h}{m^2_{p}}{x^2_{8}}\cos\left(x_{3}-x_{4}\right)
\]
\[
+I_{h}I_{p}l_{2}{m^2_{p}}{x^2_{8}}\cos\left(x_{3}+x_{4}\right)-{F_{l_1}}{l^2_{2}}m_{h}{m^2_{p}}{x^2_{2}}\cos\left(2x_{3}\right)-{F_{l_1}}{l^2_{2}}{m^2_{h}}m_{p}{x^2_{2}}\cos\left(2x_{3}\right)
\]
\[
+I_{h}Mg{l^2_{2}}{m^2_{p}}\cos\left(x_{3}\right)-{F_{l_1}}M{l^2_{2}}{m^2_{p}}{x^2_{2}}\cos\left(2x_{3}-2x_{4}\right)-2{F_{l_1}}I_{p}m_{h}m_{p}{x^2_{2}}\cos\left(2x_{3}\right)
\]
\[
+2I_{p}Mg{m^2_{h}}{x^2_{2}}\cos\left(x_{3}\right)+2I_{p}Mg{m^2_{p}}{x^2_{2}}\cos\left(x_{3}\right)+I_{h}M{l^2_{2}}{m^2_{p}}x_{2}{x^2_{7}}+2I_{h}I_{p}Mgm_{h}\cos\left(x_{3}\right)
\]
\[
+2I_{h}I_{p}Mgm_{p}\cos\left(x_{3}\right)+2I_{h}g{l^2_{2}}m_{h}{m^2_{p}}\cos\left(x_{3}\right)+2I_{h}g{l^2_{2}}{m^2_{h}}m_{p}\cos\left(x_{3}\right)+{F_{\theta_2}}I_{h}l_{2}m_{h}m_{p}\sin\left(x_{3}+x_{4}\right)
\]
\[
+2I_{h}I_{p}{m^2_{h}}x_{6}x_{7}\sin\left(2x_{3}\right)+2I_{h}I_{p}{m^2_{p}}x_{6}x_{7}\sin\left(2x_{3}\right)+2I_{h}I_{p}Mm_{h}x_{2}{x^2_{7}}
\]
\[
+2I_{h}I_{p}Mm_{p}x_{2}{x^2_{7}}+I_{h}{l^2_{2}}m_{h}{m^2_{p}}x_{2}{x^2_{7}}+I_{h}{l^2_{2}}{m^2_{h}}m_{p}x_{2}{x^2_{7}}+2I_{h}M{l^2_{2}}{m^2_{p}}x_{6}x_{7}\sin\left(2x_{3}-2x_{4}\right)
\]
\[
+2I_{h}I_{p}Ml_{2}m_{p}{x^2_{8}}\cos\left(x_{3}-x_{4}\right)+4I_{h}I_{p}m_{h}m_{p}x_{6}x_{7}\sin\left(2x_{3}\right)+I_{h}I_{p}l_{2}m_{h}m_{p}{x^2_{8}}\cos\left(x_{3}-x_{4}\right)
\]
\[
-2{F_{\theta_2}}Ml_{2}m_{h}m_{p}{x^2_{2}}\sin\left(x_{3}-x_{4}\right)+I_{h}{l^2_{2}}m_{h}{m^2_{p}}x_{2}{x^2_{7}}\cos\left(2x_{3}\right)+I_{h}{l^2_{2}}{m^2_{h}}m_{p}x_{2}{x^2_{7}}\cos\left(2x_{3}\right)
\]
\[
+I_{h}M{l^2_{2}}{m^2_{p}}x_{2}{x^2_{7}}\cos\left(2x_{3}-2x_{4}\right)+2I_{h}I_{p}m_{h}m_{p}x_{2}{x^2_{7}}\cos\left(2x_{3}\right)
\]
\[
+2M{l^3_{2}}m_{h}{m^2_{p}}{x^2_{2}}{x^2_{8}}\cos\left(x_{3}-x_{4}\right)+I_{h}I_{p}l_{2}m_{h}m_{p}{x^2_{8}}\cos\left(x_{3}+x_{4}\right)
\]
\[
+2Mg{l^2_{2}}m_{h}{m^2_{p}}{x^2_{2}}\cos\left(x_{3}\right)+2Mg{l^2_{2}}{m^2_{h}}m_{p}{x^2_{2}}\cos\left(x_{3}\right)+2I_{h}Mg{l^2_{2}}m_{h}m_{p}\cos\left(x_{3}\right)
\]
\[
+2I_{h}{l^2_{2}}m_{h}{m^2_{p}}x_{6}x_{7}\sin\left(2x_{3}\right)+2I_{h}{l^2_{2}}{m^2_{h}}m_{p}x_{6}x_{7}\sin\left(2x_{3}\right)+4I_{p}Mgm_{h}m_{p}{x^2_{2}}\cos\left(x_{3}\right)
\]
\[
+2I_{h}M{l^2_{2}}m_{h}m_{p}x_{2}{x^2_{7}}+2I_{p}Ml_{2}{m^2_{p}}{x^2_{2}}{x^2_{8}}\cos\left(x_{3}-x_{4}\right)+2I_{p}Ml_{2}m_{h}m_{p}{x^2_{2}}{x^2_{8}}\cos\left(x_{3}-x_{4}\right)\Bigg],
\]

\[
\hat G_7=
\frac{1}{\Delta}\times\Bigg[{F_{\theta_1}}I_{p}{m^2_{h}}+{F_{\theta_1}}I_{p}{m^2_{p}}+{F_{\theta_1}}M{l^2_{2}}{m^2_{p}}+2{F_{\theta_1}}I_{p}Mm_{h}
\]
\[
+2{F_{\theta_1}}I_{p}Mm_{p}+{F_{\theta_1}}{l^2_{2}}m_{h}{m^2_{p}}+{F_{\theta_1}}{l^2_{2}}{m^2_{h}}m_{p}+2{F_{\theta_1}}I_{p}m_{h}m_{p}+{F_{\theta_1}}I_{p}{m^2_{h}}\cos\left(2x_{3}\right)
\]
\[
+{F_{\theta_1}}I_{p}{m^2_{p}}\cos\left(2x_{3}\right)+2{F_{\theta_1}}I_{p}m_{h}m_{p}\cos\left(2x_{3}\right)-2{F_z}I_{p}{m^2_{h}}x_{2}\cos\left(x_{3}\right)
\]
\[
-2{F_z}I_{p}{m^2_{p}}x_{2}\cos\left(x_{3}\right)+{F_{l_1}}I_{p}{m^2_{h}}x_{2}\sin\left(2x_{3}\right)+{F_{l_1}}I_{p}{m^2_{p}}x_{2}\sin\left(2x_{3}\right)
\]
\[
+2{F_{\theta_1}}M{l^2_{2}}m_{h}m_{p}+{F_{\theta_1}}{l^2_{2}}m_{h}{m^2_{p}}\cos\left(2x_{3}\right)+{F_{\theta_1}}{l^2_{2}}{m^2_{h}}m_{p}\cos\left(2x_{3}\right)
\]
\[
+{F_{\theta_1}}M{l^2_{2}}{m^2_{p}}\cos\left(2x_{3}-2x_{4}\right)-2{F_z}{l^2_{2}}m_{h}{m^2_{p}}x_{2}\cos\left(x_{3}\right)
\]
\[
-2{F_z}{l^2_{2}}{m^2_{h}}m_{p}x_{2}\cos\left(x_{3}\right)-4{F_z}I_{p}m_{h}m_{p}x_{2}\cos\left(x_{3}\right)-4I_{p}M{m^2_{h}}x_{2}x_{6}x_{7}-4I_{p}M{m^2_{p}}x_{2}x_{6}x_{7}
\]
\[
+{F_{l_1}}{l^2_{2}}m_{h}{m^2_{p}}x_{2}\sin\left(2x_{3}\right)+{F_{l_1}}{l^2_{2}}{m^2_{h}}m_{p}x_{2}\sin\left(2x_{3}\right)+{F_{l_1}}M{l^2_{2}}{m^2_{p}}x_{2}\sin\left(2x_{3}-2x_{4}\right)
\]
\[
+2{F_{l_1}}I_{p}m_{h}m_{p}x_{2}\sin\left(2x_{3}\right)-2I_{p}Mg{m^2_{h}}x_{2}\sin\left(x_{3}\right)-2I_{p}Mg{m^2_{p}}x_{2}\sin\left(x_{3}\right)
\]
\[
-2{F_{\theta_2}}Ml_{2}{m^2_{p}}x_{2}\cos\left(x_{3}-x_{4}\right)-4M{l^2_{2}}m_{h}{m^2_{p}}x_{2}x_{6}x_{7}-4M{l^2_{2}}{m^2_{h}}m_{p}x_{2}x_{6}x_{7}
\]
\[
-2M{l^3_{2}}m_{h}{m^2_{p}}x_{2}{x^2_{8}}\sin\left(x_{3}-x_{4}\right)-8I_{p}Mm_{h}m_{p}x_{2}x_{6}x_{7}-2Mg{l^2_{2}}m_{h}{m^2_{p}}x_{2}\sin\left(x_{3}\right)
\]
\[
-2Mg{l^2_{2}}{m^2_{h}}m_{p}x_{2}\sin\left(x_{3}\right)-4I_{p}Mgm_{h}m_{p}x_{2}\sin\left(x_{3}\right)-2I_{p}Ml_{2}{m^2_{p}}x_{2}{x^2_{8}}\sin\left(x_{3}-x_{4}\right)
\]
\[
-2{F_{\theta_2}}Ml_{2}m_{h}m_{p}x_{2}\cos\left(x_{3}-x_{4}\right)-2I_{p}Ml_{2}m_{h}m_{p}x_{2}{x^2_{8}}\sin\left(x_{3}-x_{4}\right)\Bigg]
\]
and
\[
\hat G_8=
 -\frac{1}{\Delta}\times\Bigg[\sin\left(x_{3}-x_{4}\right)\bigg({F_{l_1}}I_{h}l_{2}{m^2_{p}}+2{F_{l_1}}Ml_{2}{m^2_{p}}{x^2_{2}}+2{F_{l_1}}I_{h}Ml_{2}m_{p}
\]
\[
+{F_{l_1}}I_{h}l_{2}m_{h}m_{p}+2{F_{l_1}}Ml_{2}m_{h}m_{p}{x^2_{2}}\bigg)-{F_{\theta_2}}I_{h}{m^2_{h}}-{F_{\theta_2}}I_{h}{m^2_{p}}-2{F_{\theta_2}}M{m^2_{h}}{x^2_{2}}
\]
\[
-2{F_{\theta_2}}M{m^2_{p}}{x^2_{2}}-2{F_{\theta_2}}I_{h}Mm_{h}-2{F_{\theta_2}}I_{h}Mm_{p}-2{F_{\theta_2}}I_{h}m_{h}m_{p}-{F_{\theta_2}}I_{h}{m^2_{h}}\cos\left(2x_{3}\right)
\]
\[
-{F_{\theta_2}}I_{h}{m^2_{p}}\cos\left(2x_{3}\right)+{F_z}I_{h}l_{2}{m^2_{p}}\cos\left(x_{4}\right)-2{F_{\theta_2}}I_{h}m_{h}m_{p}\cos\left(2x_{3}\right)
\]
\[
+{F_z}I_{h}l_{2}{m^2_{p}}\cos\left(2x_{3}-x_{4}\right)-4{F_{\theta_2}}Mm_{h}m_{p}{x^2_{2}}-{F_{l_1}}I_{h}l_{2}{m^2_{p}}\sin\left(x_{3}+x_{4}\right)
\]
\[
+{F_z}I_{h}l_{2}m_{h}m_{p}\cos\left(x_{4}\right)+I_{h}Mgl_{2}{m^2_{p}}\sin\left(2x_{3}-x_{4}\right)+I_{h}Mgl_{2}{m^2_{p}}\sin\left(x_{4}\right)
\]
\[
+{F_z}I_{h}l_{2}m_{h}m_{p}\cos\left(2x_{3}-x_{4}\right)+2{F_{\theta_1}}Ml_{2}{m^2_{p}}x_{2}\cos\left(x_{3}-x_{4}\right)+I_{h}M{l^2_{2}}{m^2_{p}}{x^2_{8}}\sin\left(2x_{3}-2x_{4}\right)
\]
\[
-{F_{l_1}}I_{h}l_{2}m_{h}m_{p}\sin\left(x_{3}+x_{4}\right)+I_{h}Mgl_{2}m_{h}m_{p}\sin\left(2x_{3}-x_{4}\right)+4I_{h}Ml_{2}{m^2_{p}}x_{6}x_{7}\cos\left(x_{3}-x_{4}\right)
\]
\[
+I_{h}Mgl_{2}m_{h}m_{p}\sin\left(x_{4}\right)+2{F_{\theta_1}}Ml_{2}m_{h}m_{p}x_{2}\cos\left(x_{3}-x_{4}\right)+4I_{h}Ml_{2}m_{h}m_{p}x_{6}x_{7}\cos\left(x_{3}-x_{4}\right)\Bigg].
\]

\begin{figure} 
   \centerline{
    \hbox{
     \psfig{file=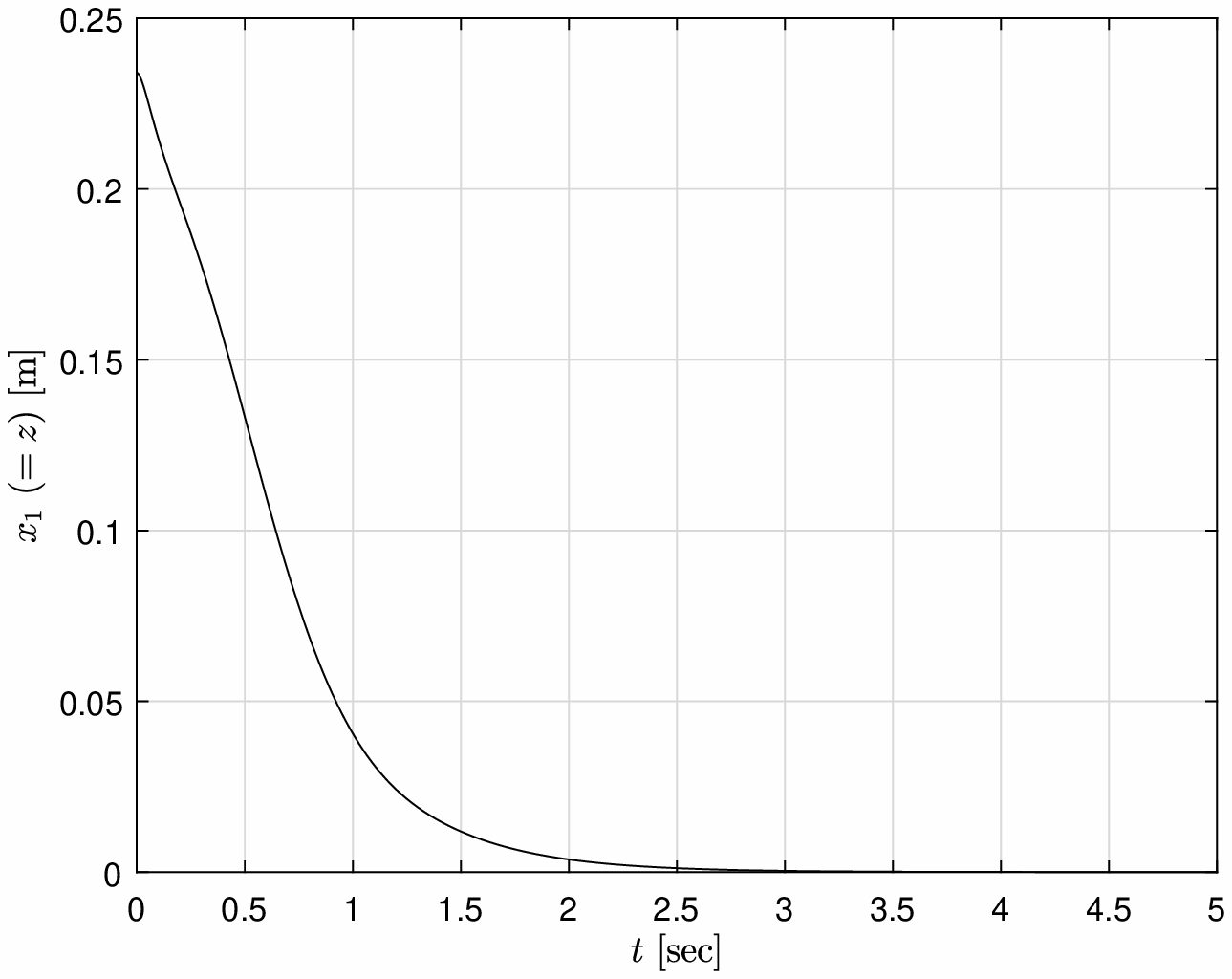,width=5.cm, clip=}
     \hspace{1.cm}
     \psfig{file=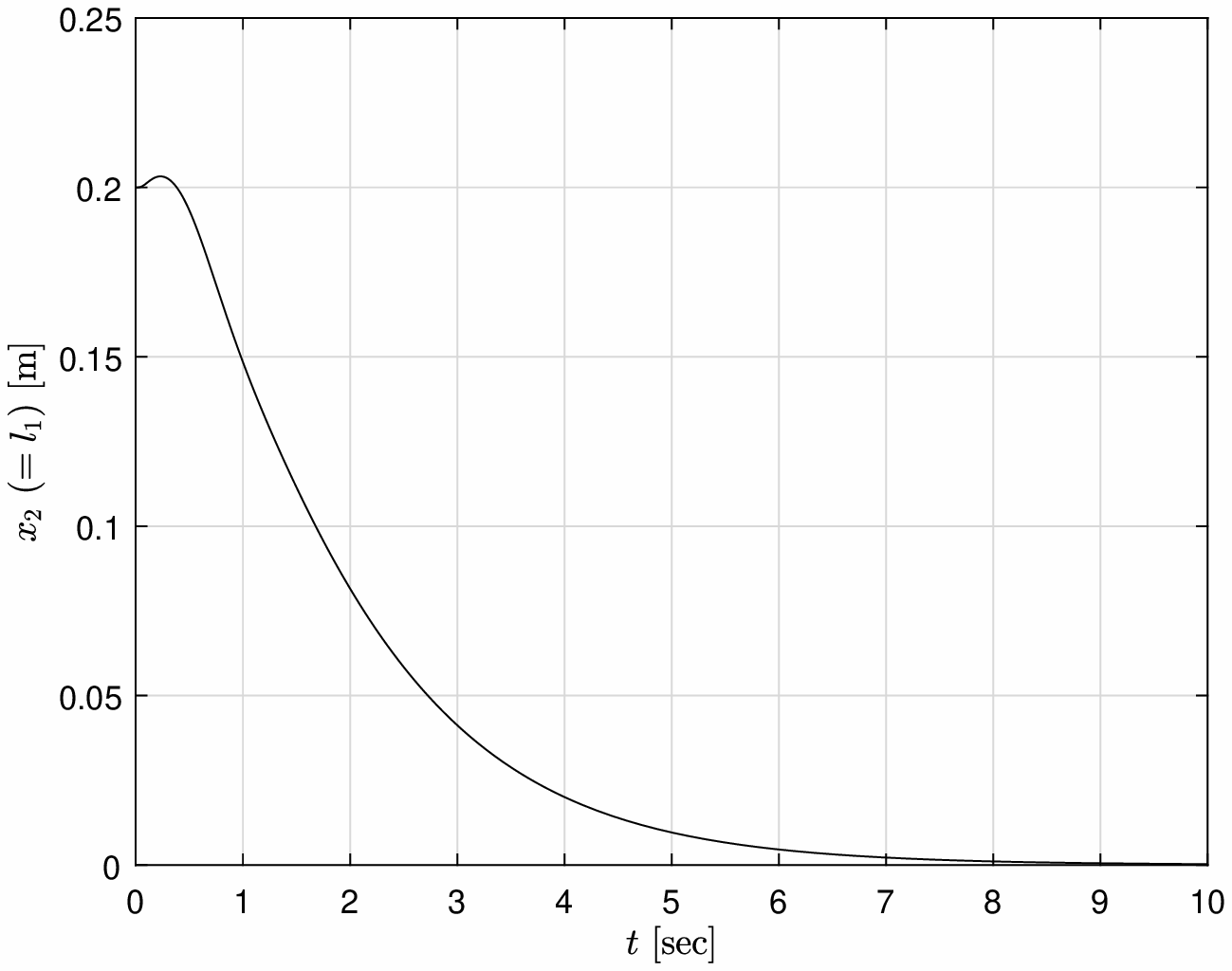,width=5cm,clip=}
    }
   }
   \vspace{0.5cm}
   \centerline{
    \hbox{
     \psfig{file=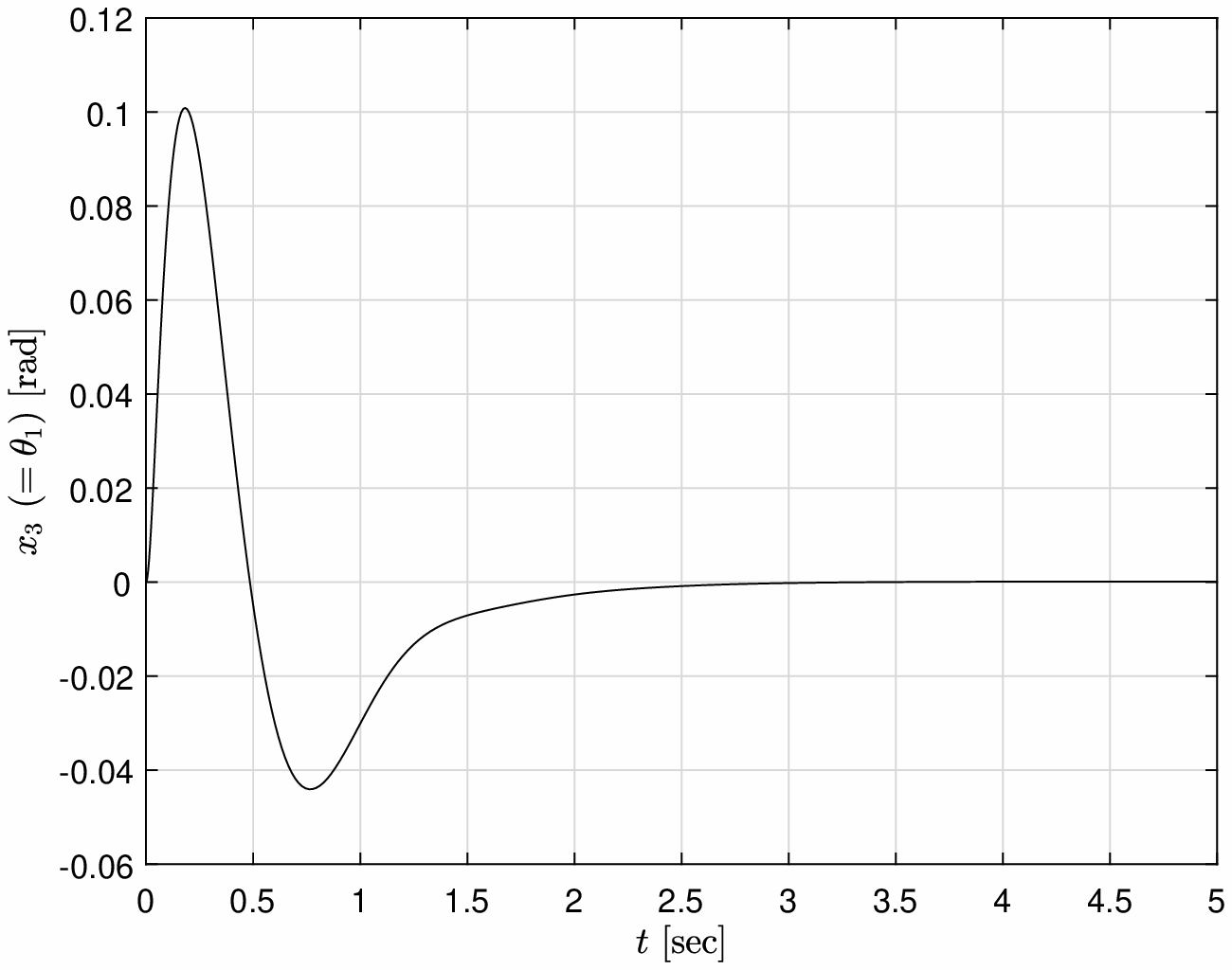,width=5.cm,clip=}
     \hspace{1.cm}
     \psfig{file=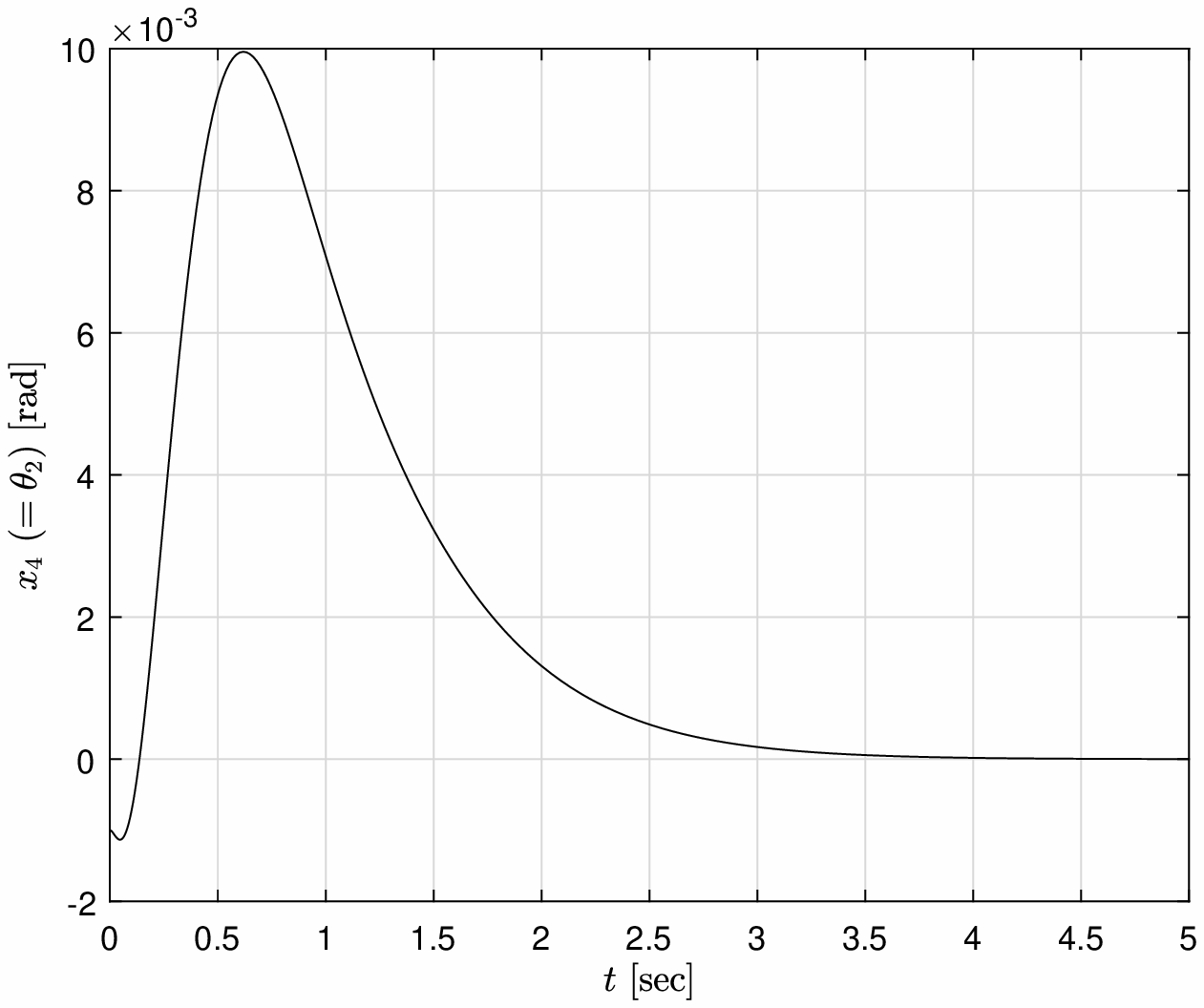,width=5.cm,clip=}
    }
   }
   \vspace{0.5cm}
   \centerline{
    \hbox{
     \psfig{file=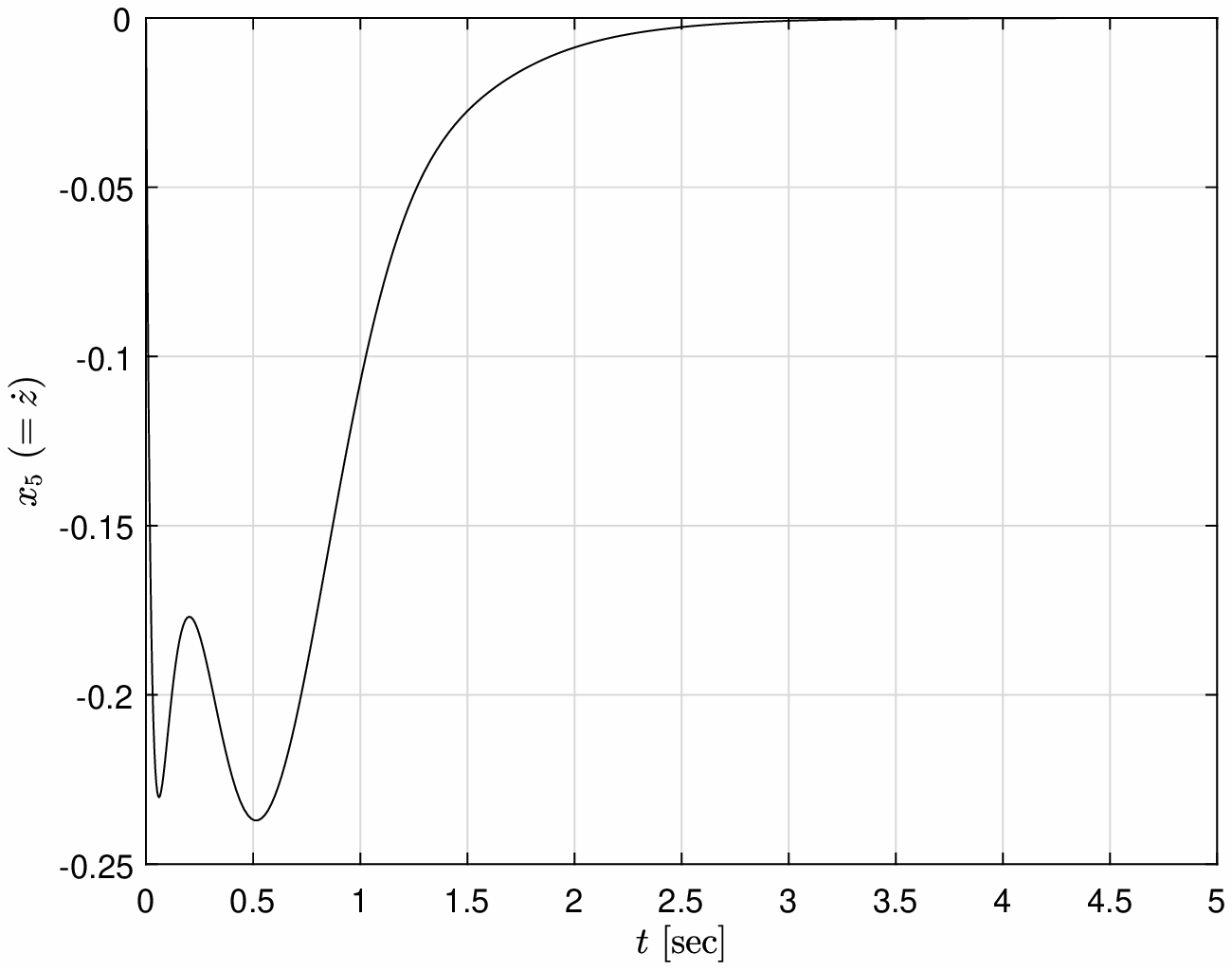,width=5.cm,clip=}
     \hspace{1.cm}
     \psfig{file=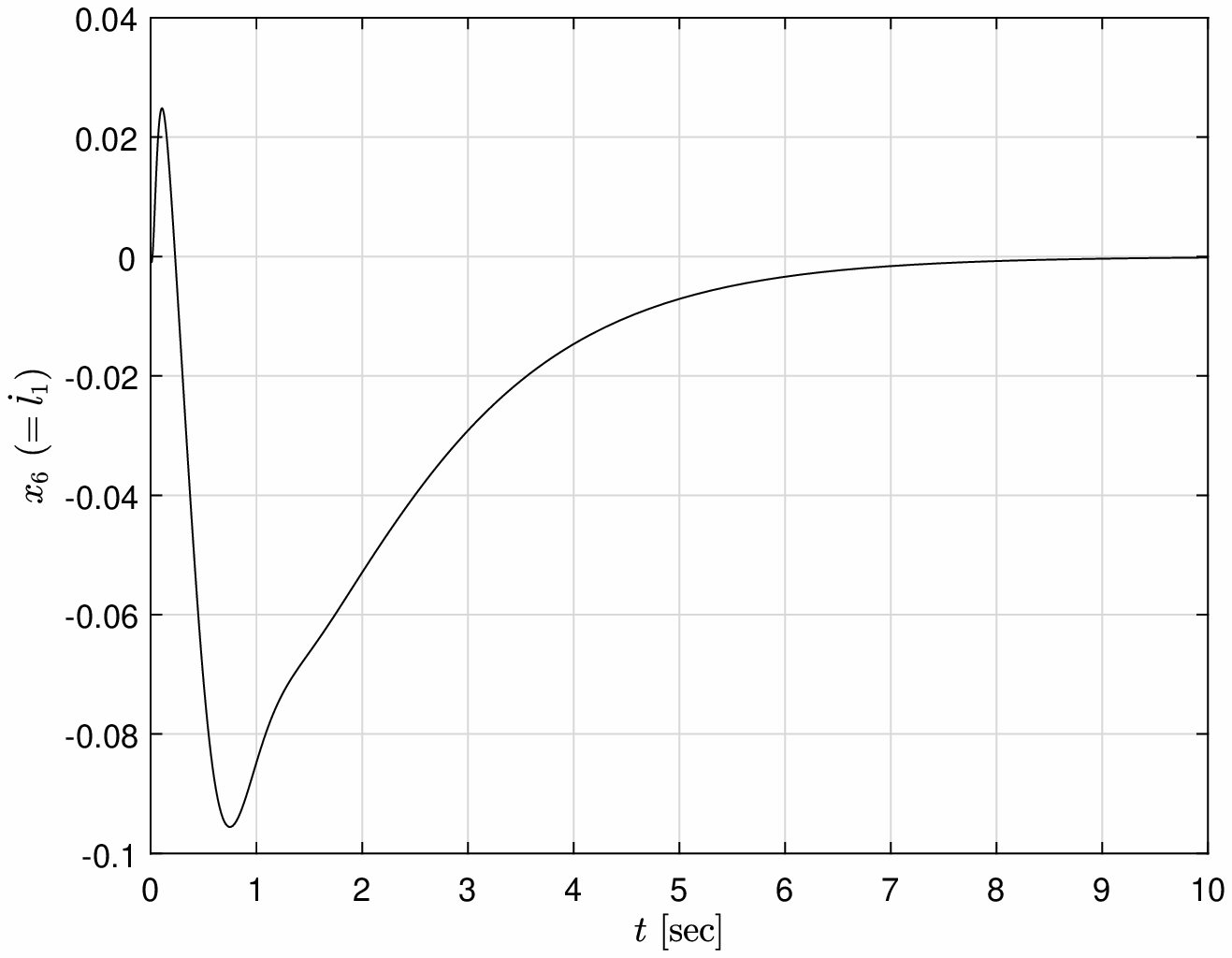,width=5.cm,clip=}
         }
   }
   \vspace{0.5cm}
   \centerline{
    \hbox{
     \psfig{file=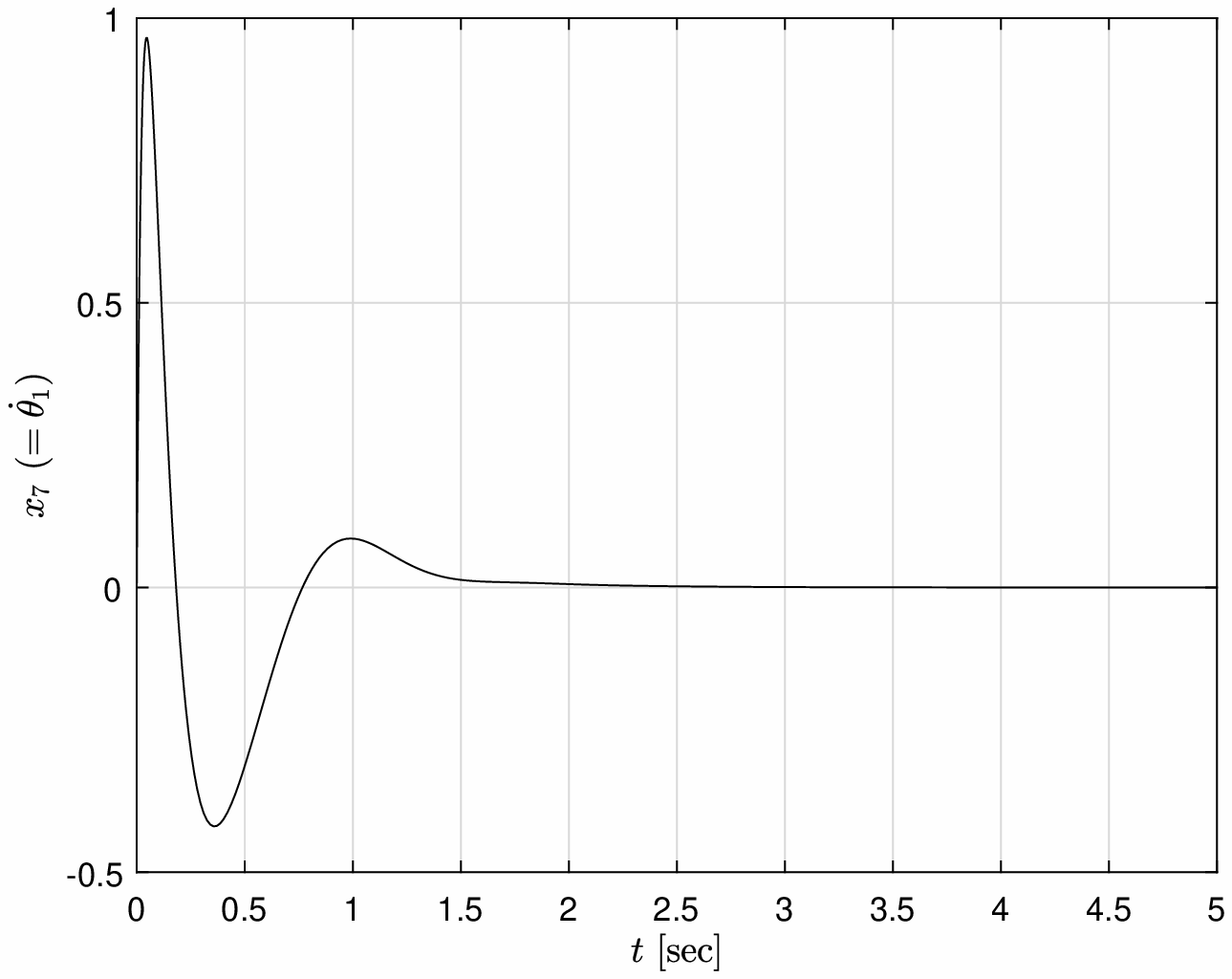,width=5.cm,clip=}
     \hspace{1.cm}
     \psfig{file=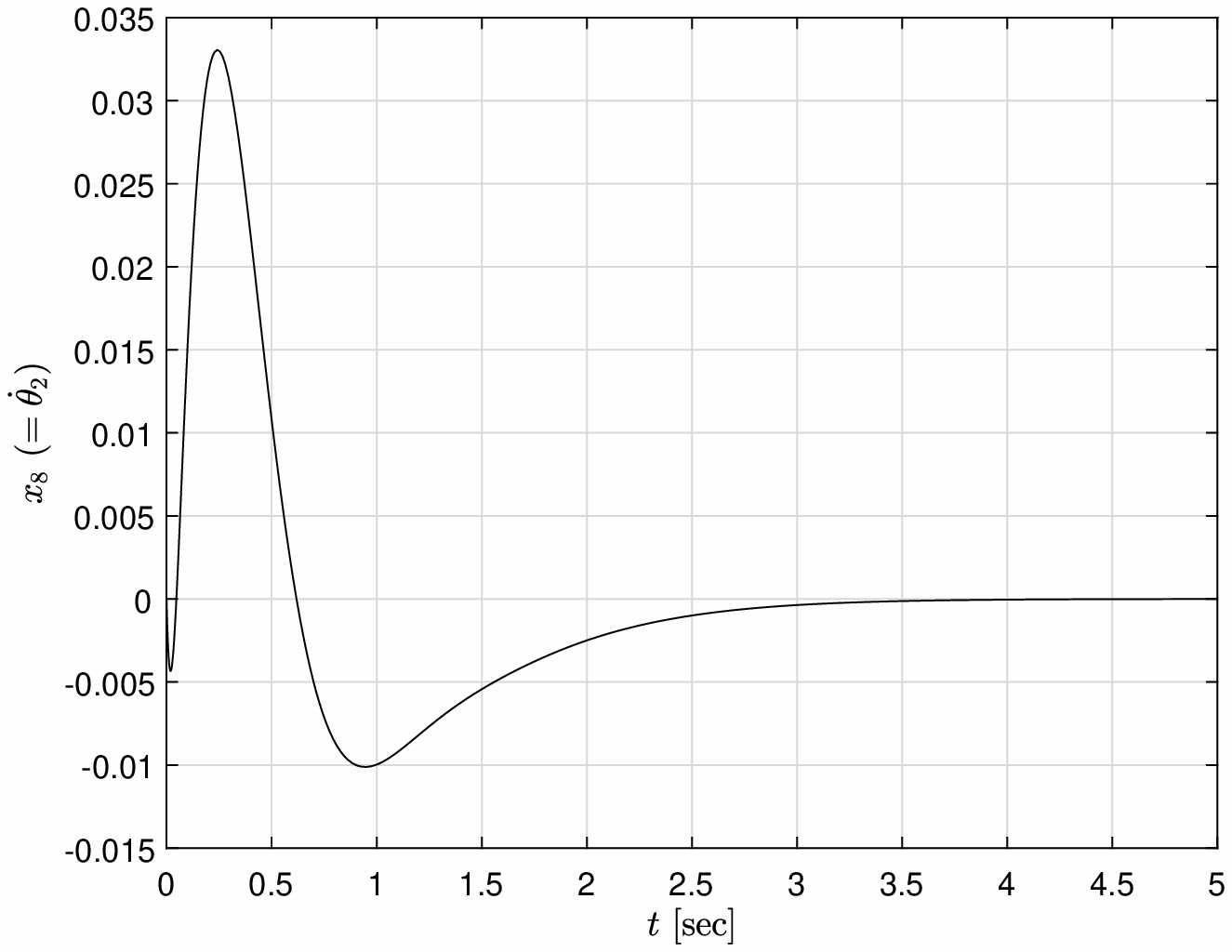,width=5.cm,clip=}
         }
   }
\caption{The solutions $ x_i,$ $i=1,\dots,8$ of the reference model {\bf with considering} the angular torque control $F_{\theta_2}$ with the initial state $(0.234\ 0.2\ 0\ -0.001\ 0\ 0\ 0\ 0)^T.$ Corresponding control forces $F_z,$ $F_{l_1},$  $F_{\theta_1}$ and $F_{\theta_2}$ are depicted in the Fig.\ref{control_forces_with_Ft2}.}
\label{solutionsFt2_xi}
\end{figure}

\begin{figure} 
   \centerline{
    \hbox{
     \psfig{file=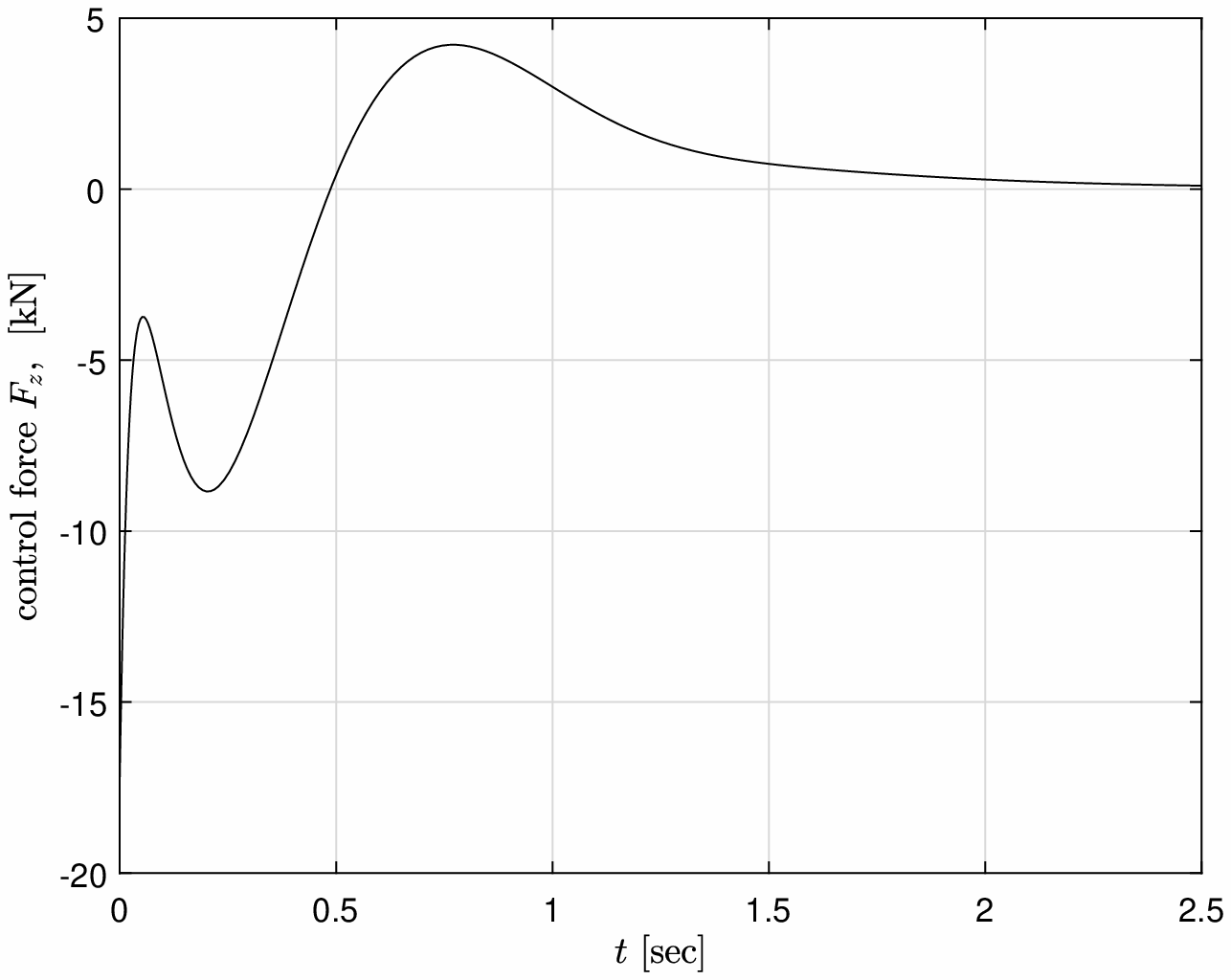,width=5.cm, clip=}
     \hspace{1.cm}
     \psfig{file=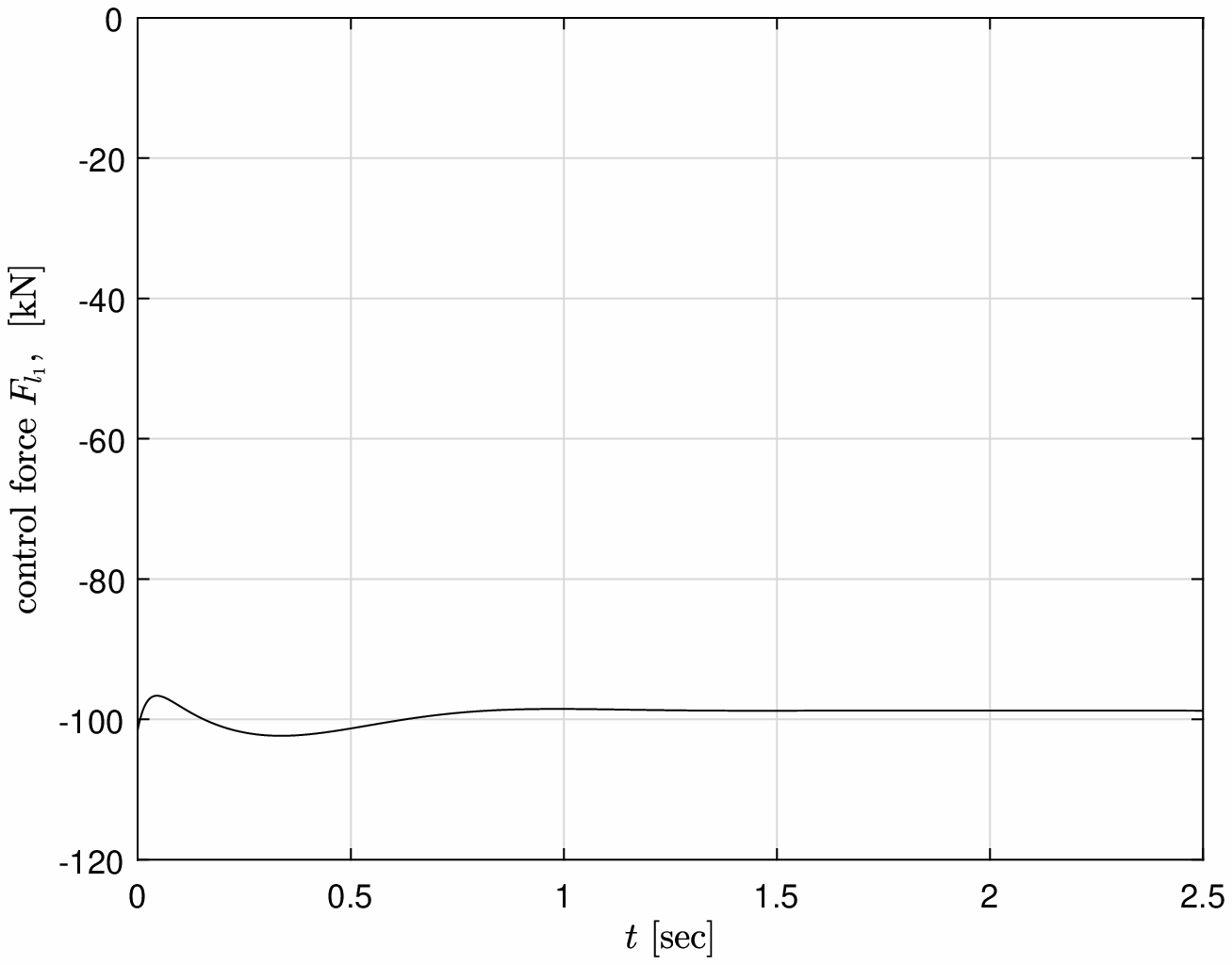,width=5cm,clip=}
    }
   }
   \vspace{0.5cm}
   \centerline{
    \hbox{
     \psfig{file=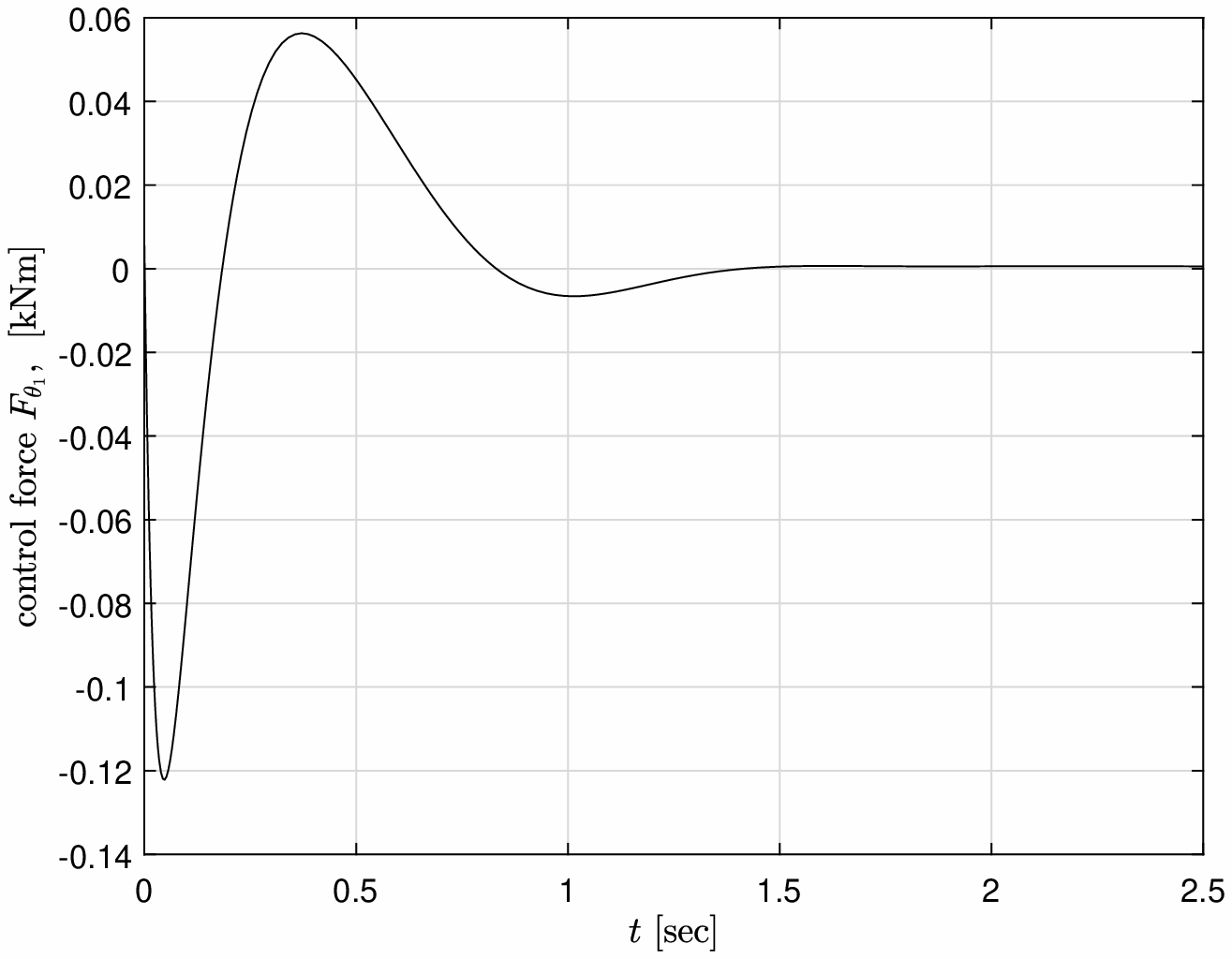,width=5.cm,clip=}
     \hspace{1.cm}
     \psfig{file=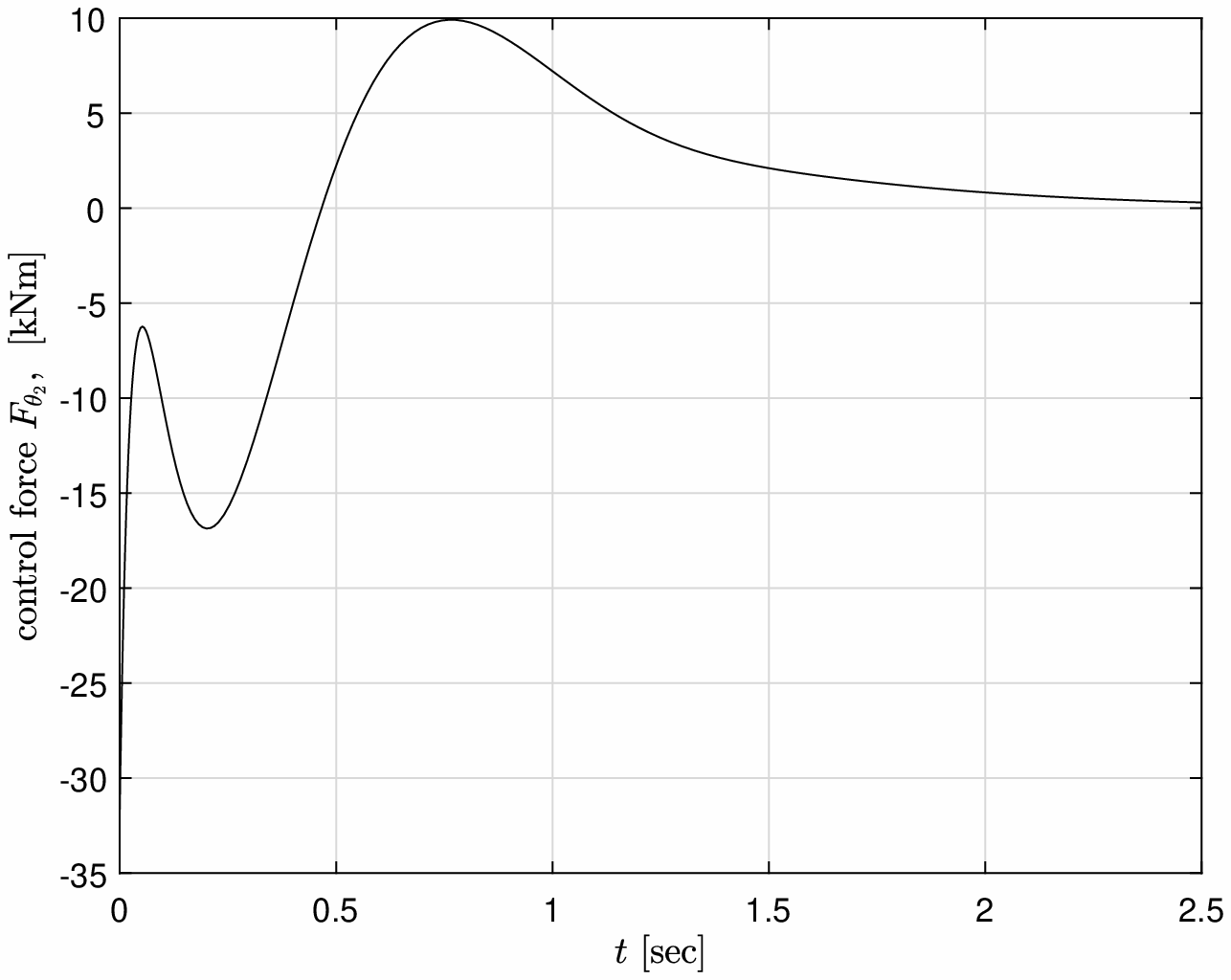,width=5cm,clip=}
         }
   }
\caption{The corresponding control forces to Fig.~\ref{solutionsFt2_xi}; $F_z(t)\rightarrow 0,$ $F_{\theta_1}(t)\rightarrow 0,$ $F_{\theta_2}(t)\rightarrow 0$ and $F_{l_1}(t)\rightarrow -g(m_h+m_p)$ for $t\rightarrow \infty.$}
\label{control_forces_with_Ft2}
\end{figure}

\begin{figure} 
   \centerline{
    \hbox{
     \psfig{file=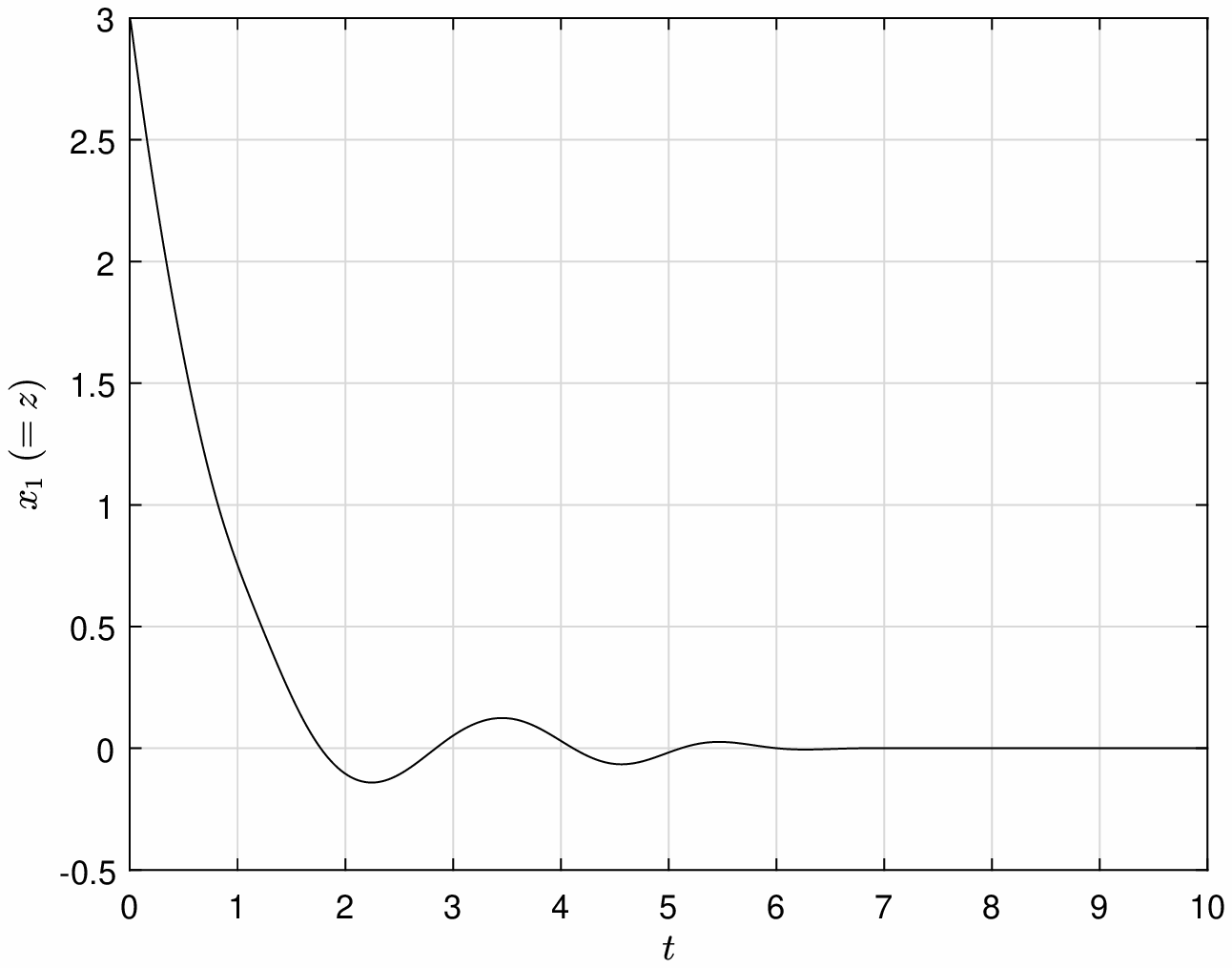,width=5.cm, clip=}
     \hspace{1.cm}
     \psfig{file=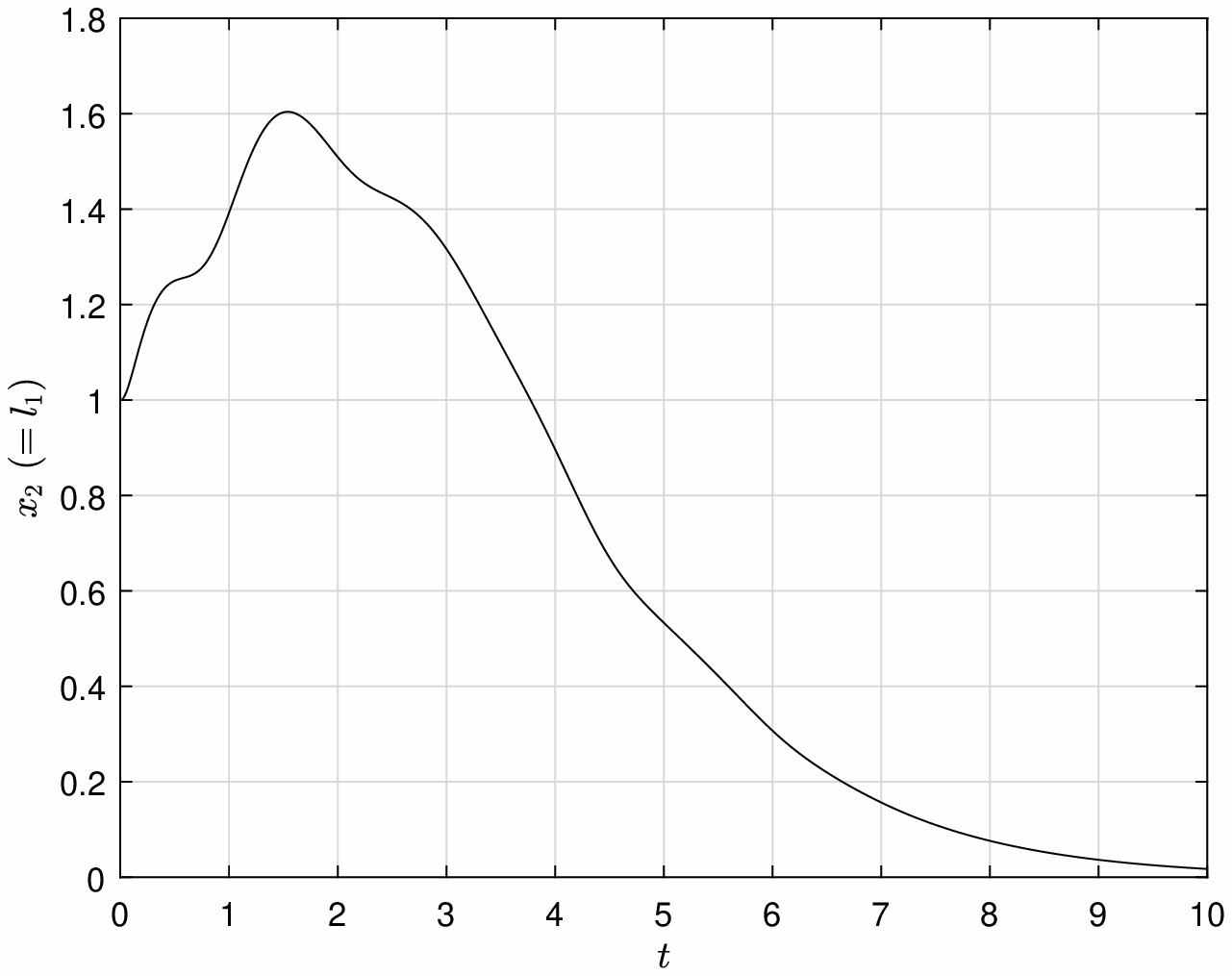,width=5cm,clip=}
    }
   }
   \vspace{0.5cm}
   \centerline{
    \hbox{
     \psfig{file=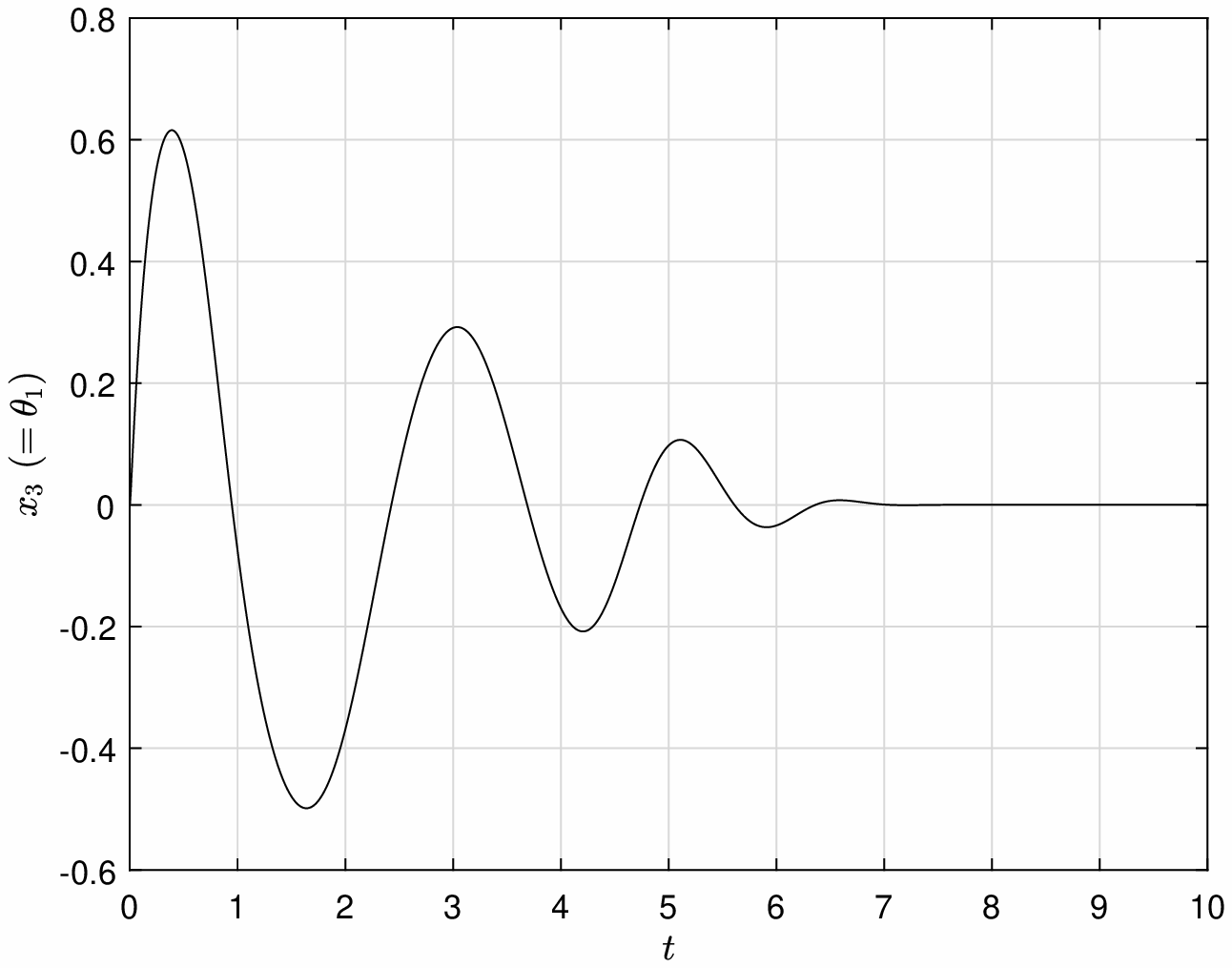,width=5.cm,clip=}
     \hspace{1.cm}
     \psfig{file=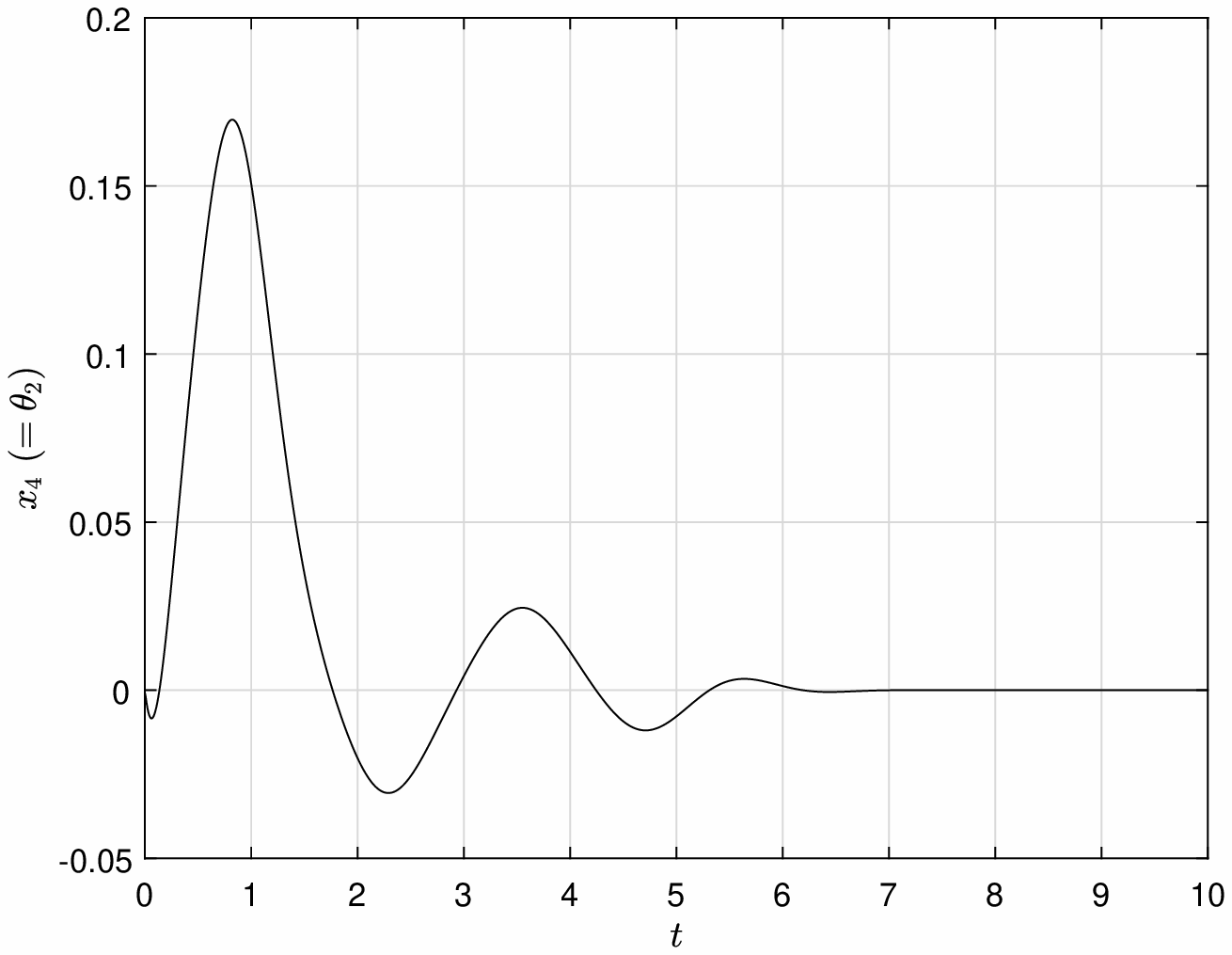,width=5.cm,clip=}
    }
   }
   \vspace{0.5cm}
   \centerline{
    \hbox{
     \psfig{file=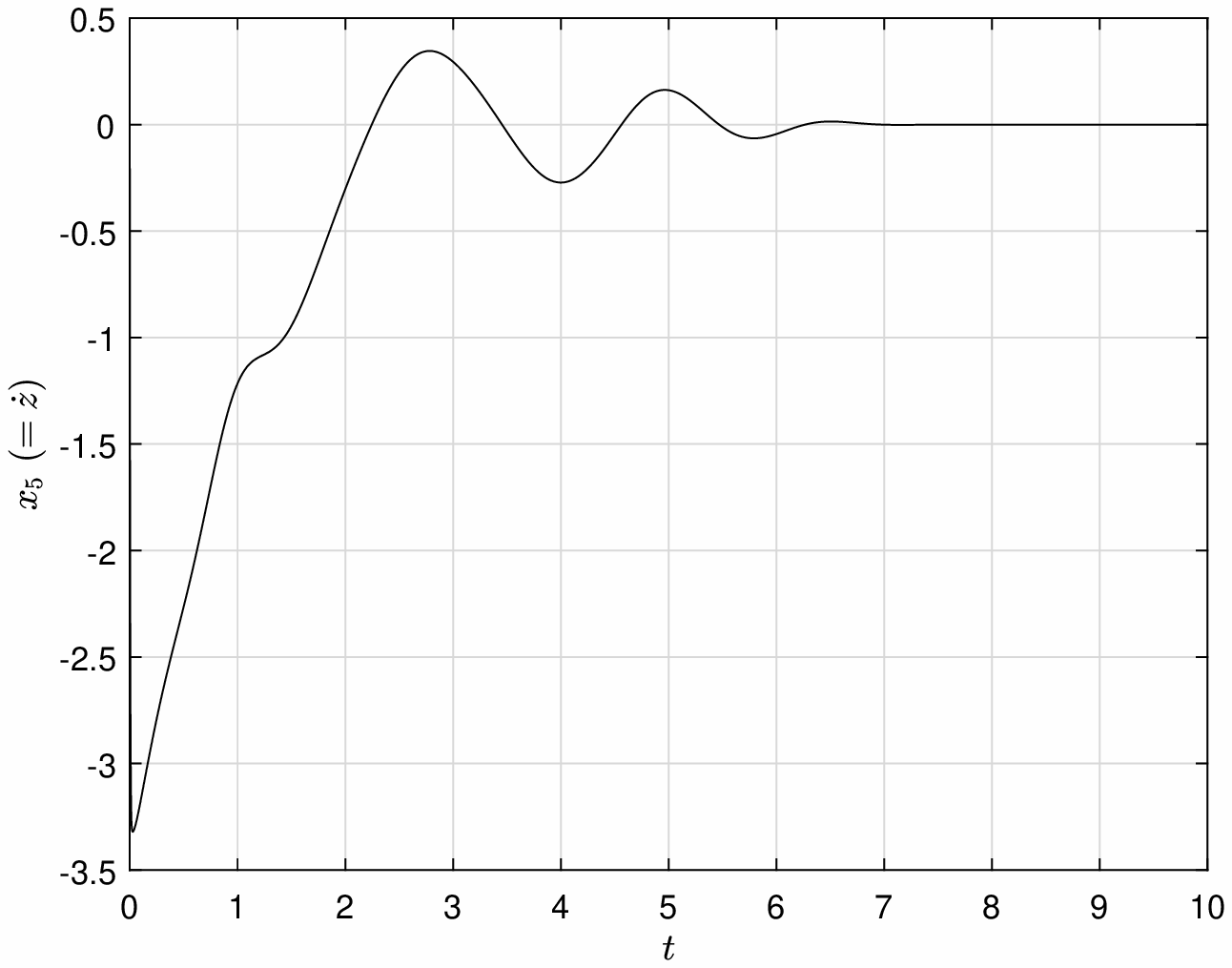,width=5.cm,clip=}
     \hspace{1.cm}
     \psfig{file=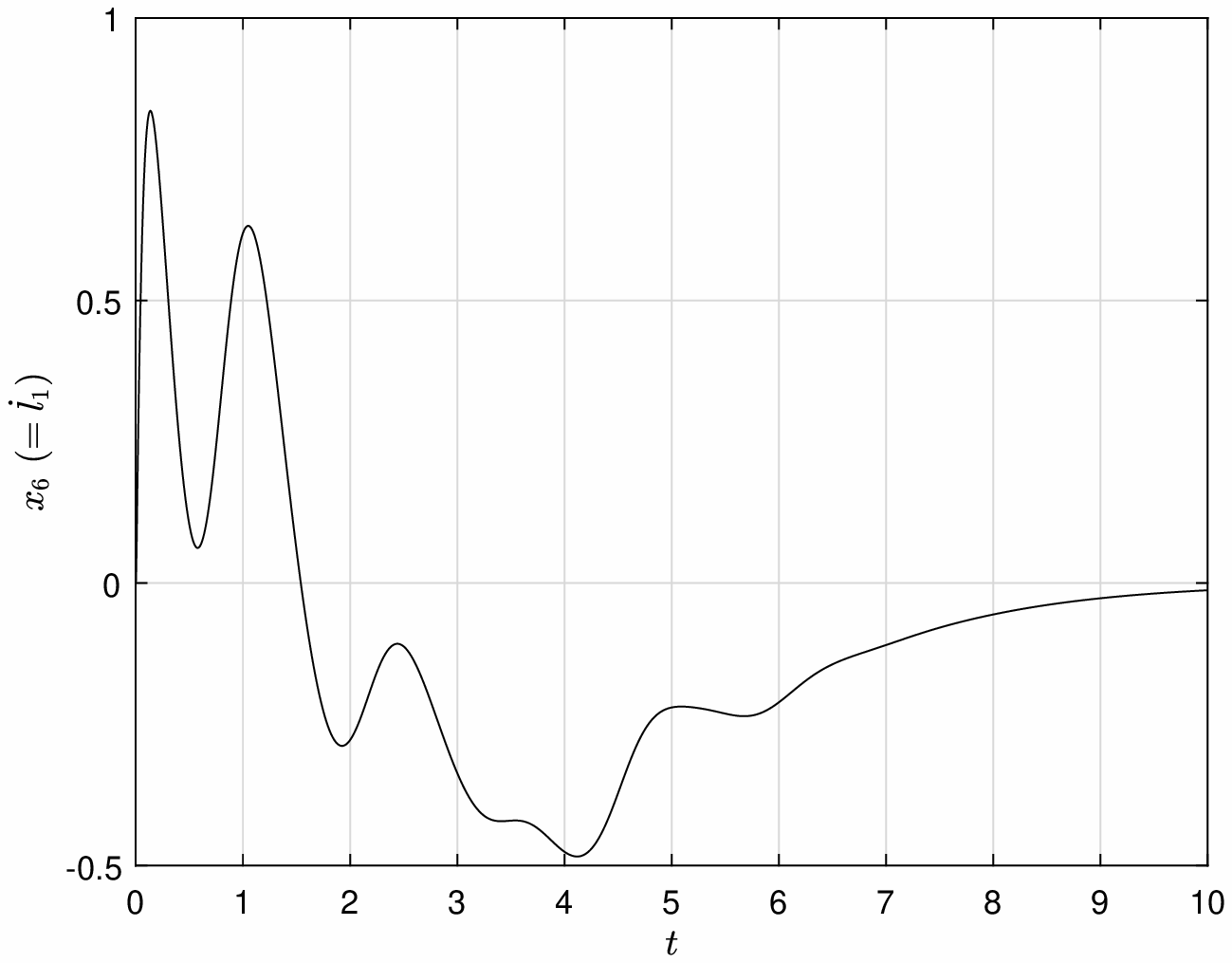,width=5.cm,clip=}
         }
   }
   \vspace{0.5cm}
   \centerline{
    \hbox{
     \psfig{file=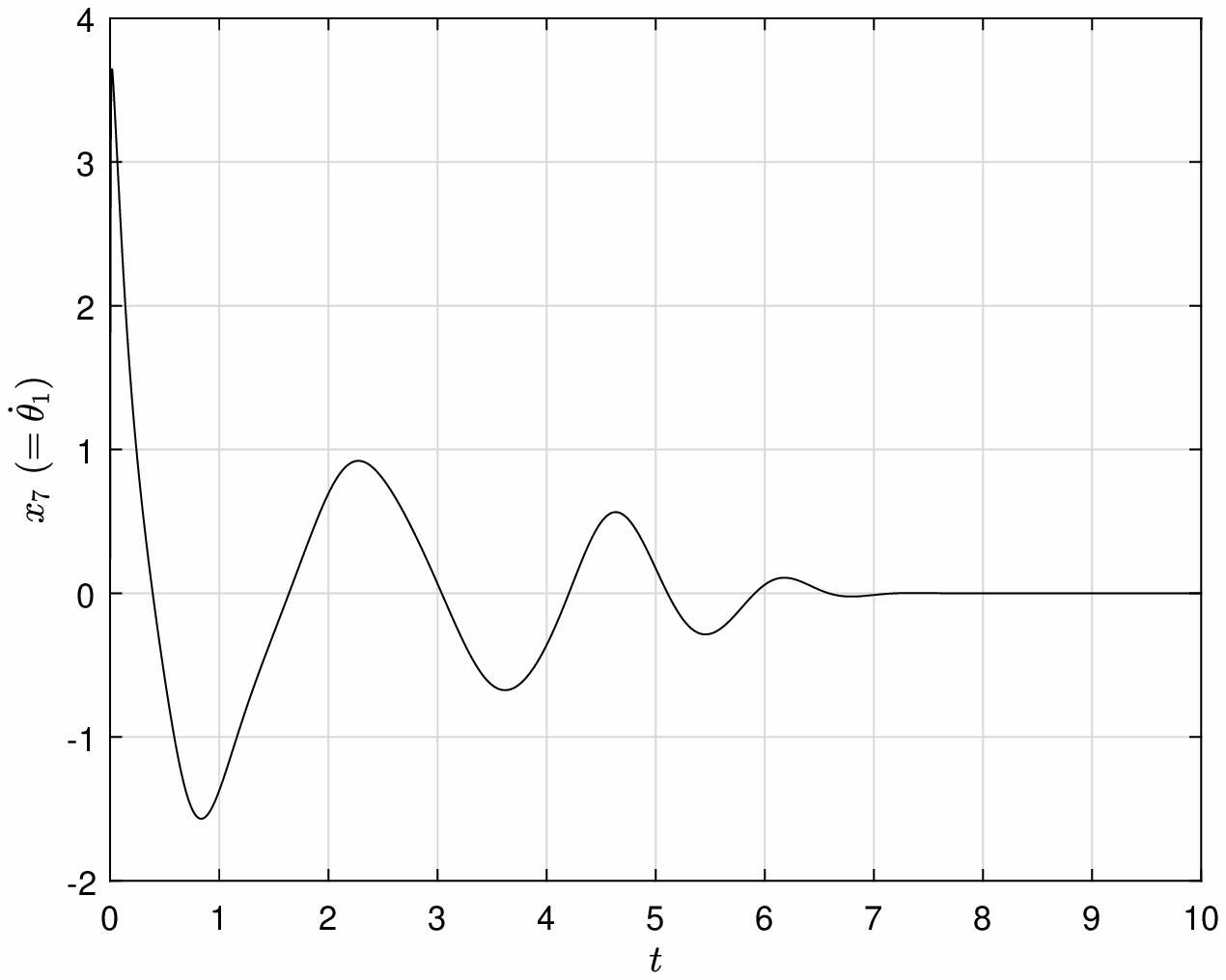,width=5.cm,clip=}
     \hspace{1.cm}
     \psfig{file=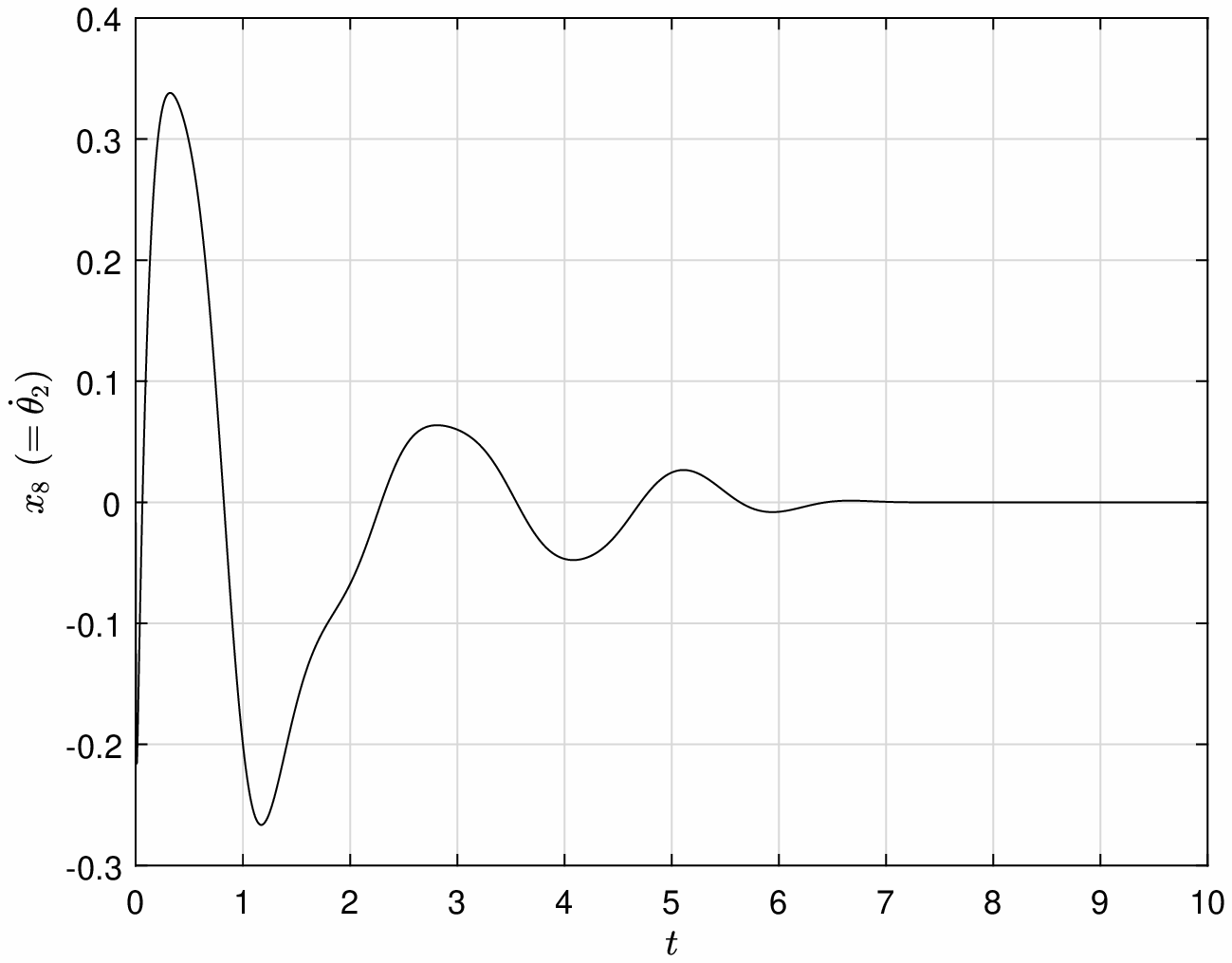,width=5.cm,clip=}
         }
   }
\caption{The solutions $ x_i,$ $i=1,\dots,8$ of the reference model {\bf with considering} the angular torque control $F_{\theta_2}$ with the initial state $(3\ 1\ 0\ -0.001\ 0\ 0\ 0\ 0)^T.$ }
\label{solutionsFt2_3_1_xi}
\end{figure}

Substituting into $\hat G_i,$ $i=5,6,7,8,$ $u_1$ instead of $F_z,$  $u_3$ instead of $F_{\theta_1},$ $u_4$ instead of $F_{\theta_2}$ and (\ref{Fl1}) instead of  $F_{l_1},$  we get the control system $\dot x=\hat G(x,u),$ $\hat G(x,u)=(x_5,x_6,x_7,x_8,\hat G_5,\hat G_6,\hat G_7,\hat G_8)^T,$ $u=(u_1,u_2,u_3,u_4),$ satisfying $\hat G(0,0)=0.$
\begin{rmk}\label{rmk_approx}
\rm The use of the control force $F_{\theta_2},$ represented by the control variable $u_4,$ admits the greater angles of the payload sway during transportation (Fig.~\ref{solutionsFt2_3_1_xi}), therefore  we must work with the complete system, and the often used small-angle approximations of the type $\sin\theta\approx \theta,$  $\cos\theta\approx 1,$ and $\dot\theta^2\approx0$ (\cite{Zhang}, \cite{sun}, \cite{sun2}, \cite{Vaughan}, \cite{Masoud}), can not be applied.
\end{rmk}
Analogously as in the previous section, first we check the controllability of the linear part. Here
\begingroup\makeatletter\def\f@size{8}\check@mathfonts
\[
\hat A=\hat G_x(0,0)=
\left(\begin{array}{cccccccc}
 0 & 0 & 0 & 0 & 1 & 0 & 0 & 0\\
 0 & 0 & 0 & 0 & 0 & 1 & 0 & 0\\
 0 & 0 & 0 & 0 & 0 & 0 & 1 & 0\\
 0 & 0 & 0 & 0 & 0 & 0 & 0 & 1\\
 0 & 0 & 0 & \frac{g{l^2_{2}}{m^2_{p}}}{I_{p}\Sigma_m+{l^2_{2}}m_{h}m_{p}+M{l^2_{2}}m_{p}} & 0 & 0 & 0 & 0\\
 0 & 0 & 0 & 0 & 0 & 0 & 0 & 0\\
 0 & 0 & 0 & 0 & 0 & 0 & 0 & 0\\
 0 & 0 & 0 & -\frac{gl_{2}m_{p}\Sigma_m}{I_{p}\Sigma_m+{l^2_{2}}m_{h}m_{p}+M{l^2_{2}}m_{p}} & 0 & 0 & 0 & 0 \end{array}\right)
\]
and
\[
\hat B=\hat G_u(0,0)=
\left(
\begin{array}{cccc} 
0 & 0 & 0 & 0\\
 0 & 0 & 0 & 0\\
 0 & 0 & 0 & 0\\
 0 & 0 & 0 & 0\\
 \frac{m_{p}{l^2_{2}}+I_{p}}{I_{p}\Sigma_m+{l^2_{2}}m_{h}m_{p}+M{l^2_{2}}m_{p}} & 0 & 0 & -\frac{l_{2}m_{p}}{I_{p}\Sigma_m+{l^2_{2}}m_{h}m_{p}+M{l^2_{2}}m_{p}}\\
 0 & \frac{1}{2I_{h}\left(m_{h}+m_{p}\right)\left(I_{p}\Sigma_m+{l^2_{2}}m_{h}m_{p}+M{l^2_{2}}m_{p}\right)} & 0 & 0\\
 0 & 0 & \frac{1}{I_{h}} & 0\\
 -\frac{l_{2}m_{p}}{I_{p}\Sigma_m+{l^2_{2}}m_{h}m_{p}+M{l^2_{2}}m_{p}} & 0 & 0 & \frac{\Sigma_m}{I_{p}\Sigma_m+{l^2_{2}}m_{h}m_{p}+M{l^2_{2}}m_{p}} 
\end{array}
\right).
\]
\endgroup
The controllability test gives
\[
\rank\left(\hat B,\hat A\hat B,\dots, \hat A^7\hat B\right)=8.
\]
For the same permissible eigenvalues as in previous section,
\[
p = [-0.5 -1 -1.5 -2 -2.5 -3 -3.5 -4],
\]
and the same values of $M,$ $m_p,$ $m_h,$ $I_p,$ $I_h$ and $l_2,$ the gain matrix $K$ is
\[
\hat K=\left(
\begin{array}{rrrrrrrr}
   74.4779  &  0.0040&   0.0051 & 193.1964 &  56.9739 &   0.0035  &  0.0039 & 125.5484\\
   -0.4298  & 66.2825 & -33.7242 &   0.6125 &  -0.1268 & 128.4893 & -26.8565 &   0.1879\\
   -0.0004  & -0.0295 &   0.0758 &   0.0005 &  -0.0001 &  -0.0238 &   0.1304 &   0.0002\\
  136.6796  &  0.0094 &   0.0120 & 233.3502 & 108.1736  &  0.0083 &   0.0092 & 277.5913
\end{array}
\right).
\]
Then the state feedback controller $(u_1,u_2,u_3,u_4)^T=-\hat K(x_1,\dots,x_8)^T,$ which defines the control forces for the reference system. The corresponding solutions $ x_i,$ $i=1,\dots,8$ are depicted on the Fig.~\ref{solutionsFt2_xi} and the control forces $F_z,$ $F_{l_1},$ $F_{\theta_1}$ and $F_{\theta_2}$ on the Fig.~\ref{control_forces_with_Ft2}. The comparison of the Figs.~\ref{solutions_xi_without_Ft2} and~\ref{solutionsFt2_xi} indicates  that the installation of an additional device for the control of the sway motion of payload results to the substantial reduction of the transportation time. The Fig.~\ref{solutionsFt2_3_1_xi} demonstrates the robustness of the designed controller with regard to the initial state data change (error), in comparison with the crane system with uncontrolled payload angle $\theta_2,$ Figs.~\ref{solutions_xi_without_Ft2} and \ref{solutions_xi_without_Ft2_n}.

\clearpage
\newpage

\section{Conclusions}

In this paper, a linear and continuous control law for fully automated double-pendulum-type overhead crane systems with the aim to suppress the sway motion and to reduce the overall time of transportation using the state feedback-based feed-forward control was proposed. We have shown that by appropriate choice of the state feedback gain matrix the crane system can be asymptotically stabilized around the desired end position. Although the technical realization of the additional device for control of the sway angle of payload requires some one-off costs of implementation, the numerical simulation indicates a substantial reduction of the transportation time (up to 80\%) in comparison with the overhead crane system with uncontrolled sway angle of payload, as demonstrate the first sub-figures on top-left in the Figs.~\ref{solutions_xi_without_Ft2} and \ref{solutionsFt2_xi}, and which may be desirable under certain circumstances -- cranes working in demanding and hazardous areas, for example. 


\end{document}